\documentclass[msom]{oo} % current default for manuscript submission
\DoubleSpacedXI

\usepackage{natbib}
 \bibpunct[, ]{(}{)}{,}{a}{}{,}%
 \def\bibfont{\small}%
\TheoremsNumberedThrough     % Preferred (Theorem 1, Lemma 1, Theorem 2)
\EquationsNumberedThrough    % Default: (1), (2), ...
\usepackage{amsmath,amssymb,amsfonts}  
\usepackage[dvipsnames]{xcolor}
\usepackage{tikz} %,tikz-3dplot,graphicx,shadings,epsfig,pstricks,animate,pgfplots
\usetikzlibrary{arrows,shapes,fit} %,positioning,calc,shadows,automata,patterns,decorations.markings,decorations.pathmorphing,plotmarks,snakes,
\usepackage{booktabs,enumitem} %,enumerate,enumitem,longtable,supertabular,lipsum,verbatim
\usepackage{comment} %,,empheq,cancel,filecontents,textpos,accents
\usepackage{url}
\usepackage{eurosym}
\usepackage{graphicx}
\usepackage{multicol}
\usepackage{multirow}
\usepackage{empheq}
\usepackage{float}
\usepackage{subcaption}
\usepackage{rotating}
\usepackage[ruled,vlined, linesnumbered, noend]{algorithm2e}
\usepackage{mathrsfs}
\usetikzlibrary{arrows.meta}
\usetikzlibrary{shapes.geometric}
\usetikzlibrary{fit}
\usepackage[colorlinks=true,linkcolor=blue,citecolor=cyan, urlcolor = black, hypertexnames=false]{hyperref}
\usepackage{lineno}
\usepackage[export]{adjustbox}
\usepackage[flushleft]{threeparttable}

\allowdisplaybreaks

\newtheorem{prop}[theorem]{Proposition}

\newcommand{\cred}{\color{black}}

\newcommand{\Z}{\mathbb{Z_{+}}}
\newcommand{\Kcycle}{K}
\newcommand{\Lchain}{L}
\newcommand{\Kix}{k}
\newcommand{\Lix}{\ell}
\newcommand{\pairsAb}{PDPs}
\newcommand{\pairAb}{PDP}
\newcommand{\pairs}{patient-donor pairs}
\newcommand{\singles}{NDDs}
\newcommand{\single}{NDD}
\newcommand{\pairset}{\mathscr{P}}
\newcommand{\singleset}{\mathscr{N}}
\newcommand{\digraph}{\mathscr{D}}
\newcommand{\Vertex}{\mathscr{V}}
\newcommand{\Arcs}{\mathscr{A}}
\newcommand{\vix}{v}
\newcommand{\vixstar}{v^{*}}
\newcommand{\uix}{u}

\newcommand{\weight}{w}
\newcommand{\cycle}{c}
\newcommand{\chain}{p}

\newcommand{\orderedSet}{\mathcal{S}}
\newcommand{\FirstEorderedSet}{{s}_{1}}

\newcommand{\digraphc}{\digraph= (\Vertex, \Arcs)}

\newcommand{\fm}{i}

\newcommand{\fbset}{\mathscr{V}^{*}}
\newcommand{\fbcopiescy}{\hat{I}}
\newcommand{\fbcopiesch}{\bar{I}}
\newcommand{\Vertexd}{\hat{\mathscr{V}}}
\newcommand{\Arcsd}{\hat{\mathscr{A}}}
\newcommand{\digraphd}{\hat{\mathscr{D}}}
\newcommand{\arcp}{(\uix, \vix)}
\newcommand{\inarcp}{(\vix, \uix)}
\newcommand{\pairarc}{\uix \vix}
\newcommand{\inpairarc}{\vix \uix}

\newcommand{\Vertexbar}{\bar{\mathscr{V}}}
\newcommand{\Arcsbar}{\bar{\mathscr{A}}}
\newcommand{\digraphbar}{\bar{\mathscr{D}}}
\newcommand{\dichainf}{\digraphbar^{\fm}= (\Vertexbar^{\fm}, \Arcsbar^{\fm})}
\newcommand{\oldfvs}{\mathcal{V}^{*}} %Ojo lo cambié
\newcommand{\fvs}{\mathscr{V}^{*}}
\newcommand{\upperc}{n}

\newcommand{\csetf}{\mathscr{\bar{C}}^{\fm}}

\newcommand{\multi}{\lambda}

\newcommand{\Lagdual}{\sigma^{LD}}
\newcommand{\objIP}{\mathcal{Z}}
\newcommand{\objLag}{\mathcal{Z}}
\newcommand{\objLagCy}{\mathcal{\hat{Z}}}
\newcommand{\objLagCh}{\mathcal{\bar{Z}}}
\newcommand{\cweight}{w_{\cycle} }
\newcommand{\pweight}{w_{\chain} }
\newcommand{\enfoarcs}{\Arcs^{*}}
\newcommand{\dualselarcs}{\mu}
\newcommand{\lagx}{\hat{x}}
\newcommand{\lagy}{\bar{x}}

\newcommand{\layer}{\mathcal{L}}
\newcommand{\pos}{\pi}
\newcommand{\dgv}{n}
\newcommand{\cydpvar}{h}

\newcommand{\zcycleb}{z_{\cycle}^{\fm}}
\newcommand{\zchainb}{z_{\chain}^{\fm}}
\newcommand{\cfeaset}{\mathscr{\hat{C}}^{\fm}_{\Kcycle}}
\newcommand{\pfeaset}{\mathscr{\bar{C}}^{\fm}_{\Lchain}}
\newcommand{\chfeasetD}{\mathscr{C}_{\Lchain}}
\newcommand{\cyfeasetD}{\mathscr{C}_{\Kcycle}}

\newcommand{\recucy}{\hat{\eta}^{\fm}}
\newcommand{\recuch}{\bar{\eta}^{\fm}}
\newcommand{\ma}{a}
\newcommand{\inn}{\delta^{\fm}_{-}}

\newcommand{\mgraph}{\mathcal{\hat{G}}^{\fm}}
\newcommand{\mnodes}{\mathcal{\hat{N}}^{\fm}}
\newcommand{\marcs}{\mathcal{\hat{A}}^{\fm}}
\newcommand{\mdd}{\mgraph = (\mnodes, \marcs)}
\newcommand{\basicMDD}{\mathcal{\hat{G}} = (\mathcal{\hat{N}}, \mathcal{\hat{A}})} 
\newcommand{\mgraphch}{\mathcal{\bar{G}}^{\fm}}
\newcommand{\mnodesch}{\mathcal{\bar{N}}^{\fm}}
\newcommand{\marcsch}{\mathcal{\bar{A}}^{\fm}}
\newcommand{\mddch}{\mgraphch = (\mnodesch, \marcsch)}

\newcommand{\lonvar}{y}
\newcommand{\lonArcs}{\mathscr{A}^{'}}
\newcommand{\lonVertex}{\mathscr{V}^{'}}
\newcommand{\lonD}{\mathscr{D}^{'}}
\newcommand{\longraph}{\mathscr{D}^{'} = (\lonVertex, \lonArcs)}
\newcommand{\allCycles}{\mathscr{C}^{'}}
\newcommand{\loncfeaset}{\mathscr{C}^{'}_{\Kcycle}}

\newcommand{\bsub}{\begin{subequations}}
\newcommand{\esub}{\end{subequations}}
\definecolor{beaublue}{rgb}{0.74, 0.83, 0.9}
\definecolor{ballblue}{rgb}{0.13, 0.67, 0.8}
\definecolor{ashgrey}{rgb}{0.7, 0.75, 0.71}

\newcommand{\spacing}{\OneAndAHalfSpacedXI}

\newcommand{\cdionne}{\color{black}}%red!85!black

\newcommand{\BP}{B\&P}

\newcommand{\matching}{M(\Kcycle, \Lchain)}
\newcommand{\matchingPartI}{M(\Kcycle, \Lchain)^{\prime}}

\newcommand{\rev}[1]{\textcolor{black}{#1}}
\begin{document}
%%%%%%%%%%%%%%%%

%\linenumbers
%\tableofcontents

%\setcounter{page}{1}
%\DOI{10.1287/xxxx.0000.0000}%

\RUNAUTHOR{Riascos-\'{A}lvarez et al.}

\RUNTITLE{B\&P algorithm enhanced by MDDs for the KEP}

% \TITLE{Bundling Information Goods of Decreasing Value}
\TITLE{A Branch-and-Price Algorithm Enhanced by Decision Diagrams for the Kidney Exchange Problem}

\ARTICLEAUTHORS{%
\AUTHOR{Lizeth Carolina Riascos-\'{A}lvarez}
\AFF{Department of Mechanical and Industrial Engineering, University of Toronto, Toronto, Ontario M5S 3G8, Canada, \EMAIL{carolina.riascos@mail.utoronto.ca}}
\AUTHOR{Merve Bodur}
\AFF{Department of Mechanical and Industrial Engineering, University of Toronto, Toronto, Ontario M5S 3G8, Canada,
\EMAIL{bodur@mie.utoronto.ca}}
\AUTHOR{Dionne M. Aleman}
\AFF{Department of Mechanical and Industrial Engineering, University of Toronto, Toronto, Ontario M5S 3G8, Canada}
\AFF{Institute of Health Policy, Management and Evaluation, University of Toronto, Toronto, Ontario M5S 3E3, Canada}
\AFF{Techna Institute at University Health Network, Toronto, Ontario M5G 1P5, Canada, \EMAIL{aleman@mie.utoronto.ca}}
}
 % end of the block

\ABSTRACT{%
\textit{\textbf{Problem definition:}} Kidney paired donation programs allow patients registered with an incompatible donor to receive a suitable kidney from another donor, as long as the latter's co-registered patient, if any, also receives a kidney from a different donor. The kidney exchange problem (KEP) aims to find an optimal collection of kidney exchanges taking the form of cycles and chains. \textit{\textbf{Methodology/results:}} % Existing exact solution methods for KEP either are designed for the case where only cyclic exchanges are considered, or can handle long chains but are scalable only if cycles are short.
We develop the first decomposition method that is able to consider long cycles and long chains for projected large realistic instances. Particularly, we propose a branch-and-price framework in which the pricing problems are solved (for the first time in packing problems in a digraph) through multi-valued decision diagrams. We present a new upper bound on the optimal value of the KEP, obtained via our master problem.
Computational experiments show superior performance of our method over the state of the art by optimally solving almost all instances in the PrefLib library for multiple cycle and chain lengths. \textit{\textbf{Managerial implications:}} Our algorithm also allows the prioritization of the solution composition, e.g., chains over cycles or vice versa, and we conclude, similar to previous findings, that chains benefit the overall matching efficiency and highly sensitized patients.
}%

%\KEYWORDS{deterministic inventory theory; infinite linear programming duality; 
%  existence of optimal policies; semi-Markov decision process; cyclic schedule}

\KEYWORDS{Kidney exchange; integer programming; branch and price; multi-valued decision diagrams}
\HISTORY{}

\maketitle
%%%%%%%%%%%%%%%%%%%%%%%%%%%%%%%%%%%%%%%%%%%%%%%%%%%%%%%%%%%%%%%%%%%%%%

% Samples of sectioning (and labeling) in MSOM
% NOTE: (1) \section and \subsection do NOT end with a period
%       (2) \subsubsection and lower need end punctuation
%       (3) capitalization is as shown (title style).
%
%\section{Introduction.}\label{intro} %%1.
%\subsection{Duality and the Classical EOQ Problem.}\label{class-EOQ} %% 1.1.
%\subsection{Outline.}\label{outline1} %% 1.2.
%\subsubsection{Cyclic Schedules for the General Deterministic SMDP.}
%  \label{cyclic-schedules} %% 1.2.1
%\section{Problem Description.}\label{problemdescription} %% 2.

% Text of your paper here

\section{Introduction}
\label{sec:Intro}

%Kidney transplantation accounts for roughly 70$\%$ of annual transplants worldwide. Although dialysis is a common treatment under kidney failure, it is undesirable. The survival rate of patients who undergo dialysis after 10 years is only about 12$\%$, their quality of life is undermined and the treatment carries high costs for healthcare systems. 

The preferred treatment for kidney failure is {\cdionne transplantation, and  demand} for deceased-donor kidneys usually outnumbers supply {\cdionne \citep{CIHI2018}}. An alternative, often desirable, is living-donor transplantation. A living donor is typically a close relative, partner or friend who is willing to donate one of their kidneys to grant a life-saving chance to a beloved one. However, biological incompatibilities, such as blood type or antibodies related discrepancies, between the patient in need of a kidney and the potential donor may exist. It is in these cases where Kidney Paired Donation Programs (KPDPs), present in multiple countries, have played a life-saving role in kidney transplantation systems. A KPDP is a centralized registry operated at a local or national level, where each patient registers voluntarily along with his or her incompatible (suboptimal) donor (\textit{paired donor}) as a pair. Patients in these \pairs\  (\pairsAb) are willing to exchange their paired donors under the promise that they will receive a suitable kidney from a different donor.  To accomplish this goal, two types of exchanges are allowed: cyclic and chain-like exchanges. 

Figure \ref{fig:PoolCycles} illustrates a pool of six \pairs\ $(p_{1}, d_{1})$, $(p_{2}, d_{2})$,...,$(p_{6}, d_{6})$ arranged in two cycles. In the cyclic exchange on the left,  donor $d_{2}$ donates to patient $p_{1}$ and donor $d_{1}$ donates to patient $p_{2}$, thereby  forming a \textit{cycle}. On the right, a cycle involving four \pairsAb\ is depicted, where donor $d_{4}$ donates a kidney to  patient $p_{3}$ and patient $p_{4}$ receives a kidney from donor $d_{6}$, after sequential donations. Due to pair drop-outs, aggravated health condition of a pair member, or last-minute detected incompatibilities, the patient in the first pair may never receive a kidney back if the donor in any subsequent pair in the cycle fails to donate.  To avoid such a risk, cycles of kidney transplants are performed simultaneously in practice, imposing limitations on $\Kcycle$, the maximum size of a cycle, where $\Kcycle \in \Z$. In the literature, a $\Kix$-way cycle refers to a cycle involving $\Kix$ transplants, with $\Kix \le \Kcycle$. Although in some countries such as the United  States, it is common to find $K = 3$ \citep{NKR9years}, other countries have reported longer cycles \citep{CaKEPFoundations, AustraliaKPD},  up to $K = 6$ in Canada \citep{CBS2019}. 

\begin{figure}[tbp]
	\centering
	\vskip 0.2cm
	\begin{adjustbox}{minipage=\linewidth,scale=0.7}
	\begin{subfigure}{.5\textwidth}
		\centering
		% include first image
		
			% Two Colored Circle Split 
			\makeatletter
			\tikzset{circle split part fill/.style  args={#1,#2}{%
					alias=tmp@name, 
					postaction={%
						insert path={
							\pgfextra{% 
								\pgfpointdiff{\pgfpointanchor{\pgf@node@name}{center}}%
								{\pgfpointanchor{\pgf@node@name}{east}}%            
								\pgfmathsetmacro\insiderad{\pgf@x}
								\fill[#1] (\pgf@node@name.base) ([xshift=-\pgflinewidth]\pgf@node@name.east) arc
								(0:180:\insiderad-\pgflinewidth)--cycle;
								\fill[#2] (\pgf@node@name.base) ([xshift=\pgflinewidth]\pgf@node@name.west)  arc
								(180:360:\insiderad-\pgflinewidth)--cycle;            
			}}}}}  
			\makeatother

			\begin{tikzpicture}
			\node[minimum size=5mm, shape=circle split, draw=ashgrey,line width=0.5mm, circle split part fill={blue!20, green!30}] (1) at (-4, -3) {$p_1$\nodepart{lower} $d_1$};
			\node[minimum size=5mm, shape=circle split, draw=ashgrey,line width=0.5mm, circle split part fill={green!30, blue!20}] (2) at (-2, -3) {$d_2$\nodepart{lower} $p_2$};
			\node[minimum size=5mm, shape=circle split, draw=ashgrey,line width=0.5mm, circle split part fill={blue!20, green!30}] (3) at (0, -4) {$p_3$\nodepart{lower} $d_3$};
			\node[minimum size=5mm, shape=circle split, draw=ashgrey,line width=0.5mm, circle split part fill={blue!20, green!30}] (4) at (0, -2) {$p_4$\nodepart{lower} $d_4$};
			\node[minimum size=5mm, shape=circle split, draw=ashgrey,line width=0.5mm, circle split part fill={green!30,blue!20}] (5) at (2, -4) {$d_5$\nodepart{lower} $p_5$};
			\node[minimum size=5mm, shape=circle split, draw=ashgrey,line width=0.5mm, circle split part fill={green!30, blue!20}] (6) at (2, -2) {$d_6$\nodepart{lower} $p_6$};
			
			\draw [-{Classical TikZ Rightarrow [scale=1.4, gray]}, gray, line width = 0.4mm] (1) to [bend right=45] (2);
			\draw [-{Classical TikZ Rightarrow [scale=1.4, gray]}, gray, line width = 0.4mm] (2) to [bend right=45] (1);
			
			\draw [-{Classical TikZ Rightarrow [scale=1.4, gray]}, gray, line width = 0.4mm] (4) to [bend right=45] (3);
			\draw [-{Classical TikZ Rightarrow [scale=1.4, gray]}, gray, line width = 0.4mm] (3) to [bend right=45] (5);
			\draw [-{Classical TikZ Rightarrow [scale=1.4, gray]}, gray, line width = 0.4mm] (5) to [bend right=45] (6);
			\draw [-{Classical TikZ Rightarrow [scale=1.4, gray]}, gray, line width = 0.4mm] (6) to [bend right=45] (4);
			
			\end{tikzpicture}	 
    		\caption{2-way cycle (left) and 4-way cycle (right)}
		\label{fig:PoolCycles}
	\end{subfigure}
	\begin{subfigure}{.5\textwidth}
		\centering
		% include second image
			% Two Colored Circle Split 
			\makeatletter
			\tikzset{circle split part fill/.style  args={#1,#2}{%
					alias=tmp@name, 
					postaction={%
						insert path={
							\pgfextra{% 
								\pgfpointdiff{\pgfpointanchor{\pgf@node@name}{center}}%
								{\pgfpointanchor{\pgf@node@name}{east}}%            
								\pgfmathsetmacro\insiderad{\pgf@x}
								\fill[#1] (\pgf@node@name.base) ([xshift=-\pgflinewidth]\pgf@node@name.east) arc
								(0:180:\insiderad-\pgflinewidth)--cycle;
								\fill[#2] (\pgf@node@name.base) ([xshift=\pgflinewidth]\pgf@node@name.west)  arc
								(180:360:\insiderad-\pgflinewidth)--cycle;            
			}}}}}  
			\makeatother  
			\begin{tikzpicture}
			\node[minimum size=5mm, shape=circle split, draw=ashgrey,line width=0.5mm, circle split part fill={blue!20, green!30}] (3) at (0, -4) {$p_7$\nodepart{lower} $d_7$};
			\node[minimum size=5mm, shape=rectangle, draw=ashgrey,line width=0.5mm, fill=yellow!20, inner sep=12pt] (4) at (0, -2) {$A$};
			\node[minimum size=5mm, shape=circle split, draw=ashgrey,line width=0.5mm, circle split part fill={green!30,blue!20}] (5) at (2, -4) {$d_8$\nodepart{lower} $p_8$};
			\node[minimum size=5mm, shape=circle split, draw=ashgrey,line width=0.5mm, circle split part fill={green!30, blue!20}] (6) at (2, -2) {$d_9$\nodepart{lower} $p_9$};
			\node[minimum size=5mm, shape=circle split, draw=ashgrey,line width=0.5mm, circle split part fill={blue!20, green!30}] (7) at (4, -2) {$p_{10}$\nodepart{lower} $d_{10}$};
			
			\draw [-{Classical TikZ Rightarrow [scale=1.4, gray]}, gray, line width = 0.4mm] (4) to [bend right=45] (3);
			\draw [-{Classical TikZ Rightarrow [scale=1.4, gray]}, gray, line width = 0.4mm] (3) to [bend right=45] (5);
			\draw [-{Classical TikZ Rightarrow [scale=1.4, gray]}, gray, line width = 0.4mm] (5) to [bend right=45] (6);
			\draw [-{Classical TikZ Rightarrow [scale=1.4, gray]}, gray, line width = 0.4mm] (6) to [bend left=45] (7);
			\end{tikzpicture}
		\caption{4-length chain}
		\label{fig:PoolChains}
	\end{subfigure}
	\end{adjustbox}
	\caption{Examples of exchanges. \pairsAb\ are represented by circle nodes, while the \single\ is square.}
	\label{fig:Exchanges}
\end{figure}

In the presence of singleton donors, chains become an exchange alternative. Some singleton donors are \textit{altruistic donors} because they are willing to donate  one of their kidneys without having a paired patient; such donors are represented as vertices with outgoing arcs to any feasible patient-donor pair vertex, and no incoming arcs. A \textit{chain} is a path that starts with a singleton donor donating to a patient in a \pairAb, after which all remaining patients receive a kidney from a paired donor (Figure \ref{fig:PoolChains}). The donor in the last pair of a chain can either donate to a patient on the transplant waitlist or become a \textit{bridge donor} (a single donor in a future chain). In general, chain-initiating donors such as altruistic, bridge, or even deceased {\cdionne donors} \citep{DDonors2017}, are referred to as \textit{non-directed donors} (\singles). Unlike cycles, {\cdionne no paired donor in a chain risks not receiving a kidney for their intended recipient, making it possible to relax the simultaneity constraint}. Since simultaneous surgeries (e.g., six surgeries for a three-way cycle) may bring capacity and logistics challenges, a simultaneity requirement poses stronger restrictions on the maximum cycle length than on the maximum  chain length \citep{Dickerson2012, Biro2021}. To the best of our knowledge, chains are unbounded in the nationwide KPDP in the United States \citep{Dickerson2012}. These very long chains can be carried out in practice through the use of bridge donors that connect multiple segments of the same chain. In most European countries, both cycles and chains are
carried out simultaneously \citep{Biro2021}, which results in shorter algorithmic chains compared to the United States, but still considered long (up to six transplants) in some European countries, e.g., Czech Republic and Italy \citep{Biro2021}. The Canadian KPDP allows long chains with up to six transplants that are not necessarily performed simultaneously \citep{CBS2019}. Thus, very often, the maximum chain length $\Lchain \in  \Z$ is as big or bigger than $\Kcycle$ ($\Lchain \ge \Kcycle$).
%\citep{CaKEPFoundations, NKR9years}

As some transplants may be more suitable or urgent from medical and logistical points of view, a score may be given to every potential transplant. A common objective in kidney exchange, although not the only one, is to match donors (paired and singleton) to patients such that the sum of the transplants' score is maximized \citep{Abraham2007, Roth2007}. The matching must satisfy that every \pairAb\ belongs to at most one cycle or chain and every \single\ to at most one chain. Finding a matching of maximum score in this context is known as the \textit{kidney exchange problem} (KEP).

Depending on the country's population and matching frequency, the average instance size may vary. For instance, \cite{Dickerson2016} used real match runs from the UNOS program which have on average $<$250 \pairsAb\ and a low number of altruists, which may be explained by the semi-weekly matching frequency. In Canada, matching frequency is every four months, with an average instance size around 150 \pairsAb\ and $\le$12 \singles\ with maximum cycles and chain lengths of 6 \citep{CBS2019}. The literature has frequently tested matching algorithms not only with respect to the current size of instances, but with respect to expected instance size increases in the future, e.g., 500 \citep{Glorie2014}, 700 \citep{Dickerson2016}, and 10,000 \pairsAb\ \citep{Abraham2007}. Additionally, problem sizes increase with collaborations among multiple KPDPs sharing participant pools with either other co-national programs or within organized large-scale international exchanges \citep{Ashlagi2021, Klimentova2021, Rees2017}, e.g., a recent collaborative pilot program between the Italian National Transplant Center and the Alliance for Paired Kidney Donation from the United States \citep{Italy2022}. These collaborations bring new operational and algorithmic challenges due to possibly larger exchange pools that may benefit from longer cycles and chains. 

Therefore, motivated by practical settings for which $K \ge 4$ along with long chains \citep{CaKEPFoundations, Ashlagi2021}, previous studies showing the value of long chains \citep{Ashlagi2012, Dickerson2012, Ding2018}, the increasing number of participants in KPDPs, and space for improvement in state-of-the-art  approaches, we develop a new solution technique for the KEP. Particularly, we make the following contributions:

\begin{enumerate}[leftmargin=0.5cm]
	\item We devise a \BP\ algorithm that incorporates cycles and long chains and solve the pricing problems via multi-valued decision diagrams (MDDs). To the best of our knowledge, this is the first study using MDDs in cycle and path packing problems in a digraph, and one of the two works in the \BP\ literature \citep{Arvind2018}. Our \BP\
    \begin{itemize}
        \item Is the first \emph{correct}\footnote{The first B\&P to ever address the cycles-and-chains variant was published in M\&SOM \citep{Glorie2014}, but was later shown to be incorrect \citep{Plaut2016b}.} B\&P algorithm for the cycles-and-chains variant \citep{Omer2022}.
        \item Provides an exact optimal solution.
        \item Can scale to instances sizes not previously possible.
        \item Allows the prioritization of cycles over chains and vice-versa, properties that are desirable in practice even when the instances are small enough for other algorithms to solve to optimality.
    \end{itemize}

	\item We present an effective three-phase solution method for the pricing problems, shifting between MDDs and linear (worst-case integer) programs.
	%\item We derive an upper bound on the optimal value of the KEP, that can be used when optimality of the column generation algorithm at the root node cannot be proven within the time limit.
    \item \rev{We present the first use of Lagrangian relaxation in this context. Additionally, we show for the first time that the dual of our Lagrangian relaxation formulation is equivalent to the disaggregated cycle formulation \citep{Klimentova2014}. This observation allowed us to propose the only other known upper bound on the optimal value of the KEP, which is tighter than the previously proposed upper bound \citep{Abraham2007}, while recent works provided no upper bound for instances that timed out \citep{Lam2020}.}
	\item We demonstrate computational improvement over state-of-the-art methods on publicly available instances.
	\item Given that multiple optima can exist in KEP \citep{Klimentova2021}, and that some of them may be preferred for logistical or fairness reasons, our algorithm allows the prioritization of chains or cycles in a single run.  Additionally, we perform an experimental analysis showing the impact of the chain/cycle composition of multiple solutions.
	\item We  empirically demonstrate the benefits of chains for different graph densities and sizes, taking into account the presence of highly sensitized patients. Our results are in agreement with previous studies \citep{Ashlagi2012, Dickerson2012, Ding2018}.
\end{enumerate}

\rev{Additionally, we note that while B\&P implementations are always bespoke frameworks for the particular application at hand, our framework could generally be used to detect large cycles and/or chains in networks, which has applications in prize-collecting traveling salesman problems and in longest-path problems on non-directed acyclic graphs. Moreover, the number of the MDDs for cycles provides a feasible solution to the feedback vertex set problem. Additionally, our framework could be used in applications for sequencing problems that can be formulated as MDDs, e.g., multicommodity pick-up-and delivery TSP \citep{Castro2020}.}

The rest of the paper is organized as follows.  In Section  \ref{sec:LitReview}, we review the relevant literature. In Section \ref{sec:PblmD}, we present a formal definition of the KEP. In Section \ref{sec:BP}, we detail our \BP \ algorithm, including the reformulation of pricing problems via MDDs. In Section \ref{sec:SolAp}, we present our general solution approach. In Section \ref{sec:Results}, we show experimental results comparing our algorithm with  the state of the art. Lastly, we draw some conclusions and point to future work in Section \ref{sec:Conclusions}.

\section{Literature Review}
\label{sec:LitReview}

A very-well studied variant of the KEP is the cycle-only version, i.e., a problem instance in which either there are no \singles\ or if present, chains are ``turned'' into cycles by adding an arc from every \pairAb\ to \singles. {\cdionne As a result, the maximum length is the same for both cycles and chains. }%, which can represent a donation from the donor in the last \pairAb\ to a patient on the waitlist. 
\rev{Different methods, mostly mixed integer programming (MIP) but also heuristics (see, e.g., the cycle-only approach by \citet{delorme2022improved}), have been used to model this variant of the KEP in the literature. While heuristics can be fast and are often worth a trade-off in optimality for speed, in the specific example of KEP, being able to match even just one more patient per matching can be life-changing for that one person. In countries where matchings are infrequent (e.g., matchings are only done every four months in Canada), having to wait just one more round of matching to receive a kidney could be very damaging to the patient's deteriorating health and quality of life. Thus, from both a practical and a philosophical point of view, if we can tractably solve large KEP instances to optimality, there is benefit to choosing an optimal approach over a heuristic, and therefore we focus on optimal MIP approaches in the literature.}

\citet{Abraham2007} and \citet{Roth2007} proposed two widely known MIP formulations: the \textit{cycle formulation}, which has an exponential number of decision variables, and the \textit{edge formulation}, which has an exponential number of constraints. \cite{Constantino2013} showed that the edge formulation scales  substantially worse than the cycle formulation, reaching more than three million constraints in instances with only 50 \pairsAb, while the cycle formulation leads to memory issues when $\Kcycle \ge 4$ in medium and high density instances with 100 \pairsAb\ or more. %to either enumerate exhaustively feasible cycles and introduce a cycle variable for each, or to add path-violating constraints in an arc-based model.%either enumerate exhaustively feasible cycles and introduce a cycle variable for each, or to add path-violating constraints in a model with arc variables on a digraph representing possible transplants between \pairsAb. These

\cite{Constantino2013}  proposed the first two MIP formulations where the number of constraints and variables are polynomial in the size of the input, referred to as \textit{compact formulations}. It was shown that their \textit{extended edge formulation} outperforms their \textit{assignment edge formulation}. Although the cycle formulation is theoretically stronger than both, the extended edge formulation is able to scale in instances where the cycle formulation requires more than three million variables. 

{\cdionne \BP\ \citep{Barnhart1998} is commonly used to overcome the exponential number of variables of the cycle formulation}, yielding the most successful solution methods for the KEP to date \citep{Abraham2007, Klimentova2014, Dickerson2016, Dickerson2019}. {\cdionne For cycle-only KEP, the state of the art is the \BP-and-cut developed by \cite{Lam2020}, where the cycle formulation is used as a master problem, and few cuts were added to the master problem in most runs on problem instances from the PrefLib library \citep{Mattei2013}}. {\cdionne Most instances with 2048 \pairsAb\ {\cdionne and} $K = 3$ reported total run-time of 2s and for the majority of instances with up to $1024$ \pairsAb, $< 1$s. For $K = 4$ only instances with up to $1024$ \pairsAb\ could be solved, taking up to 22min.}%Overall, a low number of cuts was added to the master problem in most runs.

In the general version of the KEP, chains are allowed, and  usually $\Lchain \ge \Kcycle$ or unbounded. \citet{Anderson2015} proposed two formulations for unbounded chains: recursion and a prize-collecting traveling salesman problem (PC-TSP) variation. Instances with up to 1179 \pairsAb\ and 62 \singles\ were tested with $K = 3$. The recursive formulation outperformed the PC-TSP formulation on a large historical dataset, although, in ``difficult'' instances, the PC-TSP formulation was more successful. This formulation can be modified to include bounds on the length of chains. However, \citet{Plaut2016b} showed that the PC-TSP is effective when unbounded chains and $K = 3$ are considered, otherwise, \BP-based algorithms outperform it. 
%{\cdionne In the general version of the KEP, chains are allowed, and  usually $\Lchain \ge \Kcycle$ or unbounded.} \cite{MakHau2017} introduced a compact formulation (denoted EE-MTZ) integrating chains and cycles through a variation of the well-known Miller-Tucker-Zemlin constraints and ideas similar to the extended edge formulation. {\cdionne An exponential-sized variant of the EE-MTZ including path-violating constraints, SPLIT-MTZ,  was also proposed. The largest instance presented includes 250 \pairsAb\ and 6 \singles with graph density of 5\%. EE-MTZ outperforms SPLIT-MTZ, although both formulations do not scale successfully.}

\citet{MakHau2015} introduced two formulations (EE-MTZ and SPLIT-MTZ), inspired by \cite{Abraham2007, Constantino2013} and the well-known  Miller-Tucker-Zemlin constraints. The largest instance presented includes 250 \pairsAb\ and 6 \singles\ with arc density of 5\%. Overall, EE-MTZ optimally solves more instances than SPLIT-MTZ.
%n the first formulation, infeasible cycles resulting from candidate integer solutions are ruled out on the fly, whereas in the second formulation subtours are removed  in a branch-and-cut algorithm.

\citet{Dickerson2016} proposed three formulations: a compact formulation for the cycle-only version, PIEF, a compact variation of PIEF allowing long chains, HPIEF, and an exponential-sized formulation also allowing long chains, PICEF. The first two formulations are adapted from the extended edge formulation and the last one is inspired by both the extended edge formulation and the cycle formulation. Instances from real match runs and from a realistic simulator were used to test seven algorithms \citep{Abraham2007,Klimentova2014, Anderson2015, Glorie2014, Plaut2016, Dickerson2016}, with $\Kcycle = 3$ and a time limit of 1h. On real match runs from the United Kingdom, the algorithm by \citet{Klimentova2014} outperformed the others when no \singles\ were included. However, for larger (simulated) instances with up to 700 \pairsAb\ and 171 \singles, HPIEF and PICEF were the most effective solution methods. \rev{A direct search approach for the pricing problem \citep{Abraham2007} enumerates all cycles in the worst case, which is exponential when also considering chains (the search for chains is $\mathcal{O}(\lvert N \rvert \lvert P \rvert^{L})$ and the search for cycles is $\mathcal{O}(\lvert P \rvert^{K})$). Unlike the dynamic programming algorithms in \citet{Glorie2014} and \citet{Plaut2016}, which were shown incorrect for the cycles-and-chains variant \citep{Plaut2016b}, our solution to the pricing problems occurs in acyclic graphs (our MDDs), which are based on our cycle and chain copies of the input graph. By construction of the copies, there is no positive cycle or chain in the compatibility graph if none is found in the corresponding copy. Thus, our recursive function simply becomes the computation of the longest path in an acyclic network, which can be performed efficiently.}

More recently, \citet{Duncan2019} proposed the position-indexed traveling-salesman problem formulation (PI-TSP) as part of a robust optimization model in which cycles are handled as in the cycle formulation, and chains are expressed by combining the indexing scheme presented in PICEF and ideas from the PC-TSP formulation. Experimentally, robust solutions were compared to deterministic solutions obtained by PICEF. %A significant reduction in the number of variables to model long chains is achieved. 

Other variations of the KEP include finding robust solutions \citep{Dickerson2014R, Duncan2019, Carvalho2020}, solutions that maximize the expected number of transplants \citep{Klimentova2016, Alvelos2019} and a matching that optimizes the best outcome over time \citep{Awasti2009, Monteiro2020}. Most of these approaches determine a matching that maximizes the number of transplants without distinguishing among several optima. However, real-world applications often require multiple criteria to break ties \citep{Glorie2014, Biro2019}. Tie-breaking criteria has been studied as an ordered list of objectives that are solved sequentially by adding rationality constraints \citep{Glorie2014}. While our main goal remains to maximize the number of transplants, we consider the prioritization of an exchange type (cycles or chains) as a second objective. We achieve this goal by strategically including columns in our column generation algorithm, allowing us to optimize both objectives in a single run.

\section{Problem Description}
\label{sec:PblmD}

Given a set of \pairsAb\ $\pairset$, a set of \singles\ $\singleset$, and positive integers $\Kcycle$ and $\Lchain$, the KEP can be defined in a digraph $\digraphc$. A vertex $\vix \in \Vertex$ is defined for each \pairAb\ and \single, i.e., $\Vertex = \pairset \cup \singleset$. For $\vix_i,\vix_j \in \Vertex$, there exists an arc $(\vix_i,\vix_j) \in \Arcs$ if the donor in vertex $\vix_i$ is compatible with the patient in vertex $\vix_j$, e.g., see Figure \ref{fig:Exchanges}, \rev{which illustrates possible chain and cycle solutions}. Note that, $\Arcs \subseteq \{(\vix_i,\vix_j) \mid \vix_i \in \Vertex , \vix_j \in \pairset \}$, that is, there are no incoming arcs to \singles, neither from \pairsAb\ nor from the other \singles. Each arc $(\vix_{i}, \vix_{j})$ is assigned a score $\weight_{ij}$ representing the suitability/priority of that transplant, \rev{where suitability refers to the quality of the match from a medical perspective}. \rev{With these weights, priorities for transplants can represent varying priority rules, e.g., preference to match highly sensitized patients, patients that have been on the waitlist longest, etc.}

Chains and cycles correspond to simple paths and cycles of $\digraph$, respectively, formally defined as follows:

\begin{definition}
	{\cdionne A chain $\chain = (\vix_1, ...,\vix_\Lix)$ is feasible if and only if: (1)		 $(\vix_i, \vix_{i+1}) \in \Arcs$  $ \forall i = 1,...,\Lix -1$, (2) 
	 $\vix_1 \in \singleset$ and $\vix_i \in \pairset$ $ \forall i = 2,...,\Lix$, and 
		(3) $\Lix \le \Lchain + 1$. Note that for a chain to include $\Lchain$ transplants, $\Lchain + 1$ vertices are required.}
	%\begin{itemize}
	%	\item $(\vix_i, \vix_{i+1}) \in \Arcs  \text{, }  \forall 1 \le i \le \Lix -1$
	%	\item $\vix_1 \in \singleset$ and $\vix_i \in \pairset \text{, } \forall 2 \le i \le \Lix$
	%	\item $\Lix \le \Lchain$
	%\end{itemize}	
\end{definition}

\begin{definition}
{\cdionne A cycle $\cycle= (\vix_1, ..., \vix_\Kix, \vix_1)$ is feasible if and only if: (1) $(\vix_i, \vix_{i+1}) \in \Arcs$  $ \forall i =1,...,\Kix -1$ and $(\vix_{\Kix}, \vix_{1}) \in \Arcs$, (2) $\vix_i \in \pairset$ $ \forall i = 1,...,\Kix$, and (3)  $\Kix \le \Kcycle$.}
%\begin{itemize}
%	\item $(\vix_i, \vix_{i+1}) \in \Arcs  \text{, } \forall 1 \le i \le \Kix -1$ and $(\vix_{\Kix}, \vix_{1}) \in \Arcs$
%	\item $\vix_i \in \pairset \text{, }  \forall 1 \le i \le \Kix$
%	\item $\Kix \le \Kcycle$
%\end{itemize}	
\end{definition}

 A feasible solution to the KEP is a matching of donors to patients. To formally introduce this definition, consider $\cyfeasetD$ and $\chfeasetD$ as the set of all feasible cycles and chains in $\digraph$, respectively. Throughout the paper, notation $\Vertex(\cdot)$ and $\Arcs(\cdot)$ will denote the set of vertices and arcs present in the argument $(\cdot)$, respectively. For instance, $\Vertex(\cycle)$ corresponds to the set of vertices in cycle $\cycle$.

\begin{definition}
	$M(\Kcycle, \Lchain) \subseteq \cyfeasetD \cup \chfeasetD$ is a feasible matching of donors to patients if $\Vertex(m_1) \cap \Vertex(m_2) = \emptyset$ for all $m_1, m_2 \in M(\Kcycle, \Lchain)$ such that $m_1 \ne m_2$.
\end{definition}

%\begin{definition}
%	$M(\Kcycle, \Lchain) = \{\cycle_{1}, \cycle_{2} \text{,..,} \cycle_{n_c} \text{,} \chain_{1}, \chain_{2} \text{,...,}\chain_{n_p - 1}, \chain_{n_p}\}$ is a feasible matching of donors to patients if and only if  $\Vertex(\cycle_{1}) \cap \Vertex(\cycle_{2}) \text{,..,} \Vertex(\chain_{n_p -1}) \cap \Vertex(\chain_{n_p}) = \emptyset$, where $\Vertex(\cycle_i) \subseteq \Vertex$ and $\Vertex(\chain_i) \subseteq \Vertex$ denote the set of vertices in cycle $\cycle_i \in \loncfeaset \text{, } \forall 1 \le i \le n_c$  and  chain $\chain_i \in \lonchfeaset \text{, } \forall 1 \le i \le n_p$, respectively.
%\end{definition}

That is, a matching in the KEP corresponds to a collection of feasible chains and cycles where every patient/donor participates in at most one transplant and type of exchange. Thus, the objective of the KEP is to find a matching, whose set of arcs $\Arcs(M(\Kcycle, \Lchain))$ maximizes the total transplant score, i.e., $\sum_{(i,j) \in \Arcs(M(\Kcycle, \Lchain))} \weight_{ij}$.

\section{{\cdionne \BP\ Algorithm}}
\label{sec:BP}

{\cdionne Our approach to KEP} is, to the best of our knowledge, the only \BP \ for KEP valid for long{\cdionne, yet bounded chains}. For completeness, we first discuss the general motivation behind \BP.  {\cdionne We then detail our \BP \ implementation.} Particularly, we focus on solving the pricing problems via multi-valued decision diagrams, a solution method novel to cycle and path packing problems in digraphs, and only used before in a transportation scheduling problem \citep{Arvind2018},  to the best of our knowledge.

 \subsection{Background}
 
%For large instances, the cardinality of $\cfeaset$ and $\pfeaset$ becomes prohibitive up to the point that we cannot exhaustively state all decision variables in \eqref{objDisIDCF} or constraints in \eqref{objLag2opt}. 
Instead of considering the full set of variables in a linear program, \textit{column generation}  works with a small subset of variables (\textit{columns}), forming the well-known \textit{restricted master problem} (RMP). If by duality theory a column is missing in the RMP, a cycle or chain with positive-reduced cost must be found and added to the RMP to improve its current objective value. To find such a column(s) one can solve tailored \textit{pricing problems}, which return either a ``positive-price" cycle (chain) or a certificate that none exists. RMP and pricing problems are solved iteratively, typically, until strong duality conditions are satisfied. Since the solution of the RMP may not be integer, column generation is embedded into a branch-and-bound algorithm to obtain an optimal solution, yielding a \BP \ algorithm.
 
\subsection{Restricted master problem}

We use the linear programming relaxation of the disaggregated cycle formulation \eqref{objDisIDCF},  proposed in \cite{Klimentova2014} as our RMP, \rev{which starts with zero columns}. In this formulation, the input digraph $\digraphc$ is disaggregated into copies of two types: cycles and chains. The goal is to obtain a feasible cycle and a feasible chain from each of such copies, starting with a unique vertex. \rev{This unique vertex in a chain copy is an \single\, and its copy of the graph consists of the chains triggered by it. The unique vertex is called a \emph{feedback vertex} for cycle copies, and the graph copy consists of all cycles to which this feedback vertex belongs.  We define $\fvs \subseteq \Vertex$ as the set of feedback vertices in the digraph. As an example, consider the digraph shown in Figure \ref{fig:PCompleteGraphv0}. If we set $\fvs = \{4,5\}$, then two cycle copies (Figures \ref{fig:FirstCopyv0} and \ref{fig:SecondCopyv0}) include all cycles of the input digraph. More details on the use of feedback vertices follow in Section \ref{BuildingMDDs}.}

The number of cycle copies is governed by an upper bound $\vert \fvs \vert$ on the number of cycles in a feasible KEP solution, and the number of chain copies is given by the number of \singles. Shortly, our goal will become to find those graph copies (for cycles and chains) and their associated MDD.

Thus, for every cycle copy indexed by $ \fbcopiescy = \{1,...,\vert \fvs \vert\}$ and for every chain copy indexed by $\fbcopiesch = \{1,...,\vert \singleset \vert\}$, there exists a cycle (chain) decision variable $\zcycleb$ ($\zchainb$) associated to the $\fm$th graph copy. \rev{Letting $\cfeaset$ be the set of cycles in the $i$th graph copy limited to cycles with up to $K$ exchanges, and $\pfeaset$ defined similarly for chains, the disaggregated cycle formulation of the RMP is}
\begin{subequations}
	\label{sub:discyclefv0} 
	\begin{align}
	\max \quad&  \sum_{\fm \in \fbcopiescy}  \sum_{\cycle \in \cfeaset} \cweight\zcycleb + \sum_{\fm \in \fbcopiesch}  \sum_{\chain \in \pfeaset} \pweight\zchainb  \label{objDisIDCF}\tag{IDCF}\\
	&\sum_{\fm \in \fbcopiescy} \sum_{\cycle \in \cfeaset: \vix \in \Vertex(\cycle)} \zcycleb + \sum_{\fm \in \fbcopiesch} \sum_{\chain \in \pfeaset: \vix \in \Vertex(\chain)}  \zchainb \le 1 && \vix \in \Vertex \label{eq:onepervertexIDCF}\\
	&  \zcycleb \ge 0 && \fm \in \fbcopiescy, \cycle \in \cfeaset\\
	&  \zchainb \ge 0 && \fm \in \fbcopiesch, \chain \in \pfeaset \label{eq:BinaryChainIDCF}
	 \end{align}
\end{subequations}

The objective value maximizes the weighted sum of feasible cycles and chains, whereas Constraints \eqref{eq:onepervertexIDCF} restrict every vertex to be selected in at most one cycle or chain, and thus, in at most one graph copy. 

\begin{figure}[tbp]
	\vskip 0.2cm
	\begin{subfigure}{.30\textwidth}
		\centering
		% include first image
			\tikzstyle{place}=[circle,draw=blue!50,fill=blue!20,thick,
			inner sep=0pt,minimum size=4mm]
			\tikzstyle{transition}=[circle,draw=black!50,fill=black!20,thick,
			inner sep=0pt,minimum size=4mm]
			\begin{tikzpicture}
			\node[place] (waiting) at (0,2) {$4$}; 
			\node[transition] (critical) at (0,1) {$3$};
			\node[place] (semaphore) at (0,0) {$5$};
			\node[transition] (leave critical) at (1,1)  {$1$}; 
			\node[transition] (enter critical) at (-1,1) {$2$};
			%\draw [->] (enter critical) to (critical);
			\draw [-{Classical TikZ Rightarrow}] (waiting) to [bend left=45] (leave critical);
			\draw [-{Classical TikZ Rightarrow}] (waiting) to [bend right=45] (enter critical);
			\draw [-{Classical TikZ Rightarrow}] (enter critical) to [bend right=45] (semaphore);
			\draw [-{Classical TikZ Rightarrow}] (semaphore) to [bend right=45] (waiting);
			\draw [-{Classical TikZ Rightarrow}] (waiting) to [bend right=45] (semaphore);
			\draw [-{Classical TikZ Rightarrow}] (semaphore) to [bend right=45] (leave critical);
			\draw [-{Classical TikZ Rightarrow}] (waiting) to [bend right=45] (critical);
			\draw [-{Classical TikZ Rightarrow}] (critical) to [bend right=45] (waiting);
			\draw [-{Classical TikZ Rightarrow}] (semaphore) to [bend right=45] (critical);
			\draw [-{Classical TikZ Rightarrow}] (critical) to [bend right=45] (semaphore);
			\draw [-{Classical TikZ Rightarrow}] (leave critical) to  (semaphore);
			\end{tikzpicture}
		\caption{$\digraphc$}	
		\label{fig:PCompleteGraphv0}
	\end{subfigure}
	\begin{subfigure}{.30\textwidth}
		\centering
		% include second image
			\centering
			\tikzstyle{place}=[circle,draw=blue!50,fill=blue!20,thick,
			inner sep=0pt,minimum size=4mm]
			\tikzstyle{transition}=[circle,draw=black!50,fill=black!20,thick,
			inner sep=0pt,minimum size=4mm]
			\begin{tikzpicture}
			\node[place] (waiting) at (0,2) {$4$}; 
			\node[transition] (critical) at (0,1) {$3$};
			\node[transition] (semaphore) at (0,0) {$5$};
			\node[transition] (leave critical) at (1,1)  {$1$}; 
			\node[transition] (enter critical) at (-1,1) {$2$};
			%\draw [->] (enter critical) to (critical);
			\draw [-{Classical TikZ Rightarrow}] (waiting) to [bend left=45] (leave critical);
			\draw [-{Classical TikZ Rightarrow}] (waiting) to [bend right=45] (enter critical);
			\draw [-{Classical TikZ Rightarrow}] (enter critical) to [bend right=45] (semaphore);
			\draw [-{Classical TikZ Rightarrow}] (semaphore) to [bend right=45] (waiting);
			\draw [-{Classical TikZ Rightarrow}] (waiting) to [bend right=45] (semaphore);
			\draw [-{Classical TikZ Rightarrow}] (semaphore) to [bend right=45] (leave critical);
			\draw [-{Classical TikZ Rightarrow}] (waiting) to [bend right=45] (critical);
			\draw [-{Classical TikZ Rightarrow}] (critical) to [bend right=45] (waiting);
			\draw [-{Classical TikZ Rightarrow}] (semaphore) to [bend right=45] (critical);
			\draw [-{Classical TikZ Rightarrow}] (critical) to [bend right=45] (semaphore);
			\draw [-{Classical TikZ Rightarrow}] (leave critical) to  (semaphore);
			\end{tikzpicture}
			\caption{$\digraphd^{1} = (\Vertexd^{1}, \Arcsd^{1})$}
		    \label{fig:FirstCopyv0}
	\end{subfigure}
	\begin{subfigure}{.30\textwidth}
		\centering
		% include third image
			\centering
			\tikzstyle{place}=[circle,draw=blue!50,fill=blue!20,thick,
			inner sep=0pt,minimum size=4mm]
			\tikzstyle{transition}=[circle,draw=black!50,fill=black!20,thick,
			inner sep=0pt,minimum size=4mm]
			\tikzstyle{whites}=[circle,draw=white,fill=white,thick,
			inner sep=0pt,minimum size=4mm]
			\begin{tikzpicture}
			\node[whites] (waiting) at (0,2) {}; 
			\node[transition] (critical) at (0,1) {$3$};
			\node[place] (semaphore) at (0,0) {$5$};
			\node[transition] (leave critical) at (1,1)  {$1$}; 
			\node[whites] (enter critical) at (-1,1) {};
			\draw [-{Classical TikZ Rightarrow}] (semaphore) to [bend right=45] (leave critical);
			\draw [-{Classical TikZ Rightarrow}] (semaphore) to [bend right=45] (critical);
			\draw [-{Classical TikZ Rightarrow}] (critical) to [bend right=45] (semaphore);
			\draw [-{Classical TikZ Rightarrow}] (leave critical) to  (semaphore);
			\end{tikzpicture}	 
			\caption{$\digraphd^{2} = (\Vertexd^{2}, \Arcsd^{2})$}
		    \label{fig:SecondCopyv0}
	\end{subfigure}
	\caption{Graph disaggregation. Blue vertices correspond to feedback vertices. \rev{(a) Complete digraph; (b) vertex 4 is the feedback vertex; (c) vertex 5 is the feedback vertex}}
	\label{fig:Copiesv0}
\end{figure}

\subsection{MDDs for the Pricing Problem}

%Finding an assignment of vertices to a graph copy while minimizing penalization is equivalent to finding a positive-price cycle or chain. Therefore, whether we take the primal or dual problem of our RMP, subproblems can take the form of \eqref{sprbm:cycles} and \eqref{sprbm:chains}. Notice that subproblems can be formulated differently (see, e.g., \cite{Anderson2015, Duncan2019}), or not being solved through MIP at all, provided that there is a faster algorithm. 
The pricing problem in this context corresponds to finding a feasible cycle or a feasible chain with positive reduced cost. If the reduced cost of a cycle or chain is positive (thus, a positive-price column), then such column should be added to the RMP. The reduced cost of a cycle and chain can be defined, for example, by the objective function of MIP formulations $\eqref{sprbm:cycles}$ and $\eqref{sprbm:chains}$, respectively, (see Online Appendix \ref{sec:Lagrange}). Solution methods to the pricing problem in the literature for the cycle-only version include heuristics \citep{Abraham2007},  MIP  \citep{Roth2007}, or a combination of exact and  heuristic algorithms \citep{Klimentova2014, Lam2020}. Also, \cite{Glorie2014} and \cite{Plaut2016} tackled this version by using a modified Bellman-Ford algorithm in a reduced graph, {\cdionne for the cycle-only version \citep{Plaut2016b}}. Pricing algorithms including long chains have been less studied. \cite{Dickerson2019} studied heuristics for generating positive-price short cycles and arbitrarily large chains when their expected utility is maximized. Our goal is to solve the pricing problems efficiently for both long cycles and long chains by means of MDDs.
%{\cdionne \cite{Dickerson2019} improved the current solver of  a transplant community in the United States, UNOS, based on the previous work of \cite{Abraham2007}. The authors studied heuristics for generating positive-price short cycles and arbitrarily large chains when their expected utility is maximized. Tested instances had at most 1000 \pairsAb\ and 100 \singles. Their results show that 38 out of 128 instances with 900 \pairsAb\ and 90 \singles\ could be solved within an hour.}  Next, we show how to solve the pricing problems by means of MDDs for the general version of the KEP.

\subsubsection{Decision Diagrams} A decision diagram is a graphical data structure,  used in optimization to represent the solution set of a given problem. Particularly, a decision diagram is a rooted, directed and acyclic layered graph, where (if \textit{exact}), every path from the root node $\mathbf{r}$ to a terminal node $\mathbf{t}$ has a one-to-one correspondence with a solution in the feasible space of the optimization problem. If the objective function is arc separable, then a shortest-path-style algorithm can be used to find the optimal objective value. When the domain of decision variables represented in the diagram includes three or more values, the decision diagram is an MDD. \cite{Cire2013} proposed solving sequencing problems using MDDs and showed primary applications in scheduling and routing, (see also \citet{Kinable2017, Castro2020}). The only work we are aware of, in which MDDs are used to solve pricing problems, is by \cite{Arvind2018}. They tackle a last-mile transportation problem in which passengers reach their final destination by using a last-mile service system linked to the terminal station of a mass transit system.  %It is assumed that a desired arrival time of each passenger is known in advance. For each destination, the authors build MDDs for the subset of passengers who share the same destination. They distinguish between one-arcs and zero-arcs, to refer to a passenger being aboard a vehicle or not, respectively. 
A path in these MDDs represents a partition of passengers sharing trip characteristics. 
%Instead of a partition of elements, we find feasible cycles or chains in each graph copy.

%\textcolor{red}{Decision diagrams have been widely used in a variety of applications (see, e.g., \cite{Castro2020, Kinable2017}) showing successful results.  Although in integer programming, decision diagrams are widely used as a cut generation method, our interest is to used them as a solution framework for pricing problems in column generation.} 

\subsubsection{Construction of MDDs for the KEP}

\label{BuildingMDDs}
  Decision diagrams can be constructed by finding the \textit{state transition graph} of a dynamic programming (DP) formulation  and reducing it afterwards (see \cite{Cire2013, Hooker2013} for more details). We model the pricing problems by formulating two DP models; for cycles and chains. Particularly, a DP model is formulated for each cycle copy $\fm \in \fbcopiescy$, an for each chain copy $\fm \in \fbcopiesch$. A \textit{state} in these models represents the vertices $\vix \in \Vertex$ visited at previous \textit{stages}, where a stage corresponds to the position of a vertex in a cycle or chain. \textit{Controls} $\hat{\cydpvar}^{\vert \Kcycle \vert}$ and $\bar{\cydpvar}^{\vert \Lchain \vert} $ take the index value of a vertex  $\vix \in \pairset$ and vertex $\vix \in \Vertex$, indicating that it is assigned to a cycle or chain, respectively, in the position given by the control index, i.e., in position $\Kix \le \Kcycle$ of a cycle or position $\Lix \le \Lchain$ of a chain. Since the domain of  $\hat{\cydpvar}$ and $\hat{\cydpvar}$ variables contains three or more values, the resulting diagram is an MDD.
 
\begin{algorithm}[tbp]
    \spacing
	\SetAlgoLined
	\KwResult {Feedback vertex set $\fbset \subseteq \Vertex$ and decision diagrams $\mdd \ \forall \fm \in \fbcopiescy$}
	%\begin{algorithmic}[1]
		\KwIn {Digraph $\digraphc$, vertex ordering $\orderedSet$, maximum length $\Kcycle$}
		$\fbset = \emptyset$,
		$\fbcopiescy = \emptyset$,
		$i = 0$\\
		\While{$\mid \orderedSet \mid \ge 1$}{
			$\vix^{*} \gets \FirstEorderedSet \in \orderedSet$\\
			$V \gets \Vertex \setminus \fvs $\\ % \mid d_{\vix^{*},\vix} + d_{\vix,\vix^{*}} \le \Kcycle
			$A \gets \{(\uix, \vix) \in \Arcs: \uix, \vix \in V\}$\\ %  d_{\vix^{*},\uix} + 1 + d_{\vix,\vix^{*}} \le \Kcycle
			$\basicMDD \gets \texttt{BuildMDD}(\vix^{*}, D = (V , A))$\\
			\If{$\mathcal{\hat{A}} \neq \emptyset$}{
				$i \gets i + 1$, 
				$\fbcopiescy \gets \fbcopiescy \cup \{i\}$\\
				$\mdd \gets \basicMDD$\\
				$\fbset \gets \fbset \cup \{\vix^{*}\}$\\
			}
			$\orderedSet \gets \orderedSet\setminus \{\vix^{*}\}$\\
		}
		%\KwOut {$\fbset, \digraphd^{\fm} \ \forall \fm \in \fbcopiescy$}
		\caption{$K$-limited FVS heuristics}
		\label{FVSAlg0}
	%\end{algorithmic}
\end{algorithm}

While it is easy to derive the number of MDDs for chains (one for every \single), that is not the case for cycles. As a result, we start by considering $\fvs$ as a ``$K$-limited'' \textit{feedback vertex set} (FVS), i.e., a set of vertices whose removal from the graph removes all feasible cycles. We refer to vertices in this set as \textit{feedback vertices}.  Note that when $\Kcycle$ is sufficiently large, a $K$-limited FVS corresponds to an FVS of $\digraph$. The goal is then to find a ``small'' cardinality $K$-limited FVS to build as few MDDs for cycles as possible; one for each vertex in $\fvs$. Since finding a minimum FVS and a smallest set covering are NP-hard problems \citep{Karp1972}, we introduce a heuristics as an alternative.
	
Algorithm \ref{FVSAlg0} provides a mechanism to obtain a $\Kcycle$-limited FVS  $\fbset = \{\vixstar_{1},...,  \vixstar_{ \vert \fbset \vert}\} \subseteq \pairset $ and cycle MDDs   $\mdd$, indexed by $ \fbcopiescy = \{1,...,\vert \fvs \vert\}$.  The procedure takes as input a digraph $\digraph$ and an ordered set $\orderedSet$ of \pairsAb\ in $\digraph$, sorted according to a vertex-ordering rule (e.g., maximum in-degree of a vertex). While there is at least one element in $\orderedSet$, the  first vertex in that set denoted by $\FirstEorderedSet$, is selected as the next potential feedback vertex $\vix^{*}$. Subsequently, the function \texttt{BuildMDD} creates a transition graph starting at $\vix^{*}$ and reduces it to an MDD.  If the set of arcs $\mathcal{\hat{A}}$ is not empty, then a new cycle MDD is created. Vertex $\vix^{*}$ is removed from $\orderedSet$ afterwards, and a new iteration starts. 

 The MDD $\mdd$ for the $\fm$th cycle copy of $\digraph$ has its node set $\mnodes$ partitioned into $\Kcycle$ layers, $\layer_{1},...,\layer_{\Kcycle}$, corresponding to the decision of which \pairAb\ belongs to the $\Kix$th  position in a cycle, denoted by $\pos_{\Kix}$, plus two terminal layers $\layer_{\Kcycle + 1}$ and $\layer_{\Kcycle + 2}$ representing the completion of a feasible cycle. Layers $\layer_{1}$,  $\layer_{2}$, $\layer_{\Kcycle + 1}$ and $\layer_{\Kcycle + 2}$ have a single node each.  Every arc $\ma \in \marcs$ has an associated label $val(\ma) \in \Vertexd^{\fm}$ such that $\pos_{\Kix} = val(\ma)$ corresponds to assigning vertex $val(\ma)$  to the $\Kix$th position of a cycle. Since a cycle in the $\fm$th cycle copy starts with the $\fm$th feedback vertex, then $\pos_{1} = \vixstar_\fm$. Thus, an arc-specified  path ($\ma_1, ..., \ma_{\Kix}$) from nodes $\mathbf{r}$ to $\mathbf{t}$ defines the \pairAb\ sequence $(\pos_{1},,...,\pos_{\Kix}) = (val(\ma_1), ..., val(\ma_{\Kix}))$, equivalent to a feasible cycle in $\digraph$. On the other hand, the set of nodes $\mnodesch$ for the $\fm$th MDD of a chain copy, $\mddch$,  is partitioned into $\Lchain + 2$ layers, with a single node in the first and last layer, representing the start and end of a chain, respectively. Recall that for a chain to involve $\Lchain$ transplants, $\Lchain + 1$ vertices $\vix \in \Vertex$ are required. Likewise, an arc $\ma \in \marcsch$ has a label $\pos_{\Lix} = val(\ma)$ indicating the vertex $\vix \in \Vertexbar^{\fm}$ at the $\Lix$th position in a chain, noticing that $\pos_{1}  = \uix_{\fm}$. A path starting at the root node $\mathbf{r}$ and ending at a node on the third layer or higher represents a feasible chain, since its length is at least one. 

Figure \ref{sub:exactMDD} depicts all possible sequences of \pairsAb\ $\pos_{1}, \pos_{2},...,\pos_{\Kix}$ encoding feasible cycles on the graph copy of Figure \ref{fig:FirstCopyv0}. A value $\pos_{\Kix}$ placed on an arc corresponds to the $\Kix$th vertex $\vix \in \pairset$ in a cycle covered by vertex $4$. For instance, the path $(\mathbf{r},\dgv_1, \dgv_3, \dgv_6, \mathbf{t})$ in Figure \ref{sub:exactMDD} encodes the cycle $\cycle = \{4,3,5,4\}$ in Figure \ref{fig:FirstCopyv0}. %{\cdionne See \citet{Cire2013} for details on MDD construction.}

\begin{figure}[tbp]
	\vskip 0.2cm
	\centering
	\begin{adjustbox}{minipage=\linewidth,scale=0.8}
		\begin{subfigure}{.45\textwidth}
			\centering
			\begin{tikzpicture}
			%\node[circle, fill=black, minimum  size=2pt, label=above:{$r$}] (0) at (3, 9) {};
			
			\begin{scope}[every node/.style={circle, fill=beaublue, inner sep=1pt,minimum size=2pt}]   % style={circle, fill=black}
			%\node [label=above:{$\mathbf{r}$}](0) at (3, 8) {$\{\emptyset\}$};
			\node (0) at (3,8.5) {$\mathbf{r}$};
			\node (1) at (3,7) {$\dgv_{1}$};
			\node (2) at (1,5.5) {$\dgv_{2}$};
			\node (3) at (2,5.5) {$\dgv_{3}$};
			%\node (4) at (4,6) {$\dgv_{4}$};
			\node (5) at (5,5.5) {$\dgv_{4}$};
			%\node (10) at (4,4.5) {$\dgv_{6}$};
			\node (11) at (5,4) {$\dgv_{5}$};
			\node (20) at (3,2.5) {$\dgv_{6}$};
			\node (21) at (3,1) {$\mathbf{t}$};
			\end{scope}
			
			\node (P1) at (-0.25, 7.75) {$\pos_1$};
			\node (P2) at (-0.25, 6.25) {$\pos_2$};
			\node (P3) at (-0.25, 4.75) {$\pos_3$};
			\node (P4) at (-0.25, 3.25) {$\pos_4$};

			%position auto left right above under
			\begin{scope}[>={Stealth[black]},
			every node/.style={scale=0.8}],
			every edge/.style={draw=black, thick}],
			\path [->] (0) edge node[auto] {$4$} (1);
			\path [->] (1) edge[bend right = 20] node[above] {$5$} (2);
			\path [->] (1) edge[bend right = 10] node[above] {$3$} (3);
			\path [->] (1) edge[bend right=20] node[below] {$2$} (5);
			\path [->] (1) edge[bend left = 20] node[above] {$1$} (5);
			\path [->] (3) edge[bend right=20] node[left] {$5$} (20);
			\path [->] (5) edge[bend left = 10] node[auto] {$5$} (11);
			\path [->] (1) edge[bend right=20] node[left] {$5$} (20);
			\path [->] (1) edge[bend left=20] node[right] {$3$} (20);

			%\path [-] (2) edge[bend right=30, draw=beaublue, ultra thick, line width=1.7mm, opacity=0.6] node[auto, below] {} (20);
			\path [->] (2) edge[bend right=35] node[left] {$3$} (20);
			\path [->] (5) edge[bend right = 10]  node[right, above] {$5$} (20);
			\path [->] (11) edge[bend left = 20] node[right, below] {$3$} (20);
			\path [->] (20) edge node[auto] {$4$} (21);
			\end{scope}
			\end{tikzpicture}
			\caption{{\cdionne Vertex 4 is a \pairAb. Exact decision diagram, $\Kcycle = 4$}}
			\label{sub:exactMDD}
			%\end{figure}	
		\end{subfigure}
		\begin{subfigure}{.5\textwidth}
			\centering
			\begin{tikzpicture}
			\node[circle, fill=beaublue, inner sep=1pt,minimum size=2pt] (0) at (3, 8.5){$\mathbf{r}$};
			%\node[circle, fill=beaublue, inner sep=1pt,minimum size=2pt] (21) at (3, 1){$\mathbf{t}$};
			
			\begin{scope}[every node/.style={circle, fill=beaublue, inner sep=1pt,minimum size=2pt}]   % style={circle, fill=black}
			%\node (0) at (3,8.5) {$\mathbf{r}$};
			\node (1) at (3,7) {$\dgv_{1}$};
			\node (5) at (3,5.5) {$\dgv_{4}$};
			\node (11) at (3.5,4) {$\dgv_{5}$};
			\node (20) at (3,2.5) {$\mathbf{t}$};
			%\node (21) at (3,1) {$\mathbf{t}$};
			\end{scope}
			
			\node (L0) at (4.75,8.4) {$\layer_1$};
			\node (L1) at (4.75, 7) {$\layer_2$};
			\node (L2) at (4.75, 5.5) {$\layer_3$};
			\node (L3) at (4.75, 4) {$\layer_4$};
			\node (L4) at (4.75, 2.5) {$\layer_5$};
			\node (L5) at (4.75, 1.2) {$\layer_6$};

			%position auto left right above under
			\begin{scope}[>={Stealth[black]},
			every node/.style={scale=0.8}],
			every edge/.style={draw=black, thick}],
			\path [->] (0) edge node[auto] {$4$} (1);
			\path [->] (1) edge[bend right=20] node[left] {$2$} (5);
			\path [->] (1) edge[bend left = 20] node[right] {$1$} (5);
			\path [->] (5) edge[bend left = 10] node[auto] {$5$} (11);
			\path [->] (5) edge[bend right = 25]  node[right, left] {$5$} (20);
			\path [->] (11) edge[bend left = 15] node[right, right] {$3$} (20);
			%\path [->] (20) edge node[auto] {$\tau$} (21);
			\end{scope}
			\end{tikzpicture}
			\caption{{\cdionne Vertex 4 is an \single. Restricted decision diagram, $\Lchain = 3$}}
			\label{subfig:MDDchains}
		\end{subfigure}
	\end{adjustbox}
	\caption{MDDs for the example of Figure \ref{fig:FirstCopyv0}.}
	\label{fig:MDDs}
\end{figure}

Since exact MDDs can grow exponentially large, it might be necessary to limit {\cdionne their} size, turning them into \textit{restricted} decision {\cdionne diagrams. A} decision diagram is called restricted if the set of solutions corresponding to its $\mathbf{r}$-$\mathbf{t}$ paths is a subset of the entire feasible set of solutions. Figure \ref{subfig:MDDchains} shows a restricted MDD as if vertex $4$ in Figure \ref{fig:FirstCopyv0} were an \single. In this example, the MDD is restricted to have chains including only two out of the four vertices receiving an arc from vertex $4$; namely, vertices $1$ and $2$. 

\subsubsection{Finding a positive-price column}
\label{sec:PosPriceCol}
Since the reduced cost of a cycle (or chain) in the $\fm$th MDD, $\hat{r}^{\fm}_{\cycle}$ (or $\bar{r}^{\fm}_{\chain}$), is arc separable \citep{Glorie2014}, we proceed to show how to find positive-price columns via MDDs:

Let $\inn(\dgv_{s} ) \subset \marcs$ be the set of incoming arcs to a node $\dgv_{s} $ in $\mnodes$ and $\ell(\ma)$ the layer index of the source node of arc $\ma$, e.g., in Figure \ref{sub:exactMDD} $\inn(\dgv_4) = \{(\dgv_1, \dgv_4)^{1}, (\dgv_1, \dgv_4)^{2}\}$ and $\ell((\dgv_1, \dgv_4)^{1}) = 2$, where the superscripts distinguish the two arcs coming into $\dgv_{4}$, one with $\pos_{2} = {1}$ and the other with $\pos_{2} = {2}$. Moreover, let $\multi$ be the dual variables of Constraints \ref{eq:onepervertexIDCF}. We define the recursive function values of an arc $\ma = (\dgv_{s}, \dgv_{s'})$ for the $\fm$th MDD for cycles and chains, $\recucy(\ma)$ and $\recuch(\ma)$, respectively, as the maximum reduced cost of all paths ending at $\ma$:
\begin{subequations}
	\label{eq:recu1a}
	\begin{empheq}[left={\recucy(\ma) = \recuch(\ma) = \empheqlbrace\,}]{align}
	& 0 &  \dgv_{s} =  \mathbf{r}
	\label{eq:1a} \\
	& \max_{\ma' \in \inn(\dgv_{s})}  \left \{  \recucy(\ma') + w_{val(\ma'), val(\ma)} - \multi_{val(\ma')} \right \}  & \mbox{otherwise}
	\label{eq:2a} 
	\end{empheq}
\end{subequations}

The recursive function \eqref{eq:recu1a} is valid by Bellman's principle of optimality since the reduced cost of a cycle or chain is arc separable and the portion taken by every arc only depends on the previous \pairAb\ or \single\ in the sequence. Thus, the maximum reduced cost of a cycle and chain, $\recucy$ and $\recuch$, respectively, is given by 
\begin{subequations}
\label{sub:recucych}
	\begin{align}
	\recucy= &\max  \left \{0, \recucy((n,\mathbf{t})) \right \}&   \fm \in \fbcopiescy
	\label{eq:1b} \\
	\recuch= & \max  \left \{0, \max_{\ma \in \marcsch: \ell(\ma) \ge 3} \recuch(\ma)  - \multi_{val(\ma)} \right \}& \fm \in \fbcopiesch
	\label{eq:2b} 
    \end{align}
\end{subequations}
\rev{where for a terminal node $\mathbf{t}$ in a cycle, $n$ is the only node on the $\Kcycle + 1$ layer with an arc pointing to $\mathbf{t}$ (e.g., node $n_6$ in Figure \ref{sub:exactMDD}), providing a terminal condition for the recursion for cycles.}

Recursion \eqref{eq:1b} computes the maximum reduced cost of a cycle at the terminal node $\mathbf{t}$, since all paths need to reach $\mathbf{t}$ to close it up. For chains, on the other hand, any portion of a path in $\mgraphch$ is a feasible path, for which the longest path can be found at any layer of the MDD where the length of a chain (in terms or arcs in $\Arcsbar^{\fm}$) is at least one. The term subtracted in {\cdionne Equation} $\eqref{eq:2b}$ captures the dual variable of the last pair in a chain.  Since  $\multi_{\vix} \in \mathbb{R}_{+}$ for all $\vix \in \Vertex$ (see \eqref{objLag2opt} in the Online Appendix), thus, if the dual variable of a given vertex $\vix \in \Vertex$ is large enough to lead to a negative-price path at node $\mathbf{t}$, we may need to cut that path short at some previous vertex in the sequence to obtain a positive-price chain. For instance, in Figure \ref{subfig:MDDchains} consider $\multi_{2} = \multi_{3} = 5$, all the other dual variables set to zero and {\cdionne $w_{\pairarc} = 1$} for all $\arcp \in \Arcs$. Sequence $(4,1,5)$ representing a 2-length chain in $\digraphbar^{\fm}$ is contained in $(4,1,5, 3)$. Clearly, the former  yields the highest reduced cost of a chain in Figure \ref{subfig:MDDchains}. Thus, $\recuch = \recuch((\dgv_4, \dgv_5)) = 2$. Next, we show a series of results on the complexity of computing a positive-price column via MDDs. The corresponding proofs are presented in the Online Appendix \ref{sec:Proofs}.

\setcounter{theorem}{3}
\begin{prop}
	Given the reduced costs $\hat{r}_{\cycle}^{\fm}$ and $\bar{r}^{\fm}_{\chain}$ expressed as an arc-separable function for all $(\dgv_{s}, \dgv_{s'})  \in \marcs$ and $(\dgv_{s}, \dgv_{s'})  \in \marcsch$, a positive-price cycle, if one exists, can be found in  time $\mathcal{O}(\sum_{\fm \in \fbcopiescy} \sum_{(\dgv_{s}, \dgv_{s'})  \in \marcs}  \vert \inn(\dgv_{s}) \vert)$. Similarly, a positive-price chain can be found in $\mathcal{O}(\sum_{\fm \in \fbcopiesch} \sum_{(\dgv_{s}, \dgv_{s'})  \in \marcsch}  \vert \inn(\dgv_{s}) \vert)$.
\end{prop}

%\proof{Proof.}
%For every arc $(\dgv_{s}, \dgv_{s'}) \in \marcs$ and $(\dgv_{s}, \dgv_{s'}) \in \marcsch$, $\vert \inn((\dgv_{s}, \dgv_{s'})) \vert$ comparisons need to be performed to obtain $\recucy((\dgv_{s}, \dgv_{s'}))$ and $\recuch((\dgv_{s}, \dgv_{s'}))$, respectively in \eqref{eq:recu1a}. Therefore, for the $\fm$th MDD of a cycle copy, $\sum_{(\dgv_{s}, \dgv_{s'}) \in \marcs} \vert \inn(\dgv_{s}) \vert$ comparisons are required to compute $\recucy$, whereas for the $\fm$th MDD of a chain copy,  there are the same number of comparisons plus $\vert \marcsch \vert$ comparisons  of all arcs, in \eqref{eq:2b}, before obtaining $\recuch$. Because there are $\vert \fbcopiescy \vert$ cycle MDDs and $\vert \fbcopiesch \vert$ chain MDDs, it follows that the time complexity is as shown above. \hfill$\square$

\begin{prop}
The size of the input $\sum_{\fm \in \fbcopiescy} \sum_{(\dgv_{s}, \dgv_{s'}) \in \marcs} \vert \inn(\dgv_{s}) \vert$ grows as $\vert \Vertexd^{\fm} \vert ^{\Kcycle + 1}$ does.
\end{prop}

\begin{prop}
	The size of the input $\sum_{\fm \in \fbcopiesch} \sum_{(\dgv_{s}, \dgv_{s'}) \in \marcsch} \vert \inn(\dgv_{s}) \vert$ grows as $\vert \Vertexbar^{\fm} \vert ^{\Lchain + 2}$ does for bounded chains and as $\vert  \Vertexbar^{\fm}  \vert !$ when $\Lchain \rightarrow \infty$.
\end{prop}	
	
%\proof{Proof.}
%A similar reasoning to the previous proposition can be followed, except that a chain can be cut short if by visiting a new vertex $\vix \in \Vertexbar^{\fm}$ in a sequence of the state transition graph at least one \pairAb\ is present more than once, thereby violating the condition of being a simple path. We know that for a path to have $\Lchain$-many arcs, it is necessary to have a sequence with $\Lchain$ \pairsAb, thus, $\vert \marcs \vert$ tends to $\vert \Vertexbar^{\fm} \vert^{\Lchain}$ and $\sum_{\fm \in \fbcopiesch}  \sum_{(\dgv_{s}, \dgv_{s'}) \in \marcs}  \vert \inn(\dgv_{s}) \vert \approx \vert \Vertexd^{\fm} \vert^{\Lchain + 2}$. Therefore, for bounded chains, finding a positive-price chain can be done in time  $\mathcal{O}(\vert \Vertexd^{\fm} \vert ^{\Lchain+ 2})$. The second part follows by the fact that after we visit $\Lix$ vertices, there are still $\vert \Vertexd^{\fm} \vert - \Lix$ ways to choose the next one, until only one can be chosen, thus, the time to find a positive-price column when $\Lchain$ is unbounded is exponential. \hfill$\square$
	
Despite potentially very large diagram sizes in general, there are three reasons for which finding a positive-price column can still be done efficiently in practice. First, even though $\vert \Vertex \vert$ can be large, arc density of $\digraph$ for real KEP instances is below 50\% {\cred \citep{Saidman2006,Anderson2015, Dickerson2016}}. Second, for small values of $\Kcycle$ and $\Lchain$,  it is possible to reduce considerably the size of the input by selecting appropriately a $K$-limited FVS. Lastly, MDDs are reduced significantly after the state transition graph is obtained \citep{Cire2013}.%, e.g., see Figure \ref{sub:exactMDD}. 

\subsection{Branching scheme}
\label{sub:branching}
%Whenever the solution of the RMP  in a node of the branch-and-bound tree is fractional, a branching scheme is applied to find an integer and eventually optimal solution to the problem. Different branching schemes have been investigated in the literature: branching on cycles \citep{Abraham2007}, branching on arcs $ \arcp  \in \Arcs$ \citep{Glorie2014,Lam2020} and branching on arcs $ \arcp \in \Arcs^{l}$ \citep{Klimentova2014}, where $\Arcs^{l}$ is the set of arcs in the $l$th graph copy of $\digraph$ as defined in \citep{Constantino2013}. 
{\cdionne The search tree may have exponential depth when branching is done on possibly every cycle, thus, branching on arcs is usually preferred.} If $\digraph$ is a complete graph, there can be up to $\vert \pairset \vert \vert \pairset - 1 \vert + \vert \pairset \vert \vert \singleset \vert$ arcs in $\Arcs$. On the other hand, branching on arcs in $ \Arcsd^{\fm}$ and $\Arcsbar^{\fm}$ implies that there are up to $(\vert \fbcopiesch \vert + \vert \fbcopiescy \vert ) \vert \pairset \vert \vert \pairset - 1 \vert  + \vert \fbcopiesch \vert \vert \pairset \vert \vert \singleset \vert$ arcs across all graph copies. Among the two arc-based schemes,  branching on arcs in $\Arcs$ results in a lower depth branching tree. We therefore choose this option as our branching scheme. 

Particularly, on every fractional node of the search tree we branch on an arc $\arcp \in \Arcs$ whose {\cdionne fractional value} is closest to $0.5$. That is, we generate two children, one in which the arc is prohibited and another in which the arc is selected. When banning an arc from the RMP, we modify Equation \eqref{eq:2a}  by replacing $w_{\pairarc}$ with a sufficiently large {\cred negative} number $M$. By doing so, the length of any path in any copy traversing that arc approaches negative infinity, thereby ruling it out due to the definition of $\recucy$ and $\recuch$. On the other hand, enforcing an arc $\arcp \in \Arcs$ requires the inclusion of the {\cdionne following} constraint in the RMP:
\begin{align}
\sum_{\fm \in \fbcopiescy} \sum_{\cycle \in \cfeaset: \arcp \in \Arcs(\cycle)} \zcycleb + \sum_{\fm \in \fbcopiesch} \sum_{\chain \in \pfeaset: \arcp \in \Arcs(\chain)}  \zchainb = 1 && (\dualselarcs_{\arcp})
\label{eq:fixingarc}
\end{align}

The addition of constraint \eqref{eq:fixingarc} changes the reduced cost of a chain and cycle. If we let $\enfoarcs \subseteq \Arcs $ be the set of selected arcs in a branch-and-bound node, the reduced cost of a column in the $\fm$th MDD, is now given by
\begin{subequations}
	\label{eq:NewRedCost}
	\begin{align}
	\hat{r}_{\cycle}^{\fm}  = 
	 &\cweight  - \sum_{\vix \in \Vertex(\cycle)}\multi_{\vix} - \sum_{\arcp \in \enfoarcs \cap \Arcs(\cycle)} \dualselarcs_{\arcp} &  \cycle \in \cfeaset
	\label{eq:1rc} \\
	\bar{r}_{\chain}^{\fm}  = 
	 & \pweight  - \sum_{\vix \in \Vertex(\chain)} \multi_{\vix} - \sum_{\arcp \in \enfoarcs \cap \Arcs(\chain)} \dualselarcs_{\arcp} &  \chain \in \pfeaset
	\label{eq:2rc} 
	\end{align}
\end{subequations}

\noindent where $\dualselarcs_{\arcp}$ is the dual variable of  Constraint \eqref{eq:fixingarc}. Thus, if  $(val(\ma^{'}), val(\ma)) \in \enfoarcs$, then $\dualselarcs_{(val(\ma^{'}), val(\ma))}$ can be subtracted from recursive expression  \eqref{eq:2a} to account for Constraint \eqref{eq:fixingarc}  in the RMP.

%In our implementation, we process children nodes in the branching tree in decreasing order of their upper bound. 
 If the solution of the RMP is fractional, we branch, and then apply column generation to the resulting node. After optimally solving the RMP, the upper bound given by its objective value is compared to the best lower bound found. If the former is lower, that branch is pruned. Otherwise, a lower bound is obtained by granting integrality to the decision variables (columns) in the RMP and re-solving it. Whenever a lower bound matches the best upper bound, optimality is achieved. If, due to time limitations, it is not possible to solve the RMP to optimality, \eqref{objLag2main} provides a valid upper bound  on the optimal value  of the KEP, as presented next. %To this end, if $\enfoarcs \ne \emptyset$, the second summation of \eqref{eq:NewRedCost} is subtracted from the objective function of $\eqref{sprbm:cycles}$ and $\eqref{sprbm:chains}$. %We subtract  $\sum_{\arcp \in \enfoarcs \cap \Arcsd^{\fm} } \dualselarcs_{\arcp} \text{ } \forall \fm \in \fbcopiescy$ and $\sum_{\arcp \in \enfoarcs \cap \Arcsbar^{\fm} } \dualselarcs_{\arcp} \text{ } \forall \fm \in \fbcopiesch$ from the objective function of $\eqref{sprbm:cycles}$ and $\eqref{sprbm:chains}$, respectively.

\subsection{A New Upper Bound on the Optimal Objective Value of the KEP}
\label{sec:UpperBound}
\rev{In this section, we propose what to the best of our knowledge is the only other known upper bound on the optimal objective value of the KEP, which is tighter than the existing proposed bound \citep{Abraham2007}, while recent works provide no upper bound for instances that timed out \citep{Lam2020}. The upper bound in \citet{Abraham2007} consists of solving an unrestricted maximum weight matching problem ($K = L = \infty$) on a bipartite graph, which can be solved in polynomial time. Our new upper bound (LR-UB), based on Lagrangian relaxation, is particularly useful when the optimality of the RMP cannot be proven within the time limit:}
\rev{
\begin{subequations}
	\label{sub:LR2main} 
	\begin{align}
	\centering
		&\objLag (\multi) = \sum_{\fm \in \fbcopiesch} \objLagCy^{\fm}({\multi}) + \sum_{\fm \in \fbcopiescy} \objLagCh^{\fm} ({\multi})+ \sum_{\uix \in \Vertex} \multi_{\uix} \label{objLag2main}\tag{LR-UB}
	\end{align}
\end{subequations}
}
\noindent \rev{where $\multi_{\uix}$ corresponds to the dual variables of Constraints \eqref{eq:onepervertexIDCF} in \ref{objDisIDCF}, and $\objLagCy^{\fm}({\multi})$ and $\objLagCh^{\fm} ({\multi})$ correspond to the pricing problems \eqref{sprbm:cycles} and \eqref{sprbm:chains} in Online Appendix \ref{sec:Lagrange}, respectively. We refer the reader to Appendix \ref{sec:Lagrange} for an explanation of how $\multi_{\uix}$ are decision variables in a Lagrangian relaxation and for a formal proof on the validity of \eqref{objLag2main}.}

\section{General Solution Framework}
\label{sec:SolAp}

 In our {\cdionne solution framework}, a combination of exact and restricted MDDs is used so that once built they are stored in computer memory  and retrieved every time pricing problems need to be solved. As a result, we cannot solely rely on MDDs to prove optimality of the RMP. We introduce a three-phase solution framework consisting of a search through MDDs for both cycles and chains (Phase 1), a cutting plane algorithm to search for positive-price chains and cycles, whose final goal is to prove that no more positive-price chains exist (Phase 2), and a two-step search to find a positive-price cycle, if any (Phase 3). Figure \ref{fig:LogicDiagram} illustrates the framework.

 \begin{figure}[ht!] 
	%\FIGURE
	\begin{center}
		\resizebox{\textwidth}{!}{%
			\begin{tikzpicture}[node distance=1.9cm, font=\scriptsize]
			%fill=red!30
			%\tikzstyle{io} = [trapezium, trapezium left angle=70, trapezium right angle=110, minimum width=3cm, minimum height=1cm, text centered, draw=black, fill=blue!30]
			%\tikzstyle{process} = [rectangle, minimum width=3cm, minimum height=1cm, text centered, draw=black, fill=orange!30]
			\tikzstyle{startstop} = [rectangle, rounded corners, minimum width=1.15cm, minimum height=0.7cm,text centered, draw=black, fill = gray!40!white]
			\tikzstyle{process} = [rectangle, minimum width=0.6cm, minimum height=1cm, text centered, text width=1.45cm, draw=black]
			\tikzstyle{decision} = [diamond, aspect=1.9, text centered, text width=1.65cm, draw=black]
			\tikzstyle{arrow} = [line width=0.2mm,->,>=stealth]
			
			\small
			
			%%%%%%%definiciones%%%%%%%
			%arriba
			\node at (0, 0) (nStart) [startstop, xshift=-1.0cm, text width=1.5cm] {START};
			\node (nBuildMdds) [process, below of=nStart, xshift = 0.0cm, text width=1.5cm] {Build MDDs};
			\node (nSolveRMP) [process, right of=nBuildMdds, xshift=1.0cm, text width=1.5cm] {Solve RMP};
			\node (nPerform) [process, right of=nSolveRMP, xshift=1.0cm, text width=1.5cm] {Perform LPA on MDDs};
			\node (nGenera) [process, right of=nPerform, xshift=1.0cm, text width=1.5cm] {Generate cols.};
			\node (dPositive) [decision, right of=nGenera, xshift = 1.0cm] {$\exists \ \cycle^{+}$ or $\chain^{+}$?};
			\node (dummy1) [above of=dPositive, yshift = 0.0cm] {};
			\node (dummy2) [right of=dPositive, xshift = 0.3cm] {};
			
			%medio
			\node (nSolveLPH) [process, below of=dPositive, yshift = -0.7cm, text width=1.5cm] {Solve $\eqref{objLP}$};
			\node (nFindC) [process, left of=nSolveLPH, xshift = -1.0cm, text width=1.5cm] {Find $\cycle \in \allCycles \setminus \loncfeaset$};
			\node (dZmayor) [decision, below of=nFindC, yshift = -0.0cm] {$\objLag^{\mbox{\tiny \eqref{objLP}}} > 0$?};
			\node (nCheck) [process, left of=dZmayor, xshift = -1.0cm, text width=1.5cm] {Check sol. $\lonvar^{*} \ge 0.9$};
			\node (dPosPrice) [decision, left of=nCheck, xshift = -1.0cm] {$\exists \ \cycle^{+}$ or $\chain^{+}$?};
			\node (nAdded) [decision, above of=dPosPrice, yshift = -0.0cm] {\eqref{eq:LongSizeChain}-\eqref{eq:integChain} added?};
			\node (nsolveRMP2) [process, left of=dPosPrice, xshift = -1.0cm, text width=1.5cm] {Solve RMP};
			\node (dummy3) [above of=nAdded, yshift = -0.5cm] {};
			\node (dummy4) [right of=dZmayor, xshift = 1.0cm] {};
			\node (dummy20) [above of=nsolveRMP2, xshift = 0.0cm] {};
			
			%abajo
			\node (nPerformee) [process, below of=dZmayor, yshift = -0.3cm, text width=1.5cm] {Perform ES};
			\node (ElapsedTime) [decision, left of=nPerformee, xshift = -1.0cm] {ET $ > t^e$? };
			\node (dPosPrice2) [decision, left of=ElapsedTime, xshift = -1.0cm] {$\exists \ \cycle^{+}$?};
			\node (nSolveMCC2) [process, below of=ElapsedTime, yshift = -0.0cm, text width=1.5cm] {Solve $\tiny \eqref{objCyclesMIP}$};
			\node (dMCC) [decision, left of=nSolveMCC2, xshift = -1.0cm] { $\objLag^{\mbox{\tiny \eqref{objCyclesMIP}}}  > 0$?};
			\node (nFindPrice3) [process, left of=dMCC, xshift = -1.0cm, text width=1.5cm] {Find $\cycle^{+}$};
			\node (nfin) [startstop, below of=nPerformee, xshift = 2.5cm, text width=1.5cm] {END};
			\node (dummy5) [left of=dPosPrice2, xshift = -1.0cm] {};
			\node (dummy6) [above of=dPosPrice2, yshift = -1cm] {};
			\node (dummy8) [below of=dMCC, yshift = 0.9cm] {};
			
			%%%%%%%flechas%%%%%%%
			%arriba			
			\draw [arrow] (nStart) -- (nBuildMdds);
			\draw [arrow] (nBuildMdds) -- node[anchor=south, xshift = 0.0cm, yshift=-0.0cm] {} (nSolveRMP);
			\draw [arrow] (nSolveRMP) -- (nPerform);
			\draw [arrow] (nPerform) -- (nGenera);
			\draw [arrow] (nGenera) -- (dPositive);
			\draw [-,thick] (dPositive) -- node[anchor=west, xshift = 0.0cm, yshift=-0.0cm] {yes} (dummy1);
			\draw [arrow] (dummy1.south) -| (nSolveRMP);
			\draw [-,thick] (dPositive) -- node[anchor=north, xshift = 0.1cm, yshift=-0.0cm] {no} (dummy2);
			\draw [arrow] (dummy2.west) |- (nSolveLPH);
			
			%medio
			\draw [arrow] (nFindC) -- (nSolveLPH);
			\draw [arrow] (nAdded) -- node[anchor=south, xshift = 0.0cm, yshift=-0.0cm] {yes} (nFindC);
			\draw [-,thick] (nAdded) -- node[anchor=west, xshift = 0.0cm, yshift=-0.0cm] {no, add $\eqref{eq:LongSizeChain}\eqref{eq:integChain}$} (dummy3.center);
			\draw [arrow] (dummy3.center) -| (nSolveLPH);
			\draw [arrow] (dPosPrice) -- node[anchor=west, xshift = 0.0cm, yshift=-0.0cm] {no} (nAdded);
			\draw [arrow] (nCheck) -- (dPosPrice);
			\draw [arrow] (dZmayor) -- node[anchor=south, xshift = 0.0cm, yshift=-0.0cm] {yes} (nCheck);
			\draw [-,thick] (nSolveLPH) |- (dummy4.center);
			\draw [arrow] (dummy4.center) -- (dZmayor);
			\draw [arrow] (dPosPrice) -- node[anchor=south, xshift = 0.0cm, yshift=-0.0cm] {yes} (nsolveRMP2);
			\draw [-,thick] (nsolveRMP2) -- node[anchor=west, xshift = 0.0cm, yshift=-0.0cm] {} (dummy20.center);
			\draw [arrow] (dummy20.center) -- node[anchor=west, xshift = 0.0cm, yshift=0.2cm] {} (dummy3);
			
			%abajo
			\draw [arrow] (dZmayor) -- node[anchor=west, xshift = 0.0cm, yshift= -0.24cm] {no} (nPerformee);
			\draw [arrow] (nSolveMCC2) -- (dMCC);
			\draw [arrow] (nFindPrice3) -- (dummy5.center);
			\draw [-,thick] (dPosPrice2) -- node[anchor=south, xshift = 0.0cm, yshift=-0.0cm] {yes} (dummy5.center);
			\draw [arrow] (dummy5.center) -- (nsolveRMP2);
			\draw [-,thick] (dPosPrice2) -- node[anchor=east, xshift = 0.8cm, yshift=-0.1cm] {no} (dummy6.center);
			\draw [arrow] (dummy6.center) -| (nfin);
			\draw [arrow] (nPerformee) -- (ElapsedTime);
			\draw [arrow] (ElapsedTime) -- node[anchor=west, xshift = 0.0cm, yshift=-0.0cm] {yes} (nSolveMCC2);
			\draw [arrow] (dMCC) -- node[anchor=south, xshift = 0.2cm, yshift=-0.0cm] {yes} (nFindPrice3);
			\draw [-,thick] (dMCC) -- node[anchor=east, xshift = 0.7cm, yshift=0.04cm] {no} (dummy8.center);				
			\draw [arrow] (dummy8.center) -| (nfin);
			\draw [arrow] (ElapsedTime) -- node[anchor=south, xshift = 0.0cm, yshift=0.1cm] {no} (dPosPrice2);
			
			%%%%%%%grupos%%%%%%%				
			%arriba
			\node [thick, xshift=0.0cm, draw=black, rounded corners=2mm, dotted, fit = (nSolveRMP) (nPerform) (nGenera) (dPositive) (dummy1) ] (nPhase1) {};
			\node (p1) [above of=nPhase1, xshift = -3.7cm, yshift = -0.6cm] {\textbf{Phase 1: PPs via MDDs}};
			
			%medio
			\node [thick, xshift=0.0cm, draw=black, rounded corners=2mm, dotted, fit = (nSolveLPH) (nFindC) (dZmayor) (nCheck) (dPosPrice) (nAdded) (nsolveRMP2) (dummy3)] (nPhase2)  {};
			\node (p2) [above of=nPhase2, xshift = -4.9cm, yshift = 0.1cm] {\textbf{Phase 2: PPs via $\eqref{objLP}$}};
			
			%abajo
			\node [thick, xshift=0.0cm, draw=black, rounded corners=2mm, dotted, fit = (nPerformee) (ElapsedTime) (dPosPrice2) (nSolveMCC2) (dMCC)(nFindPrice3) (dummy6) (dummy8)] (nPhase3) {};
			\node (p3) [above of=nPhase3, xshift = -3.6cm, yshift = 0.1cm] {\textbf{Phase 3: PPs via 2SP}};
			\end{tikzpicture}
		}%
		\caption{Flow chart of column generation for $\Kcycle \ge 3$, $\Lchain \in \mathbb{Z_{+}}$. PPs stands for pricing problems, LPA for longest-path algorithm, cols. for columns, ES for exhaustive search, ET for elapsed time, $\cycle^{+}$ ($\chain^{+}$) for positive-price cycle (chain), and 2SP for the two-step procedure described in Phase 3.}
		\label{fig:LogicDiagram}
	\end{center}
\end{figure}

\subsection{Phase 1: Solving the pricing problems via MDDs}
Building the MDDs is the first  step. We parameterize some aspects  to make a reasonable usage of computer memory {\cdionne (Section \ref{sec:Results})}. Particularly, if $\Kcycle \ge 4$ and $\vert \pairset \vert > 500$,  we build restricted MDDs by considering a maximum cycle length of 3 in 90\% of the decision diagrams, while in the remaining 10\% we keep the true value of  $\Kcycle$. If $\Kcycle \le 4$ and $\vert \pairset \vert \le 500$, we build exact MDDs. When constructing a transition state graph for chains, we explore at most $20\%$ of vertices receiving an arc from $\vix \in \Vertex$, unless $\vert \singleset \vert > 250$, in which case we explore only $10\%$. In all cases, the maximum length of chains (in terms of arcs) considered in the construction of MDDs is also 3, regardless of the true value of $\Lchain$. After the MDDs are built, we store them in memory and use them to solve the pricing problems as depicted by Figure \ref{fig:LogicDiagram}, in integration with Phases 2 and 3.

\subsection{Phase 2: A longest path formulation for chains and cycles}
We use a longest path problem as a relaxation of \eqref{sprbm:chains} in which the goal is to find an $\mathbf{s}\mbox{-}\mathbf{t}$ path or a cycle, both corresponding to feasible positive-price columns.  To this end, let us define $\longraph$ as a digraph whose vertex set $\lonVertex = \Vertex \cup \{\mathbf{s},\mathbf{t}\} $ has two dummy vertices $\mathbf{s}$ and $\mathbf{t}$ such that the set of arcs $\lonArcs = \Arcs \cup \{(\mathbf{s}, \uix) \cup  (\vix, \mathbf{t}) \mid  \uix \in \singleset, \vix \in \pairset\}$ connects dummy vertex $\mathbf{s}$ to \singles\, and \pairsAb\ to  vertex $\mathbf{t}$. For arcs including a dummy vertex, their weight is set to zero.  Moreover, let $\allCycles$ be the set of all simple cycles and $\loncfeaset \subseteq \allCycles$ be the set of feasible cycles  in $\lonD$.  As defined before, $\enfoarcs$ is the set of selected arcs. A decision variable  $\lonvar_{\pairarc}$ takes the value 1 if arc $\arcp \in \lonArcs$ is selected, and 0 otherwise. Thus, a {\cdionne relaxed} longest path formulation at some node in the branching tree is
\begin{subequations}
	\label{sub:LongPath} 
	\begin{align}
	\objLag^{\mbox{\footnotesize \eqref{objLP}}} := \ \max& \sum_{\arcp \in \lonArcs: \uix \ne \mathbf{s}} (w_{\pairarc} - \multi_{\uix}) \lonvar_{\pairarc}  - \sum_{\arcp \in \enfoarcs} \dualselarcs_{\arcp} \lonvar_{\pairarc} & \label{objLP}\tag{LPH}\\		
	\text{s.t.} &\sum_{\vix:\arcp \in \lonArcs} \lonvar_{\pairarc}  - \sum_{\vix: \inarcp \in \lonArcs} \lonvar_{\inpairarc} = \left\{
	\begin{array}{@{}ll@{}ll@{}}
	1, &&  \uix= \mathbf{s}\\
	0, &&   \uix \in \Vertex\\
	-1, && \uix = \mathbf{t}
	\end{array}\right.
	\label{sub:LongBalance}\\
	&\sum_{\arcp \in \lonArcs} \lonvar_{\pairarc} \le 1 \qquad\qquad\qquad\qquad\quad \uix \in \Vertex
	\label{sub:MaxFlow}\\
	&\lonvar_{\pairarc} \in [0,1] \qquad \qquad\qquad\qquad\quad\hspace{2mm} \arcp\in\lonArcs
	\end{align}
\end{subequations}

Although a solution of \eqref{objLP} may lead to a path using more than $\Lchain$ arcs, or a solution with subtours, or a non-integer solution, we see these downsides as an opportunity to either find efficiently a positive-price chain (cycle) missed in the first phase or prove that none exists. Because \eqref{objLP} is a relaxation of  \eqref{sprbm:chains}, whenever the objective value of \eqref{objLP} is zero, so is the objective value of \eqref{sprbm:chains}. {\cdionne Note that even without subtour-elimination constraints, a solution of \eqref{objLP} guarantees a  path (may be fractional) from vertex $\mathbf{s}$ to $\mathbf{t}$ due to flow-balance Constraints \eqref{sub:LongBalance}. The solution may also have subtours representing positive-price cycles useful for the RMP. Particularly, if $\lonvar^{*}$ is an optimal solution of \eqref{objLP}, we check arcs $\arcp \in \lonArcs$ for which $\lonvar_{\pairarc}^{*} \ge 0.9$ when searching for positive-price columns. Lastly,  enforcing two dummy arcs, one going out of $\mathbf{s}$ and one coming into $\mathbf{t}$, requires us to adjust the right-hand side of Constraints \eqref{eq:LongSizeChain}, resulting in the following constraints:}
\begin{subequations}
	\label{sub:OntheFly}
\begin{align}
	\sum_{\arcp \in \cycle} \lonvar_{\pairarc}  \le&\  \vert \cycle \vert - 1 &\cycle  \in \allCycles \setminus \loncfeaset &&
	\label{eq:subcons}\\
	\sum_{\arcp \in \lonArcs} \lonvar_{\pairarc}  \le&\  \Lchain + 2&&
	\label{eq:LongSizeChain}\\
	\lonvar_{\pairarc} \in&\ \{0,1\}& \arcp \in \lonArcs&&
	\label{eq:integChain}
\end{align}
\end{subequations}

{\cdionne Therefore, the goal is to solve first a linear program without  Constraints \eqref{sub:OntheFly}, check the solution for positive-price and \textit{feasible} columns, and only then add \eqref{sub:OntheFly} if need be. Experimentally, we observed that the ``warmed-up'' dual variables resulting from Phase 1 allow us to relax integrality constraints and yet obtain an integer solution in many cases. Figure \ref{fig:LogicDiagram} shows how \eqref{objLP} $+$ {\cdionne Constraints} \eqref{sub:OntheFly} is solved via a cutting plane algorithm integrated with the other two phases during column generation.}

\subsection{Phase 3: A two-step procedure for cycles}
\label{sec:2steps}
If a positive-price cycle still exists in $\digraph$ but not found in Phases 1 and 2, we perform an exhaustive enumeration of cycles while the time limit is not exceeded. Otherwise,  the exhaustive search ends and a MIP is solved instead. Note that Phase 2 guarantees to find any positive-price chain, if does exists. Therefore, in Phase 3, we only have to search positive-price cycles.

\subsubsection{Exhaustive search}
\label{sec:exsearch}
We perform a depth-first search on every cycle copy $\digraphd^{\fm}$, where a feasible cycle is rooted at $\vixstar_{\fm} \in \fvs$. First, we sort and search over the graph copies in increasing order of their $\multi_{\vix}$ values. Next, we traverse $\digraphd^{\fm}$, and when $\vixstar_{\fm}$ is the leaf node of a path from the root and it is found at a position between $3$ and $\Kcycle + 1$, that path constitutes a feasible cycle $\cycle$. If $\hat{r}_{\cycle}^{\fm}> 0$, the cycle has a positive reduced cost  and it is sent to the RMP. \cite{Abraham2007} and \cite{Lam2020} implemented a similar search, although unlike them, our paths are rooted at vertices in FVS. Despite of searching cycles on trees rooted only at a subset of vertices, exhaustive enumeration becomes a bottleneck for large instances. Therefore, whenever a time threshold, $t^e$, is surpassed we shift to solving the MIP provided next.

\subsubsection{A MIP for cycles}
\eqref{objCyclesMIP} {\cdionne finds} a feasible cycle in $\digraph$ with maximum reduced cost at {\cdionne some branch-and-bound node. If found, the cycle is sent to the RMP. Note that Constraint \eqref{eq:CycleSizePh3} guarantees at most $\Kcycle$ selected arcs, whether they are distributed in multiple (smaller) cycles or not. Thus, there is no need for subtour-elimination constraints.  \eqref{objCyclesMIP} is formulated as follows:}
\begin{subequations}
	\label{sub:CyclesMIP} 
	\begin{align}
	\objLag^{\mbox{\footnotesize \eqref{objCyclesMIP}}}= \ \max &\sum_{\arcp \in \Arcs} (w_{\pairarc} - \multi_{\uix}) \lonvar_{\pairarc}  - \sum_{\arcp \in \enfoarcs} \dualselarcs_{\arcp} \lonvar_{\pairarc}  \label{objCyclesMIP}\tag{MCC}\\		
	\text{s.t.} &\sum_{\vix:\arcp \in \Arcs} \lonvar_{\pairarc}  - \sum_{\vix: \inarcp \in \Arcs} \lonvar_{\inpairarc} = 0& \uix \in \Vertex\\\
	&\sum_{\vix: \arcp\in \Arcs} \lonvar_{\pairarc}  \le  \left\{
	\begin{array}{@{}ll@{}}
	1, & \quad \ \uix \in \pairset\\
	0, & \quad \ \uix \in \singleset
	\end{array}\right.\\
	&\sum_{\arcp \in \Arcs} \lonvar_{\pairarc}  \le \Kcycle &\label{eq:CycleSizePh3} \\
	&\lonvar_{\pairarc} \in \{0,1\} &\arcp \in \Arcs&
	\end{align}
\end{subequations}

%Figure \ref{fig:LogicDiagram} shows the solution of the RMP via the three phases discussed before. 

\subsubsection{Algorithmic details}
\label{sec:AlgorithmicDetails}
 In the three-phase solution framework given in Figure \ref{fig:LogicDiagram}, after building the MDDs, the pricing problems are solved to optimality as follows.  During Phase 1, we compute $\recucy$ and $\recuch$ for all $\fm \in \fbcopiesch $ and $\fm \in \fbcopiescy$, respectively, and add the positive-price columns found to the RMP after every iteration. We encourage the use of chains by returning up to 15 positive-price chain columns to the RMP, whereas only one from every cycle MDD. Note that in {\cdionne both cases}, recursions \eqref{eq:1b}-\eqref{eq:2b} only need to be computed once. In Phase 2, we delay the inclusion of Constraints \eqref{sub:OntheFly}, until their addition is mandatory to find a positive-price chain column. Every time a solution is checked during Phase 2, a positive-price cycle column is searched when failing to find a chain column. As for Phase 3, we set $t^{e} = 20$s to exhaustively find a positive-price cycle or reach the end with none. If the time threshold is exceeded, we proceed to solve \eqref{objCyclesMIP} and resolve the RMP, if needed. When there are no \singles\ present in the input graph $\digraph$, we simply skip Phase 2. Likewise, since MDDs are exact when $\Kcycle \le 4$, $\vert \pairset \vert \le 500$ and no \singles, Phases 2 and 3 are also dropped.  %In other cases, pricing problems could be solved, e.g.,  via MIP, but such implementation is out of the scope of our algorithm. 
 Lastly, note that the procedure given in Figure \ref{fig:LogicDiagram} can be easily adapted to the case $\Lchain = \infty$, since it suffices to remove Constraint \eqref{eq:LongSizeChain}.

\section{Computational Experiments}
\label{sec:Results}
%In this section we present  computational experiments comparing our solution framework to the state of the art.
%\begin{enumerate} 
%	\item Accumulated percentage of optimally solved instances over the time limit, including and not including pre-processing steps
%	\item Compare number of variables among the three approaches, point out the memory usage.
%	\item For every set and K,L combination present the number of optimal solutions, feasible solutions and out of memory instances.
%	\item For every set and K,L combination present the average number of columns added, whether they were chain or cycle columns and in which phase they were found, the number of iterations at phase 2 and 3, and the average time at every optimization phase.
%\end{enumerate}

\rev{Our code base and test instances are publicly available at \url{https://github.com/d-aleman/KEP-BPMDD}}. MDDs as  well as our \BP \  algorithm (BP\_MDD) are implemented in C++ and experiments are conducted on a machine with Debian GNU/Linux as operating system and a 3.60GHz processor Intel(R) Core(TM) {\cdionne with 120 GB RAM, of which we allocated a maximum of 8 GB for BP\_MDD (including the construction of MDDs) and 60 GB for state-of-the-art solution methods.} CPLEX 12.10 is used as the LP/MIP solver. {\cdionne We investigate the extent up to which leading approaches could scale in practice, therefore, we allocated a larger RAM capacity.}  BP\_MDD is compared against the state-of-the-art PICEF and HPIEF solution methods \citep{Dickerson2016}. The PICEF and HPIEF solvers are retrieved from the original authors, where HPIEF is the model with full reduction. They call Gurobi 7.5.2 to solve MIP models. Although different LP/MIP solvers may add noise to the analysis, in the latest history of benchmark sets, Gurobi is the lead \citep{Mittelman2020}. It is worth noting that \cite{Lam2020} hold the state-of-the-art solver for the cycle-only version. They tested their algorithm on the same library as ours, PrefLib \citep{Mattei2013}, achieving total run-times up to  {\cdionne 33s} for instances with 2048 \pairsAb\ and mostly less than a second for instances with 1024 \pairsAb, when $\Kcycle = 3$. {\cdionne For $\Kcycle = 4$ with 1024 \pairsAb, their maximum run-time among the 9 instances optimally solved was 21.8min. The total run-time of our algorithm, although small, is larger than that reported by \cite{Lam2020}, particularly when $\Kcycle = 3$.} {\cdionne However, since \cite{Lam2020} did not consider larger values of $\Kcycle$, and their approach is for cycles only, we do not compare to them in this analysis.} When $\Lchain = 0$, the formulation by \cite{Duncan2019} reduces to the cycle formulation, and so does PICEF. Therefore, we compare  BP\_MDD with HPIEF and PICEF, which we refer to as ``solvers''. %As will be shown, our algorithm does have advantages besides the fact that it can also tackle long chains.

\subsection{Instances}
For our subsequent analysis, we used the PrefLib repository \citep{Mattei2013}, whose instances were generated by \cite{Dickerson2012} based on data from KPDPs in the United States. The first group, referred to as KEP\_N0, has $80$ instances, {\cdionne split into $8$ subsets, each with $10$ instances and no \singles;} and hence, only cycles are considered. In this group, $\vert \pairset \vert \in \{16, 32, 64, 128, 256, 512, 1024, 2048\}$ and their arc density varies from 10\%  to 32\%. 

{\cdionne The second group of instances, referred to as KEP\_N, is split into $23$ subsets with $10$ instances in each. \singles\ are present, thus, both cycles and chains are considered in the solution. For these instances, $ 16 \le \vert \pairset \vert \le 2048$, $1 \le \vert \singleset \vert \le 307$} and the arc density is within $23 \%$ and $46\%$. 
%Let $g_{d_N} = \frac{\vert \Arcs \vert}{\vert \pairset \vert \vert \pairset -1 \vert + \vert \pairset \vert \vert \singleset \vert}$ be the graph density for instances in KEP\_N. The first sum in the {\cdionne denominator} is the number of possible arcs among \pairsAb, and the second sum is the number of possible arcs from \singles\ to \pairsAb.  {\cdionne We have that $g_{d_N}$% The first sum in the divider is the number of possible arcs among \pairsAb, whereas the second sum is the number of possible arcs from singleton donors to \pairsAb. Recall from Section \ref{sec:PblmD} that \singles\ have no incoming arcs.
%Let the tuple $(\vert \pairset \vert, \vert \singleset \vert)$ characterize each subset. Then, each subset in KEP\_N has a one-to-one  correspondence with the set of tuples $\{(17,1), (18,2), (33,1), (35,3), (36,4), (67,3), (70,6), (73,9), (134,6), (140,12), (147,19), \\ (268,12), (281,25), (294,38), (537,25), (563,51), (588,76), (1075,51), (1126,102), (1177,153), \\ (2150,102), (2252,204), (2355,307)\}$.

\subsection{Computational performance}
PICEF, HPIEF and BP\_MDD are evaluated on the two groups of instances. Each instance in KEP\_N0 is solved for $K \in \{3,4,5,6\}$, totaling 320 runs per solver. For KEP\_N , we set $\Kcycle \in \{3,4\}$ and  $\Lchain \in \{3,4,5,6\}$. The total number of runs per solver in this second group is very large,  $10\times23\times2\times4$ = 1840 runs.  Therefore, we proceed as follows: For each run, if the RAM usage exceeds the limit, the solver stops and aborts the rest of instances in that run's subset, which we presume would lead to the same memory issues. Regardless of the group, {\cdionne for every run we set a time limit of 30min.  Only one thread was used for all the experiments and the rule to find the $\Kcycle$-limited FVS is the maximum in-degree.} {\cred An analysis on multiple criteria to obtain a $\Kcycle$-limited FVS from Algorithm \ref{FVSAlg0} and their implications on performance and scalability is presented in the Online Supplement \ref{sec:OtherResults}.}

%\subsubsection{Computational performance}

\begin{figure}[tbp]
		\centering
		% include first image
		\includegraphics[width=\linewidth]{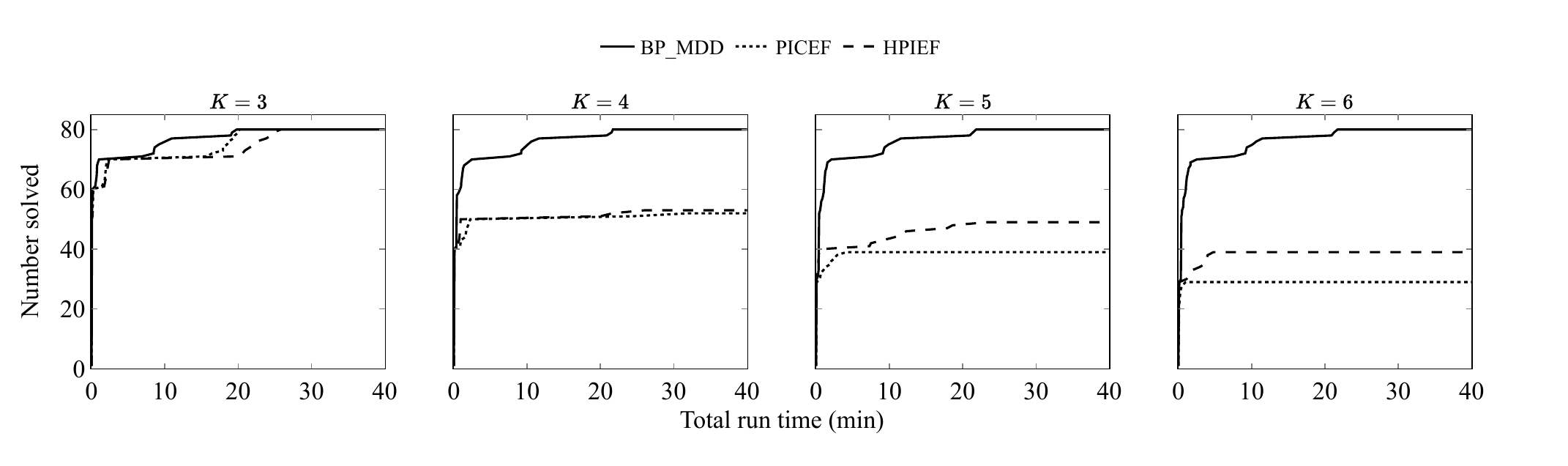}  
		\caption{Performance profiles for the set KEP\_N0.}
		\label{fig:ProfileCycles}
\end{figure}

Figure \ref{fig:ProfileCycles} plots the number of instances in KEP\_N0 solved to optimality up to discrete points in time before reaching the {\cdionne $30$-min limit}, which shows the outperformance of our algorithm over the state of the art. The time reported includes both the total MDD construction time and the \BP\ algorithm time. %The $x$-axis is extended up to $40$ minutes to account for the  MDD construction time. 

When $\Kcycle = 3$, BP\_MDD and PICEF solve all 80 instances in under 20min, followed by HPIEF under 26min. { \cdionne When $\Kcycle \in \{4,5,6\}$, BP\_MDD solves all 80 instances under 25min. Both PICEF and HPIEF perform poorly as $\Kcycle$ increases. RAM usage of BP\_MDD  did not exceed 4 GB in these experiments, while PICEF and HPIEF surpassed $60$ GB. Instances that run out of memory accumulated more than $37$ million variables.} By definition PICEF is exponential in the number of  variables, thus, this result is not surprising. However, HPIEF is polynomial in terms of the number of variables and constraints, even more so the  full-reduction HPIEF we compare our algorithm against,  and yet dimensionality is clearly a challenge.%Additional information on these results is provided in Appendix[].  as opposed to \cite{Lam2020}, whose solution method was able to solve instances up to $K = 4$ and 1024 \pairsAb

 In 94\% of runs, BP\_MDD solves instances to optimality at the root node and at most 8 nodes are explored in the branch-and-bound tree. On average, for instances with $\vert \pairset \vert \ge 1024$, 55\% of the total run-time accounts for solving the pricing problems, with a minimum and maximum percentage of 19\% and 84\%, respectively. \rev{Phase 3 is responsible on average for 27\% (range: $[5\%,72\%]$) of the pricing time when $\Kcycle \ge 4$}. On the same runs ($\vert \pairset \vert \ge 1024$ and $\Kcycle \ge 4$), the average and maximum time for building the MDDs is 109.8s and 257.5s, respectively. For only 10 instances with $\Kcycle \ge 4$, PICEF and HPIEF report a feasible solution within the time limit, meaning that for the others the RAM threshold is exceeded.  The optimality gap reported by these solvers varies between  0\% and 30\%, with respect to the best solution found among all the three solvers. Overall, PICEF is slightly better than HPIEF in terms of optimality gap but not in terms of the number of instances that could fit into memory.%The remaining time is spent among the RMP and solving the cycle formulation with columns present in the RMP. 

\begin{figure}[tbp]
	\centering
	% include first image
	\includegraphics[scale=0.4]{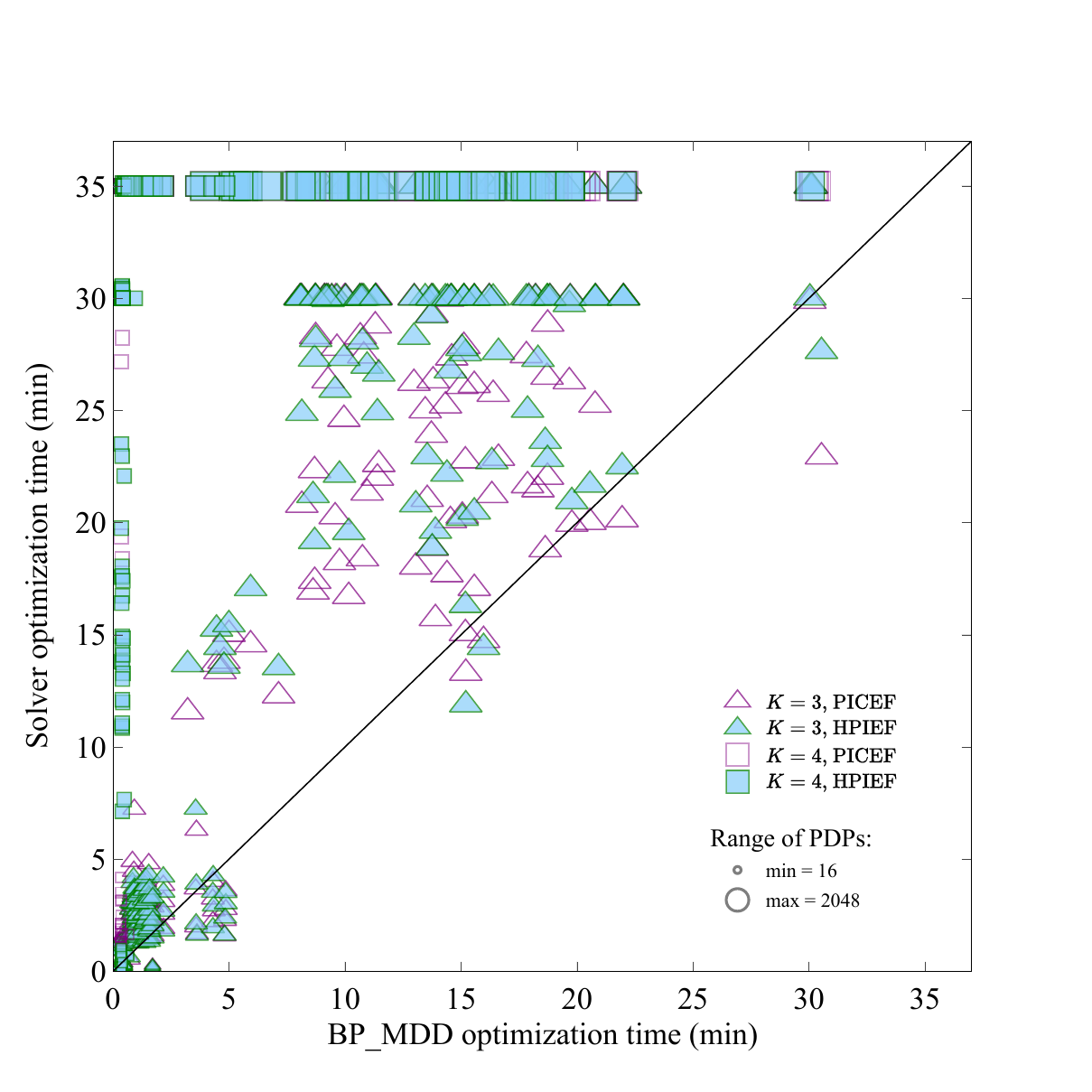}  
	\caption{Performance comparison when $\Lchain \in \{0,3,4,5,6\}$ and $\Kcycle=3$ (triangles) and $\Kcycle=4$ (squares). The size of the markers indicate the number of \pairsAb\ in the instances. Markers located at the $y$-axis value of 35 on the top of the plot indicate instances that were not solved to optimality within the time limit.}
	\label{fig:BP_RunTime.pdf}
\end{figure}

{\cdionne Figure \ref{fig:BP_RunTime.pdf} shows the solution time (without preprocessing or MDD construction)  taken for every run and the state-of-the-art solver when $K \in \{3,4\}$ for all chain lengths. Markers above the diagonal indicate BP\_MDD is faster. Particularly,  BP\_MDD is faster in $20\%$ of runs when $\vert \pairset \vert < 512$, yet the maximum run-time taken for BP\_MDD to solve these instances is $<4$s. When $\vert \pairset \vert \ge 512$, BP\_MDD is faster in $80\%$ of runs.}%and in 42\% of cases, either PICEF or HPIEF do not return an optimal solution.} %Even though the cumulative result among solvers at the end of the time limit is the same when $\Kcycle = 3$, 

\begin{figure}[tbp]
	\centering
	% include first image
	\includegraphics[width=\linewidth]{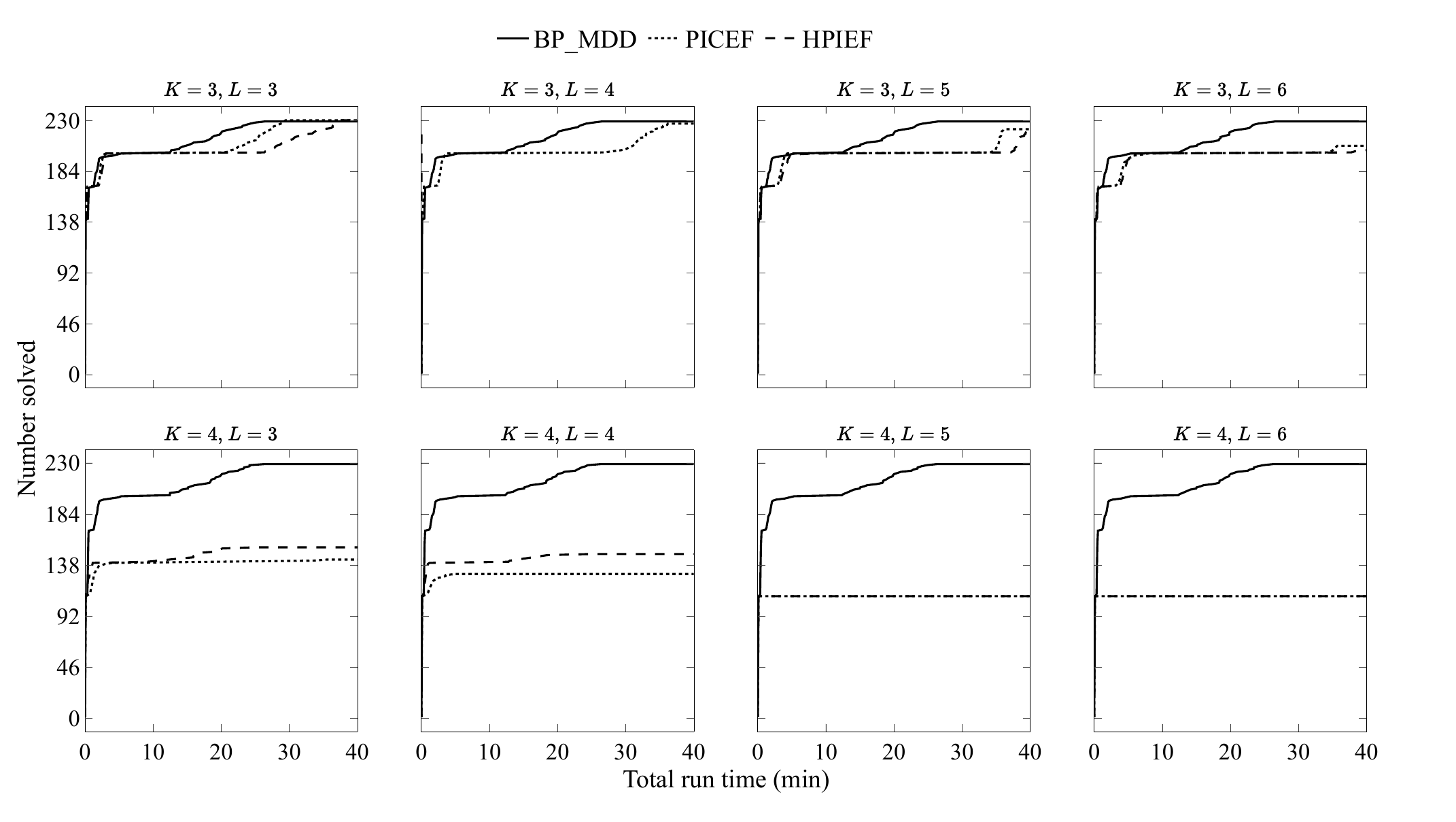}  
	\caption{Number of instances solved for set KEP\_N.}
	\label{fig:ProfileChains}
\end{figure}

Figure  \ref{fig:ProfileChains} shows the  performance profiles for the instances in KEP\_N solved to optimality by the three solvers, which demonstrates that our algorithm outperforms the others, especially when long chains are allowed. {\cdionne For $K = 3$ and $\Lchain = 3$, the performance at the time limit is similar for the three algorithms, although BP\_MDD does not solve one instance (2048 \pairsAb\ and 307 \singles) to optimality.} For that instance, BP\_MDD's optimality gap is 0.2\%. The same instance remains suboptimal across the other $\Kcycle$-$\Lchain$ combinations. {\cdionne The performance of PICEF and HPIEF decreases as $\Kcycle$ and $\Lchain$ increases.} The maximum optimality gap provided by PICEF and HPIEF, among the {\cdionne 95} suboptimal runs at the time limit is 11\%.  Among these two solvers, HPIEF provided 67\% of the {\cdionne suboptimal solutions.}

{\cdionne On average, BP\_MDD solved 94\% of runs to optimality at the root node, with a maximum of two explored branch-and-bound nodes. For instances with $\vert \pairset \vert \ge 1024$, solving the pricing problems accounts for 66\% of the total run-time, and on average 53\% of the pricing time is spent in Phase 1.} For the same instances, the time spent on building the MDDs accounts on average for 25\% of the total run-time, with a maximum of 35\%. On average, in more than 90\% of the cases, all \singles\ were used in the final solution. %{\cred The best solutions found by BP\_MDD for each setting can be found in the Online Supplement.}
Lastly, we note that the majority of columns across all runs are found in Phase 1 (Online Supplement \ref{sec:OtherResults}).
 
 {
 	\cred
\subsection{Sensitivity analysis on the solution structure}

In this section, we first investigate the benefit of long chains for different instance structures, according to the presence of highly sensitized patients.  In the second part of our analysis, we show the change in the composition of different optimal solutions when cycles are prioritized over chains and vice versa. To the best of our knowledge,  BP\_MDD is the only exact algorithm for which such a preference can be selected a priori.

\subsubsection{The value of long chains}

An important question from a practical perspective is, whether the same optimal matching value achieved when allowing chains can be obtained by relying only on cycles, and thus disregarding the need of chains. \cite{Ashlagi2012} and \cite{Dickerson2012} showed analytically that, when the graph is sufficiently sparse, chains benefit the number of matched patients. A standard methodology to study this practical question is the use of analytical models. These models, generate random graphs, e.g., vertices and arcs, by  following probability estimates on inherent characteristics to KPDP participants, such as a patient's sensitization degree and the number of highly-sensitized patients in the graph. For our further analysis, we use a similar random graph model to that of \cite{Ashlagi2012} and \cite{Ding2018}. Particularly, we also assume the existence of two patient-donor pair categories: highly-sensitized and low-sensitized. Instead of building completely random graphs, we take a sample of instances from KEP\_N and preserve arcs $(\uix, \vix) \in \Arcs$  with probability $p_h$ if the patient (and thus the pair) in vertex $\vix$ is highly-sensitized or probability $p_\ell$ if that patient is not. Our goal is to perform a sensitivity analysis on the original instances. Therefore, in Figure  \ref{fig:ValueOfChains}  we use the notation $(\sigma, p_\ell)$ to represent an instance where the proportion of low-sensitized pairs over the total number of pairs is $\sigma$ and $p_\ell$ is the compatibility probability for low-sensitized pairs. In all cases, we set $p_h = 5\%$, taken from \cite{Ding2018}.

We selected 3 small-size instances (indexed by S-181, S-182 and S-183) and 3 medium-size instances (indexed by M-201, M-202, M-203) from the group KEP\_N. The former have 256 pairs and 38 singleton donors whereas the latter have 512 pairs and 25 singleton donors. For each instance, we generated 6 random graphs according to the tuples $(\sigma, p_\ell)$ depicted in Figure  \ref{fig:ValueOfChains}. The average graph densities from left to right are $4.2\%, 6.8\%, 5.6\%, 9.2\%, 8.9\%$ and $15.7\%$. Overall, the higher the percentage of low-sensitized patients, the higher the arc density.

{\cred These instances were solved within an optimality gap of $ \leq 3\%$, with 93.5\% solved under $ \leq 1\%$ optimality gap. Overall, these instances are substantially more challenging. BP\_MDD solved optimally 57\% of them in 151.9 seconds, on average. The feasible solutions took on average 1850.3 seconds, with the largest solution time exceeding an hour}. 

\begin{figure}[tbp]
	\centering
	% include first image
	\includegraphics[width=\linewidth]{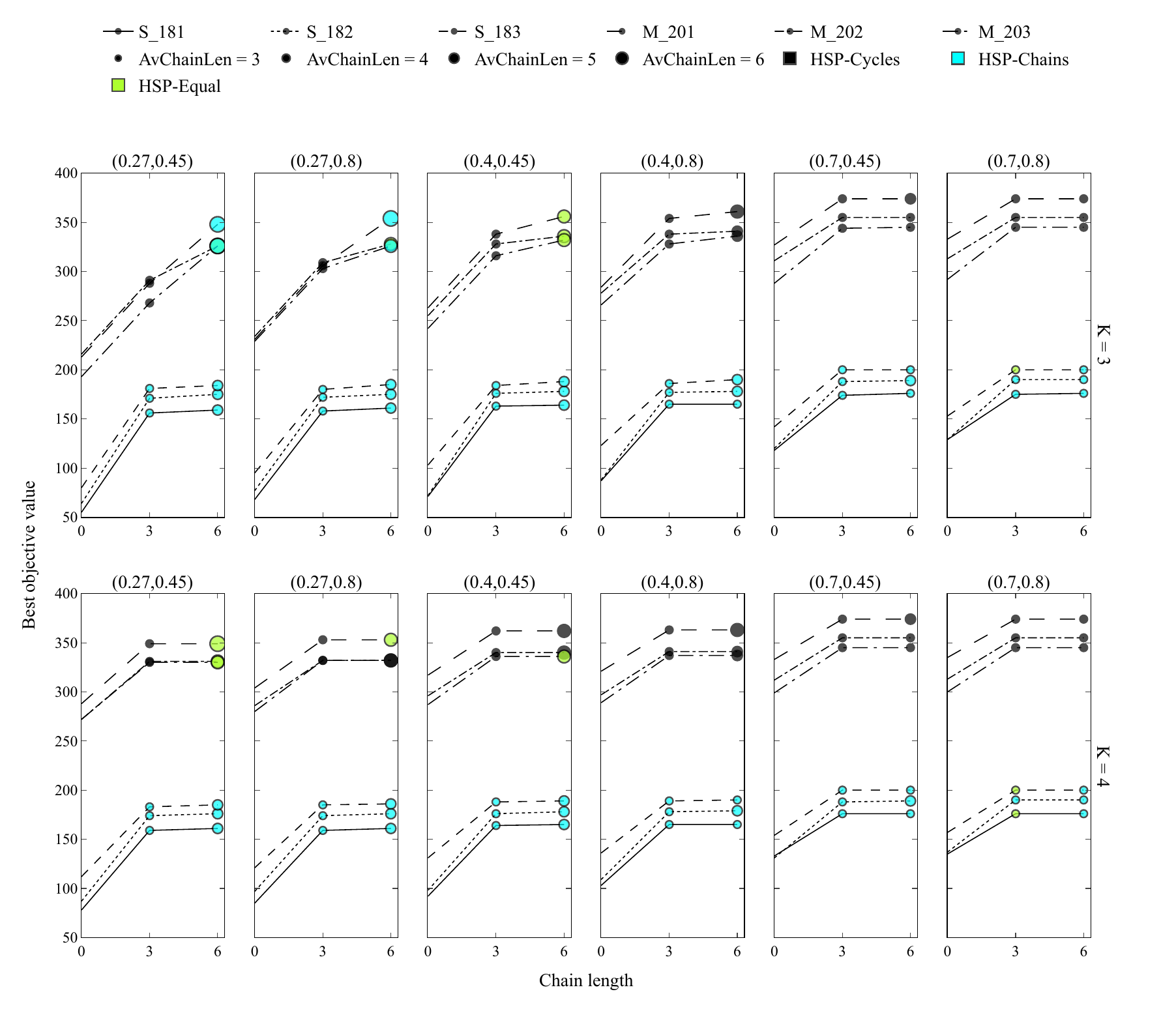}  
	\caption{\cred Change in the number of matches by sensitization category and chain length. AvChainLen: average chain length of an instance solution, rounded to the closest integer; HSP-Chains (Cycles): chains (cycles) included a higher percentage of the matched and highly sensitized patients; HSP-Equal: difference between the percentage of matched and highly-sensitized patients in chains and cycles is less than $5\%$.}
	\label{fig:ValueOfChains}
\end{figure}

The empirical results depicted in Figure \ref{fig:ValueOfChains} complement previous findings on the benefit of chains \citep{Ashlagi2012, Dickerson2012,Ding2018}. For sparse graphs, chains, particularly long chains (involving more than three transplants) are able to reach more patients when the opportunities of cyclic exchanges are scarce. In all cases there is a significant benefit from a cycle-only solution to one allowing chains of size at least three. Chains become shorter, stabilizing around three, as the percentage of low-sensitized patients increase.  This behavior is also observed for the vast majority of instances in KEP\_N, which are significantly denser than the ones presented in Figure \ref{fig:ValueOfChains}.  On average, each instance in KEP\_N sees an increase of {\cred 20\%} when passing from a solution with $\Kcycle = 3$ and no chains to one with $\Kcycle = \Lchain = 3$, and about {\cred 19\%} from no chains and $\Kcycle = 4$ to $\Kcycle = 4, \Lchain = 3$. From this point on, increasing the size of chains beyond three, seems only significant for a handful of the small and medium size instances.

For the small size instances in Figure \ref{fig:ValueOfChains}, allowing chains seem to be critical to include the difficult to match patients in a solution, even when the arc density increases. When the length of cycles is four,  chains still provide benefit, although it seems that an important number of pairs that were reached before by long chains can be matched in these long cycles. Interestingly, chains seem to have a more dramatic effect as the size of the instances increase.

\subsubsection{Solution composition under  prioritization by exchange type}
\label{sec:prioritization}

Another relevant question is, whether for the same optimal matching value, one can opt for a solution that prioritizes the presence of chains over cycles or vice versa. For instance, chains can be prioritized over cycles since even when failures occur in arcs or vertices of the graph, the original chain can simply turn into a shorter one, whereas a cyclic exchange falls apart altogether \citep{Klimentova2016, Dickerson2019}.

With this concern in mind, we added two modes under which BP\_MDD can be executed. The default  mode, referred to as ``CH'', prioritizes chains over cycles and corresponds to returning  up to 15 positive-price chain columns to the RMP, as described in Section \ref{sec:AlgorithmicDetails}. The alternative mode, referred to as ``CY'', prioritizes cycles over chains by restricting the number of chain columns returned to the RMP at every iteration of the column generation algorithm to up most one. 

\rev{We note that in CH mode, there is no guarantee that the optimal solution will use the largest number of patients matched through chains. Since the prioritization of chains is encouraged by adding positive-price chains over positive-price cycles, whenever a positive-price chain is encountered, thus, the optimal base of the RMP has more possibilities to build a KEP solution with as many chains as possible, but this solution (although optimal) may or may not match the largest number of patients through chains.} In CY mode, since the size of a K-limited FVS is generally much larger than the number of singleton donors (Online Supplement \ref{sec:OtherResults}), there is a high chance of discovering a larger number of positive-price cycle columns, as opposed to those corresponding to chains. Therefore, if possible, the resulting feasible solution consists mostly of cycles.

\rev{To compare the impact of prioritizing chains (CC) or cycles (CY), Table \ref{tab:SolStructure} shows the average chain and cycle length and number of cycles for the instances shown in Figure \ref{fig:ValueOfChains}. It is worth noting that \textit{all} singleton donors were used under both modes. Under CH mode, BP\_MDD finds longer chains and fewer cycles. In both cases, the length of cycles remains steady around three. The difference in the chain length becomes more evident as the instances become denser. The same overall behavior is observed in set KEP\_N (Online Supplement \ref{sec:OtherResults}).}

\rev{Due to the possible presence of multiple optima, the average length of chains should be analyzed in context. The prioritization of chains through the CH mode in Table \ref{tab:SolStructure} can be seen as a guidance on what the average length of chains could be if we were to implement a solution that mostly consists of chains. Conversely, when prioritizing cycles over chains, we obtain optimal solutions consisting of smaller chains and more cycles.}

In terms of computational performance, the solution time for KEP\_N instances under the CY mode, was on average 8 times slower when compared to CH. This result can be explained due to the restriction on chain columns in the RMP. Thus, forcing BP\_MDD to go over more iterations to prove optimality of the RMP. Moreover, about 90\% of the KEP\_N instances under the CY mode were solved to optimality. The average gap was 0.3\%.

\begin{table}[tbp]
    \spacing
        \caption{\rev{Solution composition by chain (CH) and cycle (CY) prioritization modes. Data is aggregated by the small and medium size instances in Figure \ref{fig:ValueOfChains} and cycle lengths $\Kcycle = \{3,4\}$.}} 
        \centering
        \label{tab:SolStructure}
        %\begin{adjustbox}{width=1\textwidth}
            \begin{tabular}{lccccccc}
  \toprule
   &  & \multicolumn{3}{c}{CH, CY mode average}            & \multicolumn{3}{c}{Average CH$/$CY ratio}          \\ 
    \cmidrule(lr){3-5} \cmidrule(lr){6-8} 
  \multicolumn{1}{c}{$(\sigma, p_\ell)$}          & \multicolumn{1}{c}{$\Lchain$}            & \multicolumn{1}{c}{Chain length} & \multicolumn{1}{c}{Cycle length} & \multicolumn{1}{c}{\# cycles} & \multicolumn{1}{c}{Chain length} & \multicolumn{1}{c}{Cycle length} &  \# cycles \\
    \cmidrule(lr){1-1}
    \cmidrule{2-2}
    \cmidrule(lr){3-3}
    \cmidrule(lr){4-4}
     \cmidrule(lr){5-5}
     \cmidrule(lr){6-6}
    \cmidrule(lr){7-7}
     \cmidrule(lr){8-8}
  (0.27,0.45)        & 3    & 2.93,\ 2.91 & 3.08,\ 3.09 & 46.25,\ 46.67        & 1.01     & 1.00 & 0.99           \\
  (0.27,0.45)       & 6    & 4.65,\ 4.45  & 2.93,\ 2.90 & 37.33,\ 40.58       & 1.06     & 1.01 & 0.81         \\
  (0.27,0.8)         & 3    & 2.92,\ 2.79  & 2.99,\ 2.78 & 50.58,\ 52.50        & 1.05     & 1.00 & 0.93            \\
  (0.27,0.8)         & 6    & 4.47,\ 4.25  & 2.88,\ 2.82  & 40.42,\ 44.75       & 1.08     & 1.02 & 0.79         \\
  (0.4,0.45)         & 3    & 2.93,\ 2.84 & 2.96,\ 2.95  & 54.75,\ 56.08       & 1.04     & 1.00 & 0.95           \\
  (0.4,0.45)         & 6    & 4.35,\ 4.05 & 2.80,\ 2.84 & 44.08,\ 48.33         & 1.09      & 0.98 & 0.83         \\
  (0.4,0.8)         & 3    & 2.89,\ 2.67  & 2.91,\ 2.86  & 57.58,\ 61.17       & 1.09      & 1.02 & 0.90          \\
  (0.4,0.8)         & 6    & 3.91,\ 3.36  & 2.84,\ 2.80 & 49.33,\ 56.67        & 1.19      & 1.01  & 0.78        \\
  (0.7,0.45)         & 3    & 2.85,\ 2.46 & 2.82,\ 2.78 & 64.08,\ 70.25        & 1.18      & 1.01 & 0.86         \\
  (0.7,0.45)        & 6    & 3.41,\ 2.73  & 2.79,\ 2.75 & 58.75,\ 68.00       & 1.26       & 1.02 & 0.78        \\
  (0.7,0.8)          & 3    & 2.85,\ 2.30  & 2.84,\ 2.82 & 64.33,\ 71.75       & 1.29      & 1.02 & 0.83        \\
  (0.7,0.8)          & 6    & 3.19,\ 2.54 & 2.83,\ 2.83  & 60.92,\ 69.17       & 1.30      & 1.00 & 0.80           \\ \bottomrule
            \end{tabular}
\end{table}

Note that  BP\_MDD can be extended to encompass other priorities.  It is known that highly-sensitized patients tend to wait longer for a match offer than those who are not. If exchanges of a given size are preferred,  the longest path for the MDDs can focus on finding, if any,  positive-price chains (cycles) among the targeted exchanges. Phases two and three also offer alternatives to choose which and how many of those exchanges to return to the RMP. These prioritization schemes work in a  similar fashion to having hierarchical objectives \citep{Glorie2014}, with the advantage that multiple runs of the algorithm are not required. 
}
\section{Conclusion}
\label{sec:Conclusions}

\rev{We addressed the problem of finding a matching in kidney exchange, considering long chains, and developed a novel \BP\ framework in which the pricing problems for both cycles and chains are primarily solved via MDDs. We note that our framework is, to the best of our knowledge, the only exact algorithm for KEP in which a preference for cycles or chains can be selected a priori by decision-makers and the first correct B\&P algorithm for the cycles-and-chains variant \citep{Omer2022}, as well proved a new Lagrangrian-based valid upper bound to the problem. This framework showed significant improved computational performance over current state-of-the-art methods, though differences in programming language (Python for \citet{Dickerson2016} and C++ for us) prevent a true apples-to-apples comparison despite using the same computational environment. While our \BP\ framework is presented specifically for KEP---as all \BP\ implementations are bespoke---it can be adapted to any problem of detecting cycles and chains in a network (e.g., variations of the travelling salesman problem), or where the pricing problem can be formulated as an MDD.}

\rev{A limitation of restricted decision diagrams is that they only include a subset
of the cycles or chains that may be needed in the restricted master problem, thus, the use of
the additional phases is still required for arbitrary cycle and chain lengths.}
\rev{For future work, we will develop a simulation to examine the performance of our approach over multiple matching epochs to assess long-term real-world utility, and consider metrics such as maximum time to transplant, with special attention on outcomes for highly-sensitized patients. We note consideration of time to transplant and other priority rules can be incorporated into the current framework as objective function weights, as the objective function is already weighted per patient.}

\rev{Finally, a significant barrier to implementation of mathematically-driven decision-making in healthcare applications is complexity of the optimization techniques. In some applications, a web-based tool or other UIUX can help adoption of mathematically complex decision-making tools, but in Canada, just like in other countries, the KPDPs operators are typically appointed to government-run institutions that have their own regulations on clinical data treatment and
other protocols that are unique to the specific needs of a KPDP, and thus, a web tool may not be highly valued from an implementation perspective. However, such a tool could be used to communicate the value of investing in high performance optimization software and personnel to administrative decision-makers; thus, development of a front-end tool is another direction of future work.}

% Acknowledgments here
\ACKNOWLEDGMENT{%
% Enter the text of acknowledgments here
}% Leave this (end of acknowledgment)

% Appendix here
% Options are (1) APPENDIX (with or without general title) or 
%             (2) APPENDICES (if it has more than one unrelated sections)
% Outcomment the appropriate case if necessary
%
% \begin{APPENDIX}{<Title of the Appendix>}
% \end{APPENDIX}
%
%   or 
%
% \begin{APPENDICES}
% \section{<Title of Section A>}
% \section{<Title of Section B>}
% etc
% \end{APPENDICES}

% References here (outcomment the appropriate case) 

% CASE 1: BiBTeX used to constantly update the references 
%   (while the paper is being written).
\bibliographystyle{informs2014} % outcomment this and next line in Case 1

{\spacing 
\bibliography{writeup_reference}}
%\bibliography{<your bib file(s)>} % if more than one, comma separated
\newpage
\begin{APPENDIX}{}
\section{Pricing Problems and New Valid Upper Bound}
\label{sec:Lagrange}
\rev{We introduce a Lagrangian decomposition of the KEP based on an improved and generalized version of an existing integer programming (IP) formulation. Particularly, we show that the Lagrangian decomposition leads to (i) an IP formulation for the pricing problems and (ii) a new valid upper bound on the optimal value of the KEP that can be used in \BP \, and is stronger than the other existing bound \citep{Abraham2007}. Lastly, we show an advantage of using the \textit{disaggregated cycle formulation} \citep{Klimentova2014} over the so-called \textit{cycle formulation} \citep{Abraham2007, Roth2007} in  \BP.} 

\subsection{IEEF: An Improved Extended Edge Formulation (EEF)}
\label{sec:IEEF}

For the cycle-only version of the KEP, \cite{Constantino2013}  {\cdionne cloned} the input digraph $\digraph$,  $\vert \pairset \vert$ times, drawn by the fact that $\vert \pairset \vert$ is a natural upper bound on the number of cycles in any {\cdionne feasible solution}. The idea is then to find a feasible cycle in each copy of $\digraph$. If the selected cycles across the copies  are vertex-disjoint, they represent a feasible matching. Thus, an IP formulation can be built based on this disaggregation of $\digraph$ to solve the problem. {\cdionne We  extend} this formulation by including long chains and reducing the number of copies. To this end, we start by showing that any \textit{feedback vertex set} (FVS) provides a valid upper bound on the number of vertex-disjoint cycles in $\digraph$, whose proof is omitted due to its simplicity. 

\begin{definition}
	Given a directed graph $\digraphc$, $\oldfvs \subseteq \Vertex$\ is a feedback vertex set of $\digraph$ if the subgraph induced by $\Vertex \setminus \oldfvs$ is acyclic.
\end{definition}	

\begin{proposition}
\label{prop:feedback}
	Given a digraph  $\digraphc$, let $\upperc$ be the cardinality of the set with the maximum number of vertex-disjoint cycles in $\digraph$, and $\oldfvs$ be an FVS. Then, $\upperc \le \vert \oldfvs \vert$.
\end{proposition}

%[OMITTED TO SAVE SPACE] \proof{Proof.} Let $\cmaxset$ be a set with the maximum number of vertex-disjoint cycles in $\digraph$, i.e., $\vert \cmaxset \vert = \upperc$. For each cycle $\cycle \in \cmaxset$, there exists at least one vertex in $\cycle$ also in $\fvs$, since by definition $\fvs$ is a FVS. Because the cycles in $\cmaxset$ are vertex disjoint, it follows that the cardinality of $\fvs$ is at least $\upperc$. \hfill$\square$

The goal is to create one cycle copy per every vertex in the FVS. We use Algorithm \ref{FVSAlg0} to this end, except that now, we replace the function \texttt{BuildMDD} in line 8 by a function that creates a digraph associated to the corresponding feedback vertex. Because our extension includes chains, we also create copies of the input graph for chains, one for every \single. We can now introduce our extension of the extended edge formulation, which will provide us with a mechanism to compute a valid upper bound on the maximum number of matches.

Let $\fbcopiesch$ denote the index set of chain copies, with $\vert \fbcopiesch \vert \le \vert \singleset \vert$. The $\fm$-th chain copy of $\digraph$, $\fm \in \fbcopiesch$, is represented by graph $\dichainf$, whose vertex set $\Vertexbar^{\fm} = \pairset \cup \{\uix_{\fm}, \tau \}$ is formed by all \pairsAb, the $\fm$-th \single, $\uix_{\fm}$, {\cdionne and a dummy vertex $\tau$, representing} a dummy patient receiving a fictitious donation from the paired donor in the  last pair of a chain. The set of arcs {\cdionne$\Arcsbar^{\fm} = ( \Arcs \setminus \{(\uix, \vix) \in \Arcs \mid  \uix \in \singleset \setminus \uix_{\fm},  \vix \in \pairset\})  \cup \{(\vix, \tau) \mid \vix \in \pairset\}$  removes arcs emanating from \singles\ other than $\uix_{\fm}$} and adds one dummy arc from every \pairAb\ to $\tau$. Thus, a chain can only be started by the \single\ $\uix_{\fm}$ on the $\fm$th chain copy. Additionally, let $\csetf$ and $\csetf_{\Kcycle} \subseteq \csetf$ be the set of all simple cycles  and the set of all feasible cycles in $\digraphbar^{\fm}$, respectively.  Lastly, let $\lagx_{\pairarc}^{\fm}$ and $\lagy_{\pairarc}^{\fm}$ be  {\cdionne decision variables} taking the value {\cdionne 1} if arc $\arcp \in \Arcsd^{\fm}$ and arc $\arcp \in \Arcsbar^{\fm}$ is selected in a cycle and chain, respectively, and {\cdionne 0} otherwise. Then, the {\cdionne improved extended edge formulation} {\cdionne is} formulated as follows:
%\objIP=\quad 
\begin{subequations}
	\label{sub:IEEF}
	\begin{align}
	 \max \quad & \sum_{\fm \in \fbcopiescy}\sum_{\arcp \in \Arcsd^{\fm}} w_{\pairarc}\lagx^{\fm}_{\pairarc} + \sum_{\fm \in \fbcopiesch}\sum_{\arcp \in \Arcsbar^{\fm}} w_{\pairarc}\lagy^{\fm}_{\pairarc} 
	\label{eq:obj} \tag{IEEF}\\
	\text{s.t.}\quad %\enskip &&&&&&&&&&&&
	&\sum_{\fm \in \fbcopiescy}\sum_{\arcp \in \Arcsd^{\fm}} \lagx^{\fm}_{\pairarc} + \sum_{\fm \in \fbcopiesch}\sum_{\arcp\in \Arcsbar^{\fm}} \lagy^{\fm}_{\pairarc}  \le 1  &&  \uix \in \Vertex
	\label{eq:maxflowIEEF}\\
	&\sum_{\vix: \arcp \in \Arcsd^{\fm}} \lagx^{\fm}_{\pairarc}  = \sum_{\vix: \inarcp\in \Arcsd^{\fm}} \lagx^{\fm}_{\inpairarc} && \fm \in \fbcopiescy, \  \uix \in \Vertexd^{\fm} 
	\label{eq:OutFlowCycleIEEF}\\
	&\sum_{\vix: \arcp \in \Arcsbar^{\fm}} \lagy^{\fm}_{\pairarc}  = \sum_{\vix: \inarcp \in \Arcsbar^{\fm}} \lagy^{\fm}_{\inpairarc} && \fm \in \fbcopiesch, \  \uix \in \Vertexbar^{\fm}  \setminus \{\uix_{\fm}, \tau\}
	\label{eq:OutFlowChainIEEF}\\
	&\sum_{\arcp \in \Arcsd^{\fm}} \lagx^{\fm}_{\pairarc}  \le \Kcycle &&  \fm \in \fbcopiescy 
	\label{eq:CycleSizeIEEF}\\
		&\sum_{\arcp \in \Arcsbar^{\fm}} \lagy^{\fm}_{\pairarc}  \le \Lchain &&  \fm \in \fbcopiesch
	\label{eq:ChainSizeIEEF}\\
	&\sum_{\vix:\arcp \in \Arcsd^{\fm}} \lagx^{\fm}_{\pairarc}  \le  \sum_{\vix:(\vix_{\fm}^{*}, \vix) \in \Arcsd^{\fm}} \lagx^{\fm}_{\vixstar_{\fm} \vix}  && \fm \in \fbcopiescy {\cred{;}} \  \uix \neq \vixstar_\fm
	\label{eq:symIEEF}\\
	&{\cred \sum_{\vix:\arcp \in \Arcsbar^{\fm}} \lagy^{\fm}_{\pairarc}  \le  \sum_{\vix:(u_{\fm}, \vix) \in \Arcsbar^{\fm}} \lagy^{\fm}_{u_{\fm} \vix}} && {\cred \fm \in \fbcopiesch; \  \uix \in \Vertexbar^{\fm}  \setminus \{\uix_{\fm}, \tau\}}
	\label{eq:symChainsIEEF}\\
	&\sum_{\arcp \in \Arcs(\cycle)} \lagy^{\fm}_{\pairarc}  \le  \vert \Vertex(\cycle) \vert - 1  && \fm \in \fbcopiesch, \  \cycle  \in \csetf \setminus \csetf_{\Kcycle}
	\label{eq:cycleElimIEEF}\\
	& \lagx^{\fm}_{\pairarc}  \in \{0,1 \}  && \fm \in \fbcopiescy, \ \arcp \in \Arcsd^{\fm}
	\label{eq:integercIEEF}\\
	& \lagy^{\fm}_{\pairarc}  \in \{0,1 \}  && \fm \in  \fbcopiesch, \ \arcp \in \Arcsbar^{\fm}
	\label{eq:integerpIEEF}
	\end{align}
\end{subequations}

The objective function maximizes the weighted sum of transplant scores. Constraints \eqref{eq:maxflowIEEF}  enforce that a donor (paired or single) donates at most one kidney. Constraints \eqref{eq:OutFlowCycleIEEF} guarantee that in a selected cycle copy, if a \pairAb\ receives a kidney, then it donates one to another pair in the same copy. Constraints \eqref{eq:OutFlowChainIEEF} ensure that in a chain if a patient in a \pairAb\ receives a kidney, his or her paired donor donates either to a \pairAb\ in the same chain or to a dummy vertex, as it is the case of the donor in the last \pairAb\ of the chain. {\cdionne Constraints \eqref{eq:CycleSizeIEEF} and \eqref{eq:ChainSizeIEEF}  limit the number of arcs in a cycle to $\Kcycle$ and the number of arcs in a chain to $\Lchain$ in a copy, respectively.} Constraints \eqref{eq:symIEEF} forbid the use of the $\fm$-th cycle copy unless the $\fm$-th vertex in the FVS, $\vixstar_{\fm}$, is selected in that copy. {\cred Likewise, constraints  \eqref{eq:symChainsIEEF} allow the use of the $\fm$-th chain copy only when the singleton donor, $u_{\fm}$, triggers a chain.} The presence of cycles (subtours) in chain copies, particularly of those with size higher than $\Kcycle$, jeopardize the correctness of the formulation. Therefore, constraints \eqref{eq:cycleElimIEEF} {\cdionne ensure} the elimination of all infeasible cycles from chain copies.  These constraints resemble those used in the recursive formulation of \cite{Anderson2015}{\cdionne, except that their arc variables are based on the input graph in a model allowing arbitrarily long chains.} {\cdionne We} note that \cite{Constantino2013} and \cite{Klimentova2014} did not consider infeasible-cycle-breaking constraints in their discussion on how to include \singles\ in the EEF, and thus can be deemed incomplete. Infeasible-cycle-breaking constraints, e.g., Constraints \eqref{eq:cycleElimIEEF}, are \textit{necessary} to preserve the correctness of the EEF, and thus that of IEEF. {\cred A counter example showing the need of such constraints in EEF to allow chains and a formal proof for Proposition \ref{prop:IEEF_Proof_FirstCopy} are provided in Section \ref{sec:Proofs}.} 
	
\begin{proposition}
	\label{prop:IEEF_Proof_FirstCopy}
{\cred IEEF is a correct formulation of the matching problem in the KEP.}
\end{proposition}
%Lastly, constraints \eqref{eq:integercIEEF} and \eqref{eq:integerpIEEF} indicate decision variables' domain.
%Due to the not complete removal of symmetric solutions, it is possible that in a feasible solution multiple cycles from the same cycle copy are selected, e.g., two different 2-way cycles can be selected from a cycle copy when $\Kcycle = 4$. However, since cycle copies are restricted to use fewer than $\Kcycle$ arcs, and exchanges are vertex-disjoint by constraints \eqref{eq:maxflowIEEF},\eqref{eq:OutFlowCycleIEEF} and \eqref{eq:CycleSizeIEEF}, such phenomenon does not compromise the correctness of the model in the absence of \singles. However, t
%Despite of this fact, not all symmetric solutions are removed, since in a solution, a cycle copy $ \digraphd^{\fm}$ can include one or more feasible cycles, where one cycle is covered by $\vix_{\fm}$ and the other ones by feedback vertices in copies with higher index. 

\subsection{Lagrangian relaxation}

Consider introducing the following valid (redundant) constraints to \eqref{eq:obj}:
 \begin{subequations}
 \label{subeqs:maxflowIEEF2}
 	\begin{align}
 	 &\sum_{\arcp \in \Arcsd^{\fm}}  \lagx^{\fm}_{\pairarc}   \le 1  &&  \fm \in \fbcopiescy, \uix \in \Vertexd^{\fm}
 	\label{eq:maxflowIEEF2cycles}\\
 	 & \sum_{\arcp \in \Arcsbar^{\fm}} \lagy^{\fm}_{\pairarc}  \le 1  && \fm \in \fbcopiesch, \uix \in \Vertexbar^{\fm}
 	\label{eq:maxflowIEEF2chains}
 	\end{align}
\end{subequations}

\noindent Before justifying these new constraints, let us first approximate the optimal objective value of  \eqref{eq:obj} by relaxing constraints guaranteeing at most one donation from every donor Constraints \eqref{eq:maxflowIEEF} and then penalizing their violation by imposing Lagrange multipliers $ (\boldsymbol{\multi})$ in the objective function. This relaxation may allow a vertex to be selected in more than one graph copy, thus, in more than one exchange. However, if a copy is selected, such a vertex can be selected at most once within that copy, {\cdionne due to Constraints \eqref{subeqs:maxflowIEEF2}}. Given $\multi \in \mathbb{R}_{+}^{\vert \Vertex \vert}$, a Lagrangian relaxation to the KEP can be formulated as follows:
\begin{subequations}
	\label{sub:LR1} 
	\begin{align}
	 \objLag (\multi) := \hspace*{-0.3cm} &&& \label{objLag}\tag{LR1}\\
	 \max \ & \sum_{\fm \in \fbcopiescy}\sum_{\arcp \in \Arcsd^{\fm}} w_{\pairarc}\lagx^{\fm}_{\pairarc} + \sum_{\fm \in \fbcopiesch}\sum_{\arcp \in \Arcsbar^{\fm}} w_{\pairarc}\lagy^{\fm}_{\pairarc}  + \sum_{\vix \in \Vertex} \multi_{\vix} \left(1 - \sum_{\fm \in \fbcopiescy}\sum_{\arcp \in \Arcsd^{\fm}} \lagx^{\fm}_{\pairarc} + \sum_{\fm \in \fbcopiesch}\sum_{\arcp \in \Arcsbar^{\fm}} \lagy^{\fm}_{\pairarc}  \right) 
	 \notag\\
	 \text{s.t.} \  %\enskip &&&&&&&&&&&&
	 &\eqref{eq:OutFlowCycleIEEF}  - \eqref{eq:integerpIEEF}, \eqref{eq:maxflowIEEF2cycles} - \eqref{eq:maxflowIEEF2chains}&&
	 \label{eq:AllOthersLag}
 	\end{align}
\end{subequations}

\noindent As the only constraints linking decision variables associated with different copies are now relaxed, \eqref{objLag} can be decomposed by graph copies as follows:
\begin{subequations}
	\label{sub:LR2} 
	\begin{align}
	\centering
		&\objLag (\multi) = \sum_{\fm \in \fbcopiesch} \objLagCy^{\fm}({\multi}) + \sum_{\fm \in \fbcopiescy} \objLagCh^{\fm} ({\multi})+ \sum_{\uix \in \Vertex} \multi_{\uix} \label{objLag2}\tag{LR2}
	\end{align}
\end{subequations}
{\cdionne where, $ \forall \fm \in \fbcopiescy $ and $\forall \fm \in \fbcopiesch$, we have the following subproblems, respectively:}
\begin{align}
	&\objLagCy^{\fm}({\multi}) := \max \left \{ \sum_{\arcp \in \Arcsd^{\fm}} (w_{\pairarc} - \multi_{\uix}) \lagx^{\fm}_{\pairarc} \mid  \eqref{eq:OutFlowCycleIEEF} , \eqref{eq:CycleSizeIEEF} \eqref{eq:symIEEF},\eqref{eq:integercIEEF}, \eqref{eq:maxflowIEEF2cycles} \right \}, 
	\label{sprbm:cycles}\tag{CC}
\end{align}
\begin{align}
	&\objLagCh^{\fm}({\multi}) := \max \left \{ \sum_{\arcp \in \Arcsbar^{\fm}} (w_{\pairarc} - \multi_{\uix}) \lagy^{\fm}_{\pairarc} \mid  \eqref{eq:OutFlowChainIEEF}, \eqref{eq:ChainSizeIEEF}, \eqref{eq:cycleElimIEEF}, \eqref{eq:integerpIEEF}, \eqref{eq:maxflowIEEF2chains} \right \}. 
	\label{sprbm:chains}\tag{CH}
\end{align}

Given a set of Lagrange multipliers, $\multi$, each subproblem {\cdionne finds} either a feasible cycle, \eqref{sprbm:cycles}, or a feasible chain, \eqref{sprbm:chains}, whose vertex assignment has minimum penalty, thus, maximum weight. Observe that the inclusion of the dummy vertex $\tau$ is useful to capture the Lagrange multiplier of the last pair in a chain in the objective function of \eqref{sprbm:chains}. The Lagrangian dual problem can be defined as $\Lagdual := \min \{ \objLag (\multi) : \multi \in \mathbb{R}^{\vert \Vertex\vert}_{+} \}$. That is, $\Lagdual$ is the smallest upper bound that can be obtained when the set of Lagrange multipliers favor an assignment of vertices to cycle and chain copies with minimum intersection. If we define $\objIP^{LP}$  as the optimal objective value of the linear programming relaxation of \eqref{eq:obj}, it is possible that $\Lagdual < \objIP^{LP}$. Moreover, we show that the quality of the bound provided by $\Lagdual$ is as tight as the one provided by the linear programming relaxation of the disaggregated cycle formulation \citep{Klimentova2014}, one of the formulations in the literature providing the tightest linear relaxation.

\begin{proposition}
	If $\objIP_{\cycle}^{LP}$ is the optimal objective value of the {\cred linear programming relaxation} of the disaggregated cycle formulation, then $\Lagdual = \objIP_{\cycle}^{LP}$.
\end{proposition}	

\proof{Proof.} Let $\cfeaset$ and $\pfeaset$ be the set of feasible cycles (including vertex $\vixstar_{\fm}$) and chains on the $\fm$th graph copy, $\fm \in \fbcopiescy, \fm \in \fbcopiesch$, respectively. For a  cycle $\cycle \in \cfeaset$ and chain $\chain \in \pfeaset$, let $\cweight = \sum_{\arcp \in \Arcs(\cycle)} w_{\pairarc}$ and $\pweight = \sum_{\arcp \in \Arcs(\chain)} w_{\pairarc}$ be the total weight of a cycle and chain, respectively. Then, \eqref{objLag2}  is reformulated as follows: 
\begin{subequations}
	\label{sub:LR2opt} 
	\begin{align}
	\min \quad & \sum_{\fm \in \fbcopiesch} \objLagCy^{\fm} + \sum_{\fm \in \fbcopiescy} \objLagCh^{\fm} + \sum_{\vix \in \Vertex} \multi_{\vix} \label{objLag2opt}\tag{LR3}\\
	&\objLagCy^{\fm} \ge \cweight  - \sum_{\vix \in \Vertex(\cycle)}\multi_{\vix} &&\fm \in \fbcopiescy, \cycle \in \cfeaset  &  (\zcycleb)&  \label{eq:Lagoptcy}\\
	&\objLagCh^{\fm} \ge \pweight  - \sum_{\vix \in \Vertex(\chain)} \multi_{\vix} &&\fm \in \fbcopiesch, \chain \in \pfeaset & (\zchainb)&  \label{eq:Lagoptch}\\
	&  \objLagCy^{\fm}  \ge 0 &&\fm \in \fbcopiescy  \label{eq:bound1}\\
	&  \objLagCh^{\fm}  \ge 0 &&\fm \in \fbcopiesch  \label{eq:bound1.1}\\
	&  \multi_{\vix} \ge 0 &&  \vix \in \Vertex \label{eq:bound2}
	 \end{align}
\end{subequations}

\eqref{objLag2opt} finds the optimal value of decision variables $\objLagCy^{\fm}$, $\objLagCh^{\fm} $ and $\multi$. The validity of (\ref{objLag2opt}) relies on the fact that the maximum weight cycle and chain is selected for every graph copy through Constraints \eqref{eq:Lagoptcy} and \eqref{eq:Lagoptch}, met in the equality by the minimization objective. Moreover, since the objective value of \eqref{sprbm:cycles} and \eqref{sprbm:chains}  can at least be zero (by selecting $\lagx = 0$ and $\lagy = 0$), Constraints \eqref{eq:bound1} and \eqref{eq:bound1.1} represent a valid lower bound on the objective value of each sub-problem. To finalize the proof, {\cdionne letting} $\zcycleb$ and $\zchainb$ be the dual variables of constraints \eqref{eq:Lagoptcy} and \eqref{eq:Lagoptch}, respectively, the dual problem of \eqref{objLag2opt}  is
\begin{subequations}
	\label{sub:discyclef} 
	\begin{align}
	\max \quad&  \sum_{\fm \in \fbcopiescy}  \sum_{\cycle \in \cfeaset} \cweight\zcycleb + \sum_{\fm \in \fbcopiesch}  \sum_{\chain \in \pfeaset} \pweight\zchainb  \label{objDis}\tag{IDCF}\\
	&\sum_{\fm \in \fbcopiescy} \sum_{\cycle \in \cfeaset: \vix \in \Vertex(\cycle)} \zcycleb + \sum_{\fm \in \fbcopiesch} \sum_{\chain \in \pfeaset: \vix \in \Vertex(\chain)}  \zchainb \le 1 && \vix \in \Vertex \label{eq:onepervertex}\\
	&  \zcycleb \ge 0 && \fm \in \fbcopiescy, \cycle \in \cfeaset\\
	&  \zchainb \ge 0 && \fm \in \fbcopiesch, \chain \in \pfeaset \label{eq:BinaryChain}
	 \end{align}
\end{subequations}

Note that Constraints $ \sum_{\cycle \in \cfeaset} \zcycleb \le 1 \text{\ \ } \forall \fm \in \fbcopiescy $ are omitted from \eqref{objDis} along with their chain counterpart  since they are implied by Constraints \eqref{eq:onepervertex}. { \cred Formulation \eqref{objDis} corresponds to the integer linear programming relaxation of the \textit{disaggregated cycle formulation} in \cite{Klimentova2014}, noting that in \eqref{objDis} chain variables are included and cycle copies are found through Algorithm \ref{FVSAlg0}}. Hence, by strong duality it follows that $\Lagdual = \objIP_{\cycle}^{LP}$.  \hfill$\square$

{ \cred Notice that by enforcing integrality  of the decision variables, \eqref{objDis} becomes a valid formulation for the matching problem in kidney exchange, similar to the so-called cycle formulation \citep{Abraham2007, Roth2007}, except that cycles and chains are not found in specific graph copies.} Thus, the $\fm$ index is dropped. It had not been shown before  whether there is an advantage {\cdionne to} using one formulation over another, particularly {\cdionne for} \BP. Note that the dual variables of Constraints \eqref{eq:onepervertex} correspond to the Lagrange multipliers $\multi$ in \eqref{objLag2}. \rev{Thus, the subproblems in the Lagrangian relaxation, \eqref{sprbm:cycles} and \eqref{sprbm:chains}, correspond to the pricing problems for the RMP based on \ref{objDis}.} Moreover, this result shows that when cycles and chains are disaggregated into graph copies, it is possible to obtain a valid upper bound on the objective value by solving \eqref{sprbm:cycles} and \eqref{sprbm:chains}, and simply plugging their result into \eqref{objLag2} afterwards. {\cdionne Even} if the set of Lagrange multipliers is not optimal, \eqref{objLag2} provides a valid upper bound, which can be as good as that of the disaggregated cycle formulation or the cycle formulation itself. Thus, even without  proving optimality of \eqref{objDis}, its dual variables can still be used to obtain a valid upper bound. To the best of our knowledge, the only previous method in the literature obtaining a valid upper bound consists of solving a relaxed problem with $\Kcycle = \Lchain = \infty$ \citep{Abraham2007}, which can be solved in polynomial time. It is easy to see that the bound provided by this special case is weaker than that of the presented Lagrangian dual problem. Since \eqref{objDis} can be used as a master problem in column generation, the goal is to use the bound provided by \eqref{objLag2} when it is not possible to prove optimality of the master problem. In Section \ref{sub:branching}, we {\cdionne show} how this new upper bound can be used not just at every node of a branch-and-bound tree.

%To the best of our knowledge this is the first bounding proof directly applicable to Branch and Price algorithms for the KEP, that is, the quality of a feasible solution can be quantified when proving optimality of the master problem is not possible within time or computer-memory limits.

%Although, one could first find a solution through the cycle formulation and then group the cycles and paths following the vertices in the FVS to obtain a valid upper bound, disaggregation allows to target the search effort on specific subproblems that can be solved separately, making \eqref{objDis} potentially more useful than the cycle formulation for large instances.

\section{On the completeness of two existing formulations}
\label{sec:Appendix}
In this section, we present the edge assignment formulation and the extended edge formulation including \singles, as proposed by \cite{Constantino2013}. We show that for these formulations to model chains correctly, the inclusion of infeasible-cycle-breaking constraints is required. 

\subsection{Edge Assignment Formulation}
\label{sec:EAF}

Following an equivalent notation to that used by \cite{Constantino2013}, we proceed to introduce some notation. Let $D = (V, A)$ be a digraph representing compatibilities among donors (single or paired) and \pairsAb. The set of vertices $V = \{1,..., \vert \pairset \vert + \vert \singleset \vert \}$ has $\vert \pairset \vert$-many \pairsAb\ and $\vert \singleset \vert$-many \singles. Let vertices $\{1,...,\vert \singleset \vert\}$ represent \singles\ and vertices $\{\vert \singleset \vert + 1,..., \vert \singleset \vert + \vert \pairset \vert\}$ represent \pairsAb. An arc $(i,j) \in A$ exists if the donor in vertex $i$ is compatible with patient in vertex $j$. Assume that a dummy patient is associated to each \single, so that paired donors $j \in \{\vert \singleset \vert + 1, ..., \vert \singleset \vert + \vert \pairset \vert\}$ are compatible with each dummy patient $i \in \{1,..., \vert \singleset\vert \}$.  Also, consider $\vert V \vert$ as an upper bound on the number of cycles and chains in any feasible solution. For each vertex $\ell \in \{\vert \singleset \vert + 1,...,\vert \pairset \vert + \vert \singleset \vert\}$, let $V^{\ell} = \{i \in V \mid i \ge \ell \}$ be the set of vertices forming cycles with index higher or equal to $\ell$, whereas for each index $\ell \in \{1,..., \vert \singleset \vert\}$ let $V^{\ell} = \{i \in \pairset\} \cup \{\ell\}$ be the set of vertices forming part of a chain started by the $\ell$-th \single. Notice that only vertices $i \ge \ell$ are included in each vertex set to remove multiplicity of solutions. Moreover, it can happen that for some $\ell \in \{1,..., \vert V \vert \}$, $V^{\ell}= \emptyset$, e.g., if all vertices pointing or receiving an arc from the vertex with the lowest index are removed. Thus, denote by $\mathscr{L}$ the set of indices for which $V^{\ell} \ne \emptyset$. Lastly, consider the following binary decision variables: $x_{ij}$ takes on value 1 if arc $(i,j)$ is selected in a cycle or chain and $y_{i}^{\ell}$ takes on value 1 if node $i$ is assigned to the $\ell$th cycle or chain. Then, the edge assignment formulation is defined as follows:
%

%$$
%x_{ij} = \left\{
%\begin{array}{ll}
%1,      & \mbox{if arc } (i,j) \mbox{ is selected in a cycle or chain} \\
%0,      & \mbox{otherwise}
%\end{array}
%\right.
%$$

%$$
%y_{i}^{\ell} = \left\{
%\begin{array}{ll}
%1,      & \mbox{if node } i \mbox{ is assigned to the } \ell\mbox{-th}\mbox{ cycle (chain)} \\
%0,      & \mbox{otherwise}
%\end{array}
%\right.
%$$
\begin{subequations}
	\label{sub:EAF}
	\begin{align}
	\max &\sum_{(i,j) \in A} w_{ij}x_{ij}
	\label{eq:ObjEAF}\\
	&\sum_{j:(j,i) \in A} x_{ji} =  \sum_{j:(i,j) \in A} x_{ij} && i \in V
	\label{eq:Balance}\\
	&\sum_{j:(i,j) \in A} x_{ij} \le 1&& i \in V
	\label{eq:MaxOutput}\\
	&\sum_{i \in \{\ell\} \cup \{\vert \singleset \vert + 1, ..., \vert \pairset \vert + \vert \singleset \vert \}} y_{i}^{\ell} \le \Lchain && \ell \in \{1,..., \vert \singleset \vert\}
	\label{eq:ChainSize}\\
	&\sum_{i \ge  \vert \singleset \vert + 1} y_{i}^{\ell} \le \Kcycle && \ell \in \{\vert \singleset \vert + 1,..., \vert \singleset \vert + \vert \pairset \vert\}
	\label{eq:CycleSize}\\
	&\sum_{\ell \in \mathscr{L}: i \in V^{\ell}} y^{\ell}_{i} = \sum_{j: (i,j) \in A} x_{ij} && i \in V
	\label{eq:activate}\\
	&y_{i}^{\ell} + x_{ij} \le 1 + y_{j}^{\ell}&& (i,j) \in A, \ell \in \mathscr{L}, i \in V^{\ell}
	\label{eq:assignment}\\
	&y^{\ell}_{i} \le y^{\ell}_{\ell}&& \ell \in \mathscr{L}, i \in V^{\ell}
	\label{eq:vlvertex}\\
	&y^{\ell}_{i} \in \{0,1\} && \ell \in \mathscr{L}, i \in V^{\ell}
	\label{eq:integrality1}\\
	&x_{ij} \in \{0,1\}&& (i,j) \in A
	\label{eq:integrality2}
	\end{align}
\end{subequations}

Constraints \eqref{eq:Balance} assure that patient $i$ receives a kidney if and only if donor $j$ donates a kidney. Constraints \eqref{eq:MaxOutput} allow at most one donation. Constraints \eqref{eq:ChainSize} and \eqref{eq:CycleSize} limit the length of chains and cycles. Constraints \eqref{eq:activate} ensure that vertex $i$ is in a cycle (chain) if and only if there is an assignment of $i$ to some $\ell$. Constraints \eqref{eq:assignment} state that if vertex $i$ is in cycle (chain) $\ell$ and donor $i$ gives a kidney to recipient $j$, then vertex $j$ must also be in the  $\ell$-th cycle (chain). Constraints \eqref{eq:vlvertex} establish that a vertex $i \in V^{\ell}$ can be assigned to the $\ell$-th cycle (chain) only if vertex $\ell$ is also assigned. Constraints \eqref{eq:integrality1} and \eqref{eq:integrality2} indicate decision variables' domain.

Now, we proceed to show a solution example satisfying \eqref{sub:EAF} and yet infeasible to the KEP. Consider Figure \ref{fig:Copies2} and assume $\Kcycle = 3$ and $\Lchain = 6$. Notice that in the solution, $y^{2}_{2} = y^{2}_{3} = y^{2}_{4} = y^{2}_{5} = y^{2}_{6} = y^{3}_{3} = y^{3}_{4} = y^{3}_{5} =y^{3}_{6} = y^{4}_{4} = y^{4}_{5} = y^{4}_{6} =  y^{5}_{5}  = y^{5}_{6} = x_{46} = x_{54} = x_{56} = x_{65} = x_{52} = x_{62} = x_{26} =  0$ and $y^{1}_{1} = y^{1}_{2} = y^{1}_{3}  = y^{1}_{4} = y^{1}_{5} = y^{1}_{6} = x_{53} = x_{36} = x_{64} = x_{45} = x_{12} = x_{21} = 1$. 

\begin{figure}[ht]
	\vspace*{0.25cm}
	\begin{subfigure}{.20\textwidth}
		\centering
		% include first image
		\caption{$D = (V,A)$}
		\centering
		\tikzstyle{place}=[circle,draw=blue!50,fill=blue!20,thick,
		inner sep=0pt,minimum size=4mm]
		\tikzstyle{transition}=[circle,draw=black!50,fill=black!20,thick,
		inner sep=0pt,minimum size=4mm]
		\tikzstyle{altru}=[rectangle,draw=black!50,fill=black!20,thick,
		inner sep=0pt,minimum size=4mm]
		\begin{tikzpicture}
		\node[transition] (waiting) at (0,2) {$5$}; 
		\node[transition] (critical) at (0,1) {$4$};
		\node[transition] (semaphore) at (0,0) {$6$};
		\node[altru] (alt) at (2,1)  {$1$}; 
		\node[transition] (leave critical) at (1,1)  {$2$}; 
		\node[transition] (enter critical) at (-1,1) {$3$};
		%\draw [->] (enter critical) to (critical);
		\draw [dashed, ->] (waiting) to [bend left=45] (alt);
		\draw [dashed, ->] (leave critical)  [bend right=45] to (alt);
		\draw [dashed, ->] (semaphore)  [bend right=45] to (alt);
		\draw [dashed, ->] (critical)  [bend left=45] to (alt);
		\draw [dashed, ->] (enter critical)  [bend left=45] to (alt);
		\draw [-{Classical TikZ Rightarrow}] (alt) to [bend right=45] (leave critical);
		\draw [-{Classical TikZ Rightarrow}] (waiting) to [bend left=45] (leave critical);
		\draw [-{Classical TikZ Rightarrow}] (waiting) to [bend right=45] (enter critical);
		\draw [-{Classical TikZ Rightarrow}] (enter critical) to [bend right=45] (semaphore);
		\draw [-{Classical TikZ Rightarrow}] (semaphore) to [bend right=45] (waiting);
		\draw [-{Classical TikZ Rightarrow}] (waiting) to [bend right=45] (semaphore);
		\draw [-{Classical TikZ Rightarrow}] (semaphore) to [bend right=45] (leave critical);
		\draw [-{Classical TikZ Rightarrow}] (waiting) to [bend right=45] (critical);
		\draw [-{Classical TikZ Rightarrow}] (critical) to [bend right=45] (waiting);
		\draw [-{Classical TikZ Rightarrow}] (semaphore) to [bend right=45] (critical);
		\draw [-{Classical TikZ Rightarrow}] (critical) to [bend right=45] (semaphore);
		\draw [-{Classical TikZ Rightarrow}] (leave critical) to  (semaphore);
		\end{tikzpicture}
		\label{fig:PCompleteGraph2}
	\end{subfigure}
	\begin{subfigure}{.20\textwidth}
		\centering
		% include second image
		\caption{$\ell = 1$. Chain 1}
		\centering
		\tikzstyle{place}=[circle,draw=blue!50,fill=blue!20,thick,
		inner sep=0pt,minimum size=4mm]
		\tikzstyle{transition}=[circle,draw=black!50,fill=black!20,thick,
		inner sep=0pt,minimum size=4mm]
		\tikzstyle{altru}=[rectangle,draw=blue!50,fill=blue!20,thick,
		inner sep=0pt,minimum size=4mm]
		\begin{tikzpicture}
		\node[transition] (waiting) at (0,2) {$5$}; 
		\node[transition] (critical) at (0,1) {$4$};
		\node[transition] (semaphore) at (0,0) {$6$};
		\node[altru] (alt) at (2,1)  {$1$}; 
		\node[transition] (leave critical) at (1,1)  {$2$}; 
		\node[transition] (enter critical) at (-1,1) {$3$};
		%\draw [->] (enter critical) to (critical);
		\draw [-{Classical TikZ Rightarrow}, ultra thick] (alt) to [bend right=45] (leave critical);
		\draw [-{Classical TikZ Rightarrow}, ultra thick] (leave critical)  to [bend right=45] (alt);
		\draw [-{Classical TikZ Rightarrow}] (waiting) to [bend left=45] (leave critical);
		\draw [-{Classical TikZ Rightarrow}, ultra thick] (waiting) to [bend right=45] (enter critical);
		\draw [-{Classical TikZ Rightarrow}, ultra thick] (enter critical) to [bend right=45] (semaphore);
		\draw [-{Classical TikZ Rightarrow}] (semaphore) to [bend right=45] (waiting);
		\draw [-{Classical TikZ Rightarrow}] (waiting) to [bend right=45] (semaphore);
		\draw [-{Classical TikZ Rightarrow}] (semaphore) to [bend right=45] (leave critical);
		\draw [-{Classical TikZ Rightarrow}] (waiting) to [bend right=45] (critical);
		\draw [-{Classical TikZ Rightarrow}, ultra thick] (critical) to [bend right=45] (waiting);
		\draw [-{Classical TikZ Rightarrow}, ultra thick] (semaphore) to [bend right=45] (critical);
		\draw [-{Classical TikZ Rightarrow}] (critical) to [bend right=45] (semaphore);
		\draw [-{Classical TikZ Rightarrow}] (leave critical) to  (semaphore);
		\end{tikzpicture}
		\label{fig:FirstCopy2}
	\end{subfigure}
	\begin{subfigure}{.20\textwidth}
		\centering
		% include second image
		\caption{$\ell = 2$. Cycle 1}
		\centering
		\tikzstyle{place}=[circle,draw=blue!50,fill=blue!20,thick,
		inner sep=0pt,minimum size=4mm]
		\tikzstyle{transition}=[circle,draw=black!50,fill=black!20,thick,
		inner sep=0pt,minimum size=4mm]
		\tikzstyle{blank}=[circle,draw=white,fill=white,thick,
		inner sep=0pt,minimum size=4mm]
		\begin{tikzpicture}
		\node[blank] (alt) at (2,1)  {}; 
		\node[transition] (waiting) at (0,2) {$5$}; 
		\node[transition] (critical) at (0,1) {$4$};
		\node[transition] (semaphore) at (0,0) {$6$};
		\node[place] (leave critical) at (1,1)  {$2$}; 
		\node[transition] (enter critical) at (-1,1) {$3$};
		%\draw [->] (enter critical) to (critical);
		\draw [-{Classical TikZ Rightarrow}] (waiting) to [bend left=45] (leave critical);
		\draw [-{Classical TikZ Rightarrow}] (waiting) to [bend right=45] (enter critical);
		\draw [-{Classical TikZ Rightarrow}] (enter critical) to [bend right=45] (semaphore);
		\draw [-{Classical TikZ Rightarrow}] (semaphore) to [bend right=45] (waiting);
		\draw [-{Classical TikZ Rightarrow}] (waiting) to [bend right=45] (semaphore);
		\draw [-{Classical TikZ Rightarrow}] (semaphore) to [bend right=45] (leave critical);
		\draw [-{Classical TikZ Rightarrow}] (waiting) to [bend right=45] (critical);
		\draw [-{Classical TikZ Rightarrow}] (critical) to [bend right=45] (waiting);
		\draw [-{Classical TikZ Rightarrow}] (semaphore) to [bend right=45] (critical);
		\draw [-{Classical TikZ Rightarrow}] (critical) to [bend right=45] (semaphore);
		\draw [-{Classical TikZ Rightarrow}] (leave critical) to  (semaphore);
		\end{tikzpicture}
		\label{fig:FirstCopy3}
	\end{subfigure}
	\begin{subfigure}{.20\textwidth}
		\centering
		% include second image
		\caption{$\ell = 3$. Cycle 2}
		\tikzstyle{place}=[circle,draw=blue!50,fill=blue!20,thick,
		inner sep=0pt,minimum size=4mm]
		\tikzstyle{transition}=[circle,draw=black!50,fill=black!20,thick,
		inner sep=0pt,minimum size=4mm]
		\tikzstyle{blank}=[circle,draw=white,fill=white,thick,
		inner sep=0pt,minimum size=4mm]
		\begin{tikzpicture}
		\node[blank] (alt) at (2,1)  {}; 
		\node[transition] (waiting) at (0,2) {$5$}; 
		\node[transition] (critical) at (0,1) {$4$};
		\node[transition] (semaphore) at (0,0) {$6$};
		\node[place] (enter critical) at (-1,1) {$3$};
		%\draw [->] (enter critical) to (critical);
		\draw [-{Classical TikZ Rightarrow}] (waiting) to [bend right=45] (enter critical);
		\draw [-{Classical TikZ Rightarrow}] (enter critical) to [bend right=45] (semaphore);
		\draw [-{Classical TikZ Rightarrow}] (semaphore) to [bend right=45] (waiting);
		\draw [-{Classical TikZ Rightarrow}] (waiting) to [bend right=45] (semaphore);
		\draw [-{Classical TikZ Rightarrow}] (waiting) to [bend right=45] (critical);
		\draw [-{Classical TikZ Rightarrow}] (critical) to [bend right=45] (waiting);
		\draw [-{Classical TikZ Rightarrow}] (semaphore) to [bend right=45] (critical);
		\draw [-{Classical TikZ Rightarrow}] (critical) to [bend right=45] (semaphore);
		\end{tikzpicture}
		\label{fig:FirstCopy4}
	\end{subfigure}
	%\bigskip
	\begin{subfigure}{0.3\textwidth}
		\centering
		% include second image
		\caption{$\ell = 4$. Cycle 3}
		\centering
		\tikzstyle{place}=[circle,draw=blue!50,fill=blue!20,thick,
		inner sep=0pt,minimum size=4mm]
		\tikzstyle{transition}=[circle,draw=black!50,fill=black!20,thick,
		inner sep=0pt,minimum size=4mm]
		\begin{tikzpicture}
		\node[transition] (waiting) at (0,2) {$5$}; 
		\node[place] (critical) at (0,1) {$4$};
		\node[transition] (semaphore) at (0,0) {$6$};
		%\draw [->] (enter critical) to (critical);
		\draw [-{Classical TikZ Rightarrow}] (semaphore) to [bend right=45] (waiting);
		\draw [-{Classical TikZ Rightarrow}] (waiting) to [bend right=45] (semaphore);
		\draw [-{Classical TikZ Rightarrow}] (waiting) to [bend right=45] (critical);
		\draw [-{Classical TikZ Rightarrow}] (critical) to [bend right=45] (waiting);
		\draw [-{Classical TikZ Rightarrow}] (semaphore) to [bend right=45] (critical);
		\draw [-{Classical TikZ Rightarrow}] (critical) to [bend right=45] (semaphore);
		\end{tikzpicture}
	\label{fig:FirstCopy5}
	\end{subfigure}
\begin{subfigure}{0.3\textwidth}
	\centering
	% include second image
	\caption{$\ell = 5$. Cycle 4}
	\centering
	\tikzstyle{place}=[circle,draw=blue!50,fill=blue!20,thick,
	inner sep=0pt,minimum size=4mm]
	\tikzstyle{transition}=[circle,draw=black!50,fill=black!20,thick,
	inner sep=0pt,minimum size=4mm]
	\begin{tikzpicture}
	\node[place] (waiting) at (0,2) {$5$}; 
	\node[transition] (semaphore) at (0,0) {$6$};
	%\draw [->] (enter critical) to (critical);
	\draw [-{Classical TikZ Rightarrow}] (semaphore) to [bend right=45] (waiting);
	\draw [-{Classical TikZ Rightarrow}] (waiting) to [bend right=45] (semaphore);
	\end{tikzpicture}
	\label{fig:FirstCopy6}
\end{subfigure}
	\begin{subfigure}{0.3\textwidth}
	\centering
	% include second image
	\caption{Complete solution}
	\centering
	\tikzstyle{place}=[circle,draw=blue!50,fill=blue!20,thick,
	inner sep=0pt,minimum size=4mm]
	\tikzstyle{transition}=[circle,draw=black!50,fill=black!20,thick,
	inner sep=0pt,minimum size=4mm]
	\tikzstyle{altru}=[rectangle,draw=blue!50,fill=blue!20,thick,
	inner sep=0pt,minimum size=4mm]
	\begin{tikzpicture}
	\node[transition] (waiting) at (0,2) {$5$}; 
	\node[transition] (critical) at (0,1) {$4$};
	\node[transition] (semaphore) at (0,0) {$6$};
	\node[altru] (alt) at (2,1)  {$1$}; 
	\node[transition] (leave critical) at (1,1)  {$2$}; 
	\node[transition] (enter critical) at (-1,1) {$3$};
	%\draw [->] (enter critical) to (critical);
	\draw [-{Classical TikZ Rightarrow}, ultra thick] (alt) to  [bend right=45] (leave critical);
	\draw [-{Classical TikZ Rightarrow}, ultra thick] (leave critical) to  [bend right=45] (alt);
	\draw [-{Classical TikZ Rightarrow}, ultra thick] (waiting) to [bend right=45] (enter critical);
	\draw [-{Classical TikZ Rightarrow}, ultra thick] (enter critical) to [bend right=45] (semaphore);
	\draw [-{Classical TikZ Rightarrow}, ultra thick] (critical) to [bend right=45] (waiting);
	\draw [-{Classical TikZ Rightarrow}, ultra thick] (semaphore) to [bend right=45] (critical);
	\end{tikzpicture}
	\label{fig:FirstCopy7}
\end{subfigure}
	\caption{Counter example where \eqref{sub:EAF} and \eqref{sub:EEF} provide an infeasible solution. An altruistic donor is represented by a square and blue vertices correspond to the vertex with lowest index in the $\ell$-th cycle or chain. Dashed arcs indicate compatibility of paired donors with a dummy patient associated to the altruistic donor, solid arcs indicate compatibilities among real donors and patients, and bold arcs represent the ones selected in a  solution.}
	\label{fig:Copies2}
\end{figure}

\subsection{EEF: Extended Edge Formulation}

Consider $\pairset + \singleset$ copies of the graph $D$, $D^{\ell} = (V^{\ell}, A^{\ell})$, where $V^{\ell}$ is as defined in Section \ref{sec:EAF} and $A^{\ell} = \{(i,j) \in \Arcs \mid i,j \in V^{\ell} \}$ is the set of arcs in the $\ell$-th copy. Recall that $\pairset$ and $\singleset$ is an upper bound on the number of cycles and chains in a feasible solution, respectively. In each copy $\ell \in \{1,..., \vert \singleset \vert\}$ chains are triggered by the $\ell$-th \single\ and at most $\Lchain$ arcs can be selected. In each copy $\ell \in \{\vert \singleset \vert + 1,..., \vert \singleset \vert + \vert \pairset \vert\}$, all cycles include the vertex with lowest index in that copy (e.g., see Figure \ref{fig:Copies2}). Let $x_{ij}^{\ell}$ be an arc variable taking the value one if the arc $(i,j) \in A^{\ell}$ is selected in the $\ell$-th copy and zero, otherwise. Similarly to EAF, consider $\mathscr{L}$ as the set of indices for which $V^{\ell} \ne \emptyset$.  Then, the extended edge formulation can be formulated as follows:
\begin{subequations}
	\label{sub:EEF}
	\begin{align}
	\max &\sum_{\ell \in \mathscr{L}}\sum_{(i,j) \in A^{\ell}} w_{ij}x_{ij}^{\ell}
	\label{eq:ObjEEF}\\
	&\sum_{j:(j,i) \in A^{\ell}} x_{ji}^{\ell} =  \sum_{j:(j,i) \in A^{\ell}} x_{ij}^{\ell} && \ell \in \mathscr{L}, i \in V^{\ell}
	\label{eq:EBalance}\\
	&\sum_{\ell \in \mathscr{L}} \sum_{j:(i,j) \in A^{\ell}} x_{ij}^{\ell} \le 1&& i \in V
	\label{eq:EMaxOutput}\\
	&\sum_{(i,j) \in A^{\ell}} x_{ij}^{\ell}  \le \Lchain && \ell \in \{1,..., \vert \singleset \vert\}
	\label{eq:EChainSize}\\
	&\sum_{(i,j) \in A^{\ell}} x_{ij}^{\ell}  \le \Kcycle && \ell \in \{\vert \singleset \vert + 1,..., \vert \pairset \vert + \vert \singleset \vert\}
	\label{eq:ECycleSize}\\
	&\sum_{j:(i,j) \in A^{\ell}} x_{ij}^{\ell}  \le \sum_{j:(\ell,j) \in A^{\ell}} x_{\ell j}^{\ell} && \ell \in \mathscr{L}, i \in V^{\ell}
	\label{eq:ESymmetry}\\
	&x_{ij}^{\ell} \in \{0,1\}&& \ell \in \mathscr{L}, (i,j) \in A^{\ell}
	\label{eq:Eintegrality}
	\end{align}
\end{subequations}

Constraints \eqref{eq:EBalance} assure that a patient in the $\ell$-th graph copy receives a kidney if his or her paired donor donates one. Constraints \eqref{eq:EMaxOutput} allow every donor (paired or singleton) to donate at most one kidney in only one copy of the graph. Constraints \eqref{eq:EChainSize} and \eqref{eq:ECycleSize} guarantee the maximum length allowed for chains and cycles in terms or arcs. Constraints \eqref{eq:ESymmetry} assure that a copy is selected, only if the vertex with the lowest index is selected. Lastly, constraints \eqref{eq:Eintegrality} define the nature of the decision variables.

The same counter example used for the EAF defined in Section \ref{sec:EAF} applies to the EEF. For the example given in Figure \eqref{fig:Copies2}, consider $\Kcycle = 3$ and $\Lchain = 6$. Note that only the chain copy is selected, but due to the presence of subtours of length superior to the cycle size limit, both formulations can provide an infeasible solution. Therefore, infeasible-cycle-breaking constraints are required for every chain copy.

\section{Proofs}
\label{sec:Proofs}

In this section,  {\cred we provide the proofs  on the validity of IEEF and the complexity of finding a positive-price column via MDDs.} For the sake of consistency, the numeration of the following propositions coincides with that used in the main body of the paper.

{\cred
\setcounter{theorem}{1}
\begin{prop}
	\label{prop:IEEF_Proof}
	IEEF is a correct formulation of the matching problem in the KEP.
\end{prop}

\proof{Proof.}

First, we need to show that a feasible matching to the KEP, $\matchingPartI$, corresponds to a feasible solution to IEEF whose objective function value equals the sum of the arc weights in  $\matchingPartI$, which is
\vspace{-5pt}
\begin{equation*}
	\sum_{\arcp \in \Arcs(\matchingPartI)} \weight_{\pairarc}.
	\vspace{5pt}
\end{equation*}	

\noindent  Formally, let $\matchingPartI = \cyfeasetD^{\prime} \cup \chfeasetD^{\prime}$ where $ \cyfeasetD^{\prime} \subseteq   \cyfeasetD$ and $\chfeasetD^{\prime} \subseteq \chfeasetD$. Since there are as many chain copies as \singles\ and a single \single\ per copy, for a chain $p = (\vix_{1},...,\vix_{\ell}) \in \chfeasetD^{\prime}$, there exists a unique copy $\fm \in \fbcopiesch$ for which the $\fm$-th \single, $\uix_{\fm}$, corresponds to  $\vix_{1}$. We then denote by $\fbcopiesch^{\prime} \subseteq \fbcopiesch$ the set of chain copies corresponding to altruistic donors in $\singleset \cap \Vertex(\matching) \ne \emptyset$, which can be empty, i.e., there are no chains in a feasible solution. Thus, we define $\lagy^{\prime}$ as
\vspace{20pt}

\begin{minipage}{0.5\textwidth}
	$$
	\lagy^{\fm \prime}_{\pairarc} = \left\{
	\begin{array}{ll}
		1      & \mbox{if } \arcp \in \Arcs(p)\\
		0      & \mbox{o.w.}
	\end{array}
	\right.
	$$
\end{minipage}
\begin{minipage}{0.5\textwidth}
	$ \forall  \fm \in \fbcopiesch^{\prime}, p = (\uix_{\fm}, \vix_{2}, ...,\vix_{\ell}) \in \chfeasetD^{\prime} $
\end{minipage}
and

\begin{minipage}{0.5\textwidth}
	\qquad\qquad\enspace $\lagy^{\fm \prime} = 0$ 
\end{minipage}
\begin{minipage}{0.5\textwidth}
	$\forall \fm \in \fbcopiesch \setminus \fbcopiesch^{\prime}.$
\end{minipage}

\vspace{10pt}
In the case of cycles, suppose $\fbset = \{\vix_{1}, ..., \vix_{\mid \fbset \mid}\}$ is a $\Kcycle$-limited FVS provided by Algorithm \ref{FVSAlg0}. Then, for a cycle $\cycle = (\vix_{1},...,\vix_{\Kix}, \vix_{1}) \in \cyfeasetD^{\prime}$ at least one of its vertices is also in $\fbset$ by the definition of a $K$-limited FVS (see Section \ref{sec:IEEF}). Let $e$ be the smallest position in $\fbset$ where a feedback vertex $\vixstar_{e} \in \Vertex(\cycle)$ is at, which coincides with the position of that copy in $\fbcopiescy$. Then, for the cycle $\cycle \in \cyfeasetD^{\prime}$ we define

\vspace{20pt}

\begin{minipage}{0.5\textwidth}
$$
\lagx^{e \prime}_{\pairarc} = \left\{
\begin{array}{ll}
	1      & \mbox{if } \arcp \in \Arcs(c)  \\
	0      & \mbox{o.w.}
\end{array}
\right.
$$
\end{minipage}
\begin{minipage}{0.5\textwidth}
		$ \forall  \arcp \in \Arcsd^{e}.$
\end{minipage}

\vspace{10pt}

We set the other components of $\lagx^{\prime}$ vector similarly. At $(\lagx^{\prime}, \lagy^{\prime})$, constraints \eqref{eq:maxflowIEEF} are satisfied since $\matchingPartI$ is feasible, i.e., every donor can donate at most once. This unitary flow satisfies flow-balance constraints in cycle copies \eqref{eq:OutFlowCycleIEEF} and chain copies \eqref{eq:OutFlowChainIEEF}. By the construction of $(\lagx^{\prime}, \lagy^{\prime})$, we know that at most one cycle/chain is selected in each copy, and they are the same cycles and chains of $\matchingPartI$. Thus, it is easy to see that the sum of arc variables for cycle \eqref{eq:CycleSizeIEEF} and chain \eqref{eq:ChainSizeIEEF} copies does not exceed $\Kcycle$ and $\Lchain$, respectively. Constraints \eqref{eq:symIEEF} and \eqref{eq:symChainsIEEF} are also satisfied since the left-hand side can take one as its maximum value, by \eqref{eq:maxflowIEEF}, and only arcs from the chain and cycle copies whose critical vertices $\uix_{\fm}$ and $\vixstar_{\fm}$ appear in $\matchingPartI$ are selected while constructing $(\lagx^{\prime}, \lagy^{\prime})$. Observe that constraints \eqref{eq:cycleElimIEEF} are applied to cycles present in chain copies. They are satisfied by  $\lagy^{\prime}$ since only the set of arcs forming a chain is selected in $\lagy^{\prime}$ for copies in $\fbcopiesch^{\prime}$, and no arc is selected from copies in $\fbcopiesch \setminus \fbcopiesch^{\prime}$.

Therefore, $(\lagx^{\prime}, \lagy^{\prime})$ is a feasible solution to IEEF. Lastly, it is easy to see that its objective value equals 
\vspace{-5pt}
\begin{equation*}
\sum_{\fm \in \fbcopiescy}\sum_{\arcp \in \Arcsd^{\fm}} w_{\pairarc}\lagx^{\fm \prime}_{\pairarc} + \sum_{\fm \in \fbcopiesch}\sum_{\arcp \in \Arcsbar^{\fm}} w_{\pairarc}\lagy^{\fm \prime}_{\pairarc} = \sum_{\arcp \in \Arcs(\matching)} \weight_{\pairarc},
\vspace{5pt}
\end{equation*}

Again, by construction, we need to show for the reverse direction that any feasible solution of IEEF corresponds to a feasible matching of the KEP, $\matching^{\star}$. Let $(\lagx^{\star}, \lagy^{\star})$ be a feasible solution of IEEF. Because IEEF can have symmetric solutions, a cycle $\cycle \in  \cyfeasetD$ can be selected in (a) the cycle copy corresponding to its associated feedback vertex as determined by Algorithm \ref{FVSAlg0},  (b) another cycle copy, or  (c) a chain copy.  However, constraints \eqref{eq:symIEEF} allow case (b) only when the feedback vertex of that copy is also selected. Similarly, constraints \eqref{eq:symChainsIEEF} allow case (c) only when the \single\ associated to a chain copy is also selected. A chain, on the other hand, can only be selected in a unique chain copy by construction of the chain copies. From \eqref{eq:maxflowIEEF}, we guarantee that the flow leaving every vertex is at most one across all copies, whether  they are cycle or chain copies. Thus, an arc traversal in IEEF is vertex-disjoint.  

 To build feasible cycles from $(\lagx^{\star}, \lagy^{\star})$ to include in $\matching^{\star}$, we start by inspecting the arc flow leaving feedback vertices in their corresponding copies, i.e., $\sum_{\vix: (\vix^{\star}_{\fm}, \vix) \in \Arcsd^{\fm}} \lagx_{\vix^{\star}_{\fm}\vix} = 1, \fm \in \fbcopiescy$. Constraints \eqref{eq:OutFlowCycleIEEF} enforce flow preservation, i.e., if a vertex  gives one unit of flow, it must also receive one. Thus, the arc traversal starting from $\vix_{\fm}^{\star}$ must be a cycle and it is also simple by \eqref{eq:maxflowIEEF}. Constraint \eqref{eq:CycleSizeIEEF} enforces the number of selected arcs in a cycle copy to be at most $\Kcycle$. Thus, under case (a), the arc traversal is guaranteed to be a feasible cycle $\cycle \in \cyfeasetD$. If $\mid \Arcs(\cycle) \mid  \le \Kcycle - 2$, then we evaluate case (b), since a cycle can still arise. We proceed to inspect the outgoing flow of other feedback vertices also present in that copy. It is easy to see that due to constraints \eqref{eq:maxflowIEEF}, \eqref{eq:OutFlowCycleIEEF} and \eqref{eq:CycleSizeIEEF} the new arc traversal  is also a simple cycle $\tilde{\cycle} \in \cyfeasetD$, such that $ \mid \Arcs(\tilde{\cycle}) \mid \le \Kcycle - \mid \Arcs(\cycle) \mid$. Thus, there can be as many cycles inside a cycle copy as long as the total number of arcs selected does not exceed $\Kcycle$. Under case (c), when an \single\ is selected in a chain copy, a subtour may also be selected. By constraints \eqref{eq:maxflowIEEF} and \eqref{eq:OutFlowChainIEEF} one can see that such subtours correspond to simple cycles. Nonetheless, their cardinality can exceed $\Kcycle$ when $\Lchain > \Kcycle$. Constraints \eqref{eq:cycleElimIEEF} remove all $\cycle  \in \csetf \setminus \csetf_{\Kcycle}, \ \fm \in \fbcopiesch$, thus, guaranteeing that if a subtour is selected it corresponds to a feasible cycle.

To include chains in $\matching^{\star}$, constraints \eqref{eq:maxflowIEEF}, \eqref{eq:OutFlowChainIEEF}  and \eqref{eq:ChainSizeIEEF} guarantee the arc traversal in a chain copy to contain a feasible chain. Note that all the cycles and chains included in $\matching^{\star}$ are vertex-disjoint due to \eqref{eq:maxflowIEEF}. Lastly, the weight of each arc variable in IEEF has a one-to-one correspondence with a transplant weight in $\matching^{\star}$. Thus, the objective value of $(\lagx^{\star}, \lagy^{\star})$ corresponds to $\sum_{\arcp \in \Arcs(\matching^{\star})} \weight_{\pairarc}$, which completes the proof.\hfill$\square$}

\setcounter{theorem}{3}
\begin{prop}
	Given the reduced costs $\hat{r}_{\cycle}^{\fm}$ and $\bar{r}^{\fm}_{\chain}$ expressed as an arc-separable function for all $(\dgv_{s}, \dgv_{s'})  \in \marcs$ and $(\dgv_{s}, \dgv_{s'})  \in \marcsch$, a positive-price cycle, if one exists, can be found in  time $\mathcal{O}(\sum_{\fm \in \fbcopiescy} \sum_{(\dgv_{s}, \dgv_{s'})  \in \marcs}  \vert \inn(\dgv_{s}) \vert)$. Similarly, a positive-price chain can be found in $\mathcal{O}(\sum_{\fm \in \fbcopiesch} \sum_{(\dgv_{s}, \dgv_{s'})  \in \marcsch}  \vert \inn(\dgv_{s}) \vert)$.
\end{prop}

\proof{Proof.}
For every arc $(\dgv_{s}, \dgv_{s'}) \in \marcs$ and $(\dgv_{s}, \dgv_{s'}) \in \marcsch$, $\vert \inn((\dgv_{s}, \dgv_{s'})) \vert$ comparisons need to be performed to obtain $\recucy((\dgv_{s}, \dgv_{s'}))$ and $\recuch((\dgv_{s}, \dgv_{s'}))$, respectively in \eqref{eq:recu1a}. Thus, for the $\fm$th MDD of a cycle copy, $\sum_{(\dgv_{s}, \dgv_{s'}) \in \marcs} \vert \inn(\dgv_{s}) \vert$ comparisons are required to compute $\recucy$, whereas for the $\fm$th MDD of a chain copy,  there are the same number of comparisons plus $\vert \marcsch \vert$ comparisons  of all arcs, in \eqref{eq:2b}, before obtaining $\recuch$. Since there are $\vert \fbcopiescy \vert$ cycle MDDs and $\vert \fbcopiesch \vert$ chain MDDs, it follows that the time complexity is as shown above. \hfill$\square$

\begin{prop}
\label{pre:sizeMDDcycles}
	The size of the input $\sum_{\fm \in \fbcopiescy} \sum_{(\dgv_{s}, \dgv_{s'}) \in \marcs} \vert \inn(\dgv_{s}) \vert$ grows as $\vert \Vertexd^{\fm} \vert ^{\Kcycle + 1}$ does.
\end{prop}

\proof{Proof.}
Without lost of generality, assume $\digraph$ is complete. As stated before, the layer of an arc $\ma: (\dgv_{s}, \dgv_{s'}) \in \marcs$ is the layer to which its source node belongs, i.e., if node $\dgv_{s'}$ is on layer $\Kix$, then $\ell(\ma)  = \Kix$. Moreover, let $\marcs_{\Kix} := \{\ma \in \marcs \mid \ell(\ma) = \Kix \}$ be the set of arcs that belong to layer $\Kix$. By construction of the MDDs, $\mathbf{r}$ has only one outgoing arc such that $val((\mathbf{r}, \dgv_1)) = \vixstar_{\fm}$. In the second layer, $\layer_{2}$, the cardinality of $\marcs_{2}$ can be up to $\vert \Vertexd^{\fm} \vert - 1$ corresponding to the vertices $\vix \in  \Vertexd^{\fm} \setminus \{\vixstar_{\fm}\}$ that can be selected in the second position of a cycle. Note that since there can be the same number of 2-way cycles,  an arc $\ma \in \marcs_{2}$ has a sink node $\dgv$ such that $(\dgv, \mathbf{t}) \in \marcs$. If the length of a cycle is larger than two, then the cardinality of $\marcs_{3}$ can be up to $\vert \Vertexd^{\fm} \vert - 2$, representing the $\vert \Vertexd^{\fm} \vert - 2$ vertices $\vix \in  \Vertexd^{\fm} \setminus \{\vixstar_{\fm}\}$ that can be selected in the third position of a cycle. The same process is repeated until in layer $\layer_\Kcycle$ there are $\vert \Vertexd^{\fm} \vert - (\Kcycle - 1)$ vertices $\vix \in \Vertexd^{\fm}$ to choose, and thus, $\vert \Vertexd^{\fm} \vert - (\Kcycle - 1)$ arcs  $\ma \in \marcs$ pointing to $\dgv$. Thus, $\vert \marcs \vert$ equals
%\vspace{-3mm}
\begin{subequations}
	\begin{align}
	\prod_{\Kix =2}^{\Kcycle} \vert \Vertexd^{\fm} \vert  - (\Kix - 1) + \sum_{\Kix = 2}^{\Kcycle - 1} (\vert \Vertexd^{\fm} \vert  - (\Kix - 1) ) <  \vert \Vertexd^{\fm}  \vert ^{\Kcycle - 1} + \Kcycle \vert \Vertexd^{\fm} \vert
	\end{align}
\end{subequations}

The second sum on the left-hand side corresponds to the number of arcs at every layer whose sink node is $n$, thus, closing up cycles using fewer than $\Kcycle$ arcs. Therefore, under a worst-case scenario in which $ \vert \inn(\dgv_{s}) \vert$ and $\vert \fbcopiescy \vert$ tend to $\vert \Vertexd^{\fm} \vert$, the complexity of finding a positive-price cycle becomes $\mathcal{O}(\vert \Vertexd^{\fm} \vert ^{\Kcycle + 1})$. \hfill$\square$

\begin{prop}
	The size of the input $\sum_{\fm \in \fbcopiesch} \sum_{(\dgv_{s}, \dgv_{s'}) \in \marcsch} \vert \inn(\dgv_{s}) \vert$ grows as $\vert \Vertexbar^{\fm} \vert ^{\Lchain + 2}$ does for bounded chains and as $\vert  \Vertexbar^{\fm}  \vert !$ when $\Lchain \rightarrow \infty$.
\end{prop}	

\proof{Proof.}
A similar reasoning to proposition \ref{pre:sizeMDDcycles} can be followed, except that a chain can be cut short if by visiting a new vertex $\vix \in \Vertexbar^{\fm}$ in a sequence of the state transition graph at least one \pairAb\ is present more than once, thereby violating the condition of being a simple path. We know that for a path to have $\Lchain$-many arcs, it is necessary to have a sequence with $\Lchain$ \pairsAb, thus, $\vert \marcs \vert$ tends to $\vert \Vertexbar^{\fm} \vert^{\Lchain}$ and $\sum_{\fm \in \fbcopiesch}  \sum_{(\dgv_{s}, \dgv_{s'}) \in \marcs}  \vert \inn(\dgv_{s}) \vert \approx \vert \Vertexd^{\fm} \vert^{\Lchain + 2}$. Therefore, for bounded chains, finding a positive-price chain can be done in time  $\mathcal{O}(\vert \Vertexd^{\fm} \vert ^{\Lchain+ 2})$. The second part follows by the fact that after we visit $\Lix$ vertices, there are still $\vert \Vertexd^{\fm} \vert - \Lix$ ways to choose the next one, until only one can be chosen, thus, the time to find a positive-price column when $\Lchain$ is unbounded is exponential. \hfill$\square$

\section{Additional Results}
\label{sec:OtherResults}

\paragraph{Columns by phase \newline}

Figure \ref{fig:cols.pdf} depicts the type and number of columns found in each phase for individual runs across all the $\Kcycle$-$\Lchain$ combinations. The {\cdionne $x$-axis} represents the total time  (in minutes) to solve the pricing problem of a single run, i.e., an instance on a $\Kcycle$-$\Lchain$ setting,  {\cdionne during the three phases. For every $x$-value} there may be multiple y-points representing the number of columns found in a specific phase (sub phase) in thousands and whether they are cycle or chain columns. Therefore, for the same {\cdionne $x$-value}, multiple {\cdionne $y$-values} may correspond to the same instance, maybe on different settings, if the markers share the same size and color, regardless of the shape. For similar total pricing-problem times, markers may overlap. As an example, consider the time interval from 2min to 10min. In the cycle subplot there are big (purple) circle markers indicating that from 10k up to 30k cycle columns were found in Phase 1, as apposed to the chain subplot where no such markers appear. This means that this subset of instances correspond to $\Lchain = 0$. In another example, some instances with 2252 \pairsAb\ and 204 \singles\ whose pricing time is  15.6min and have a green circle around 12k and 6k in the cycle and chain subplots, respectively, plus a triangle indicating one chain found through \eqref{objLP}. Thus, the total number of columns of a run is the sum over the number of columns indicated by all markers in the y-axis with respect to the same {\cdionne $x$-value}, provided that markers share color and size. {\cdionne In this particular case of 2048 \pairsAb\ and 204 \singles\ there are 17,396 columns.} Overall, most markers are circles, indicating that the majority of columns are found via MDDs across all runs, while {\cdionne Phase 2 and Phase 3} mostly  certify that no more positive-price columns exist. %The minimum number of columns (greater than zero) found jointly in PH2-3 is 1 and the maximum is 447.the time interval from 2min to 10min. In the cycle subplot there are big (purple circle) markers indicating that from 10k up to 30k cycle columns were found in PH1, as opposed to the chain subplot where no such markers appear. This means that this subset of instances corresponds to $\Lchain = 0$. In another example,
\begin{figure}[h!]
	\centering
	% include first image
	\includegraphics[scale=0.42]{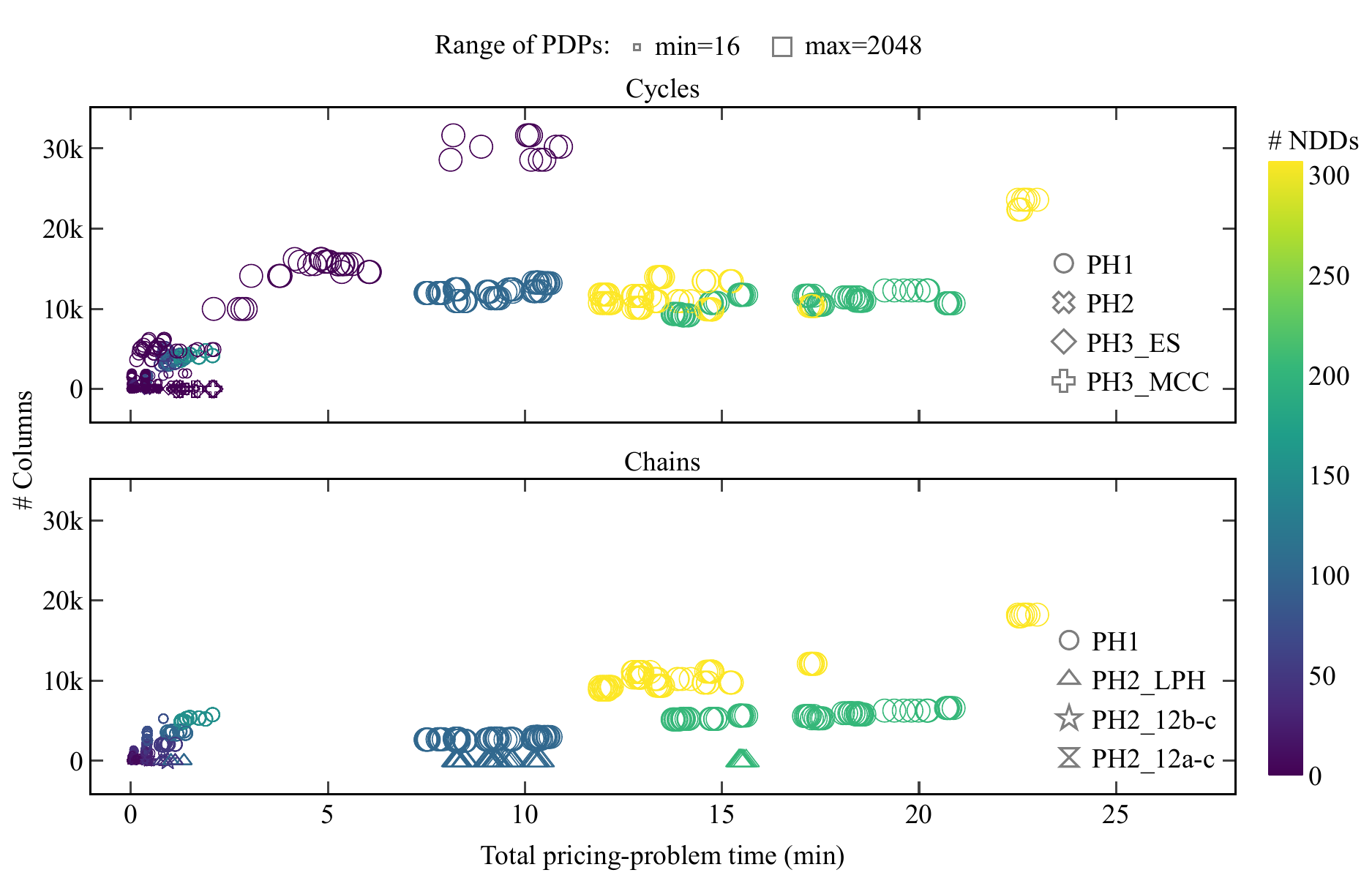}  %width=\linewidth
	\caption{\cdionne {Classification of columns by run, phase (PH1, PH2, PH3) and type. Marker types indicate different (sub) phases of column generation. Marker sizes are correlated with the number of \pairsAb, while marker colors indicate the number of \singles. ES stands for the Exhaustive Search in PH2.}}
	\label{fig:cols.pdf}
\end{figure}

{
	\cred
\paragraph{Effect of vertex-selection rules on BP\_MDD \newline}

For instances in the set KEP\_N0, Table \ref{tab:FVSAnalysis} provides average values on the performance of BP\_MDD in terms of the number of states (i.e., the number of nodes in MDDs) and the solution time when different selection rules are used in Algorithm \ref{FVSAlg0}. Selection rule Incr.Index corresponds to selecting every pair in the graph in increasing order of their index, as in \citep{Constantino2013}. The column $K$-FVS corresponds to finding the $K$-FVS with smallest cardinality  through a MIP formulation.  This formulation was given a time limit of 1hr. There is a hyphen if a feasible solution was not available after this time.

\begin{table}[h!]
\spacing
	\centering
	\caption{\cred Trade-off between performance and scalability when $\Kcycle = 3$}
	\label{tab:FVSAnalysis}
	\cred
    \begin{adjustbox}{width=0.8\textwidth}
	\begin{tabular}{lrrrrrrrr}
		\toprule
		 & \multicolumn{4}{c}{Average \# States}                                                                                           & \multicolumn{4}{c}{Solution Time (s)}                                                                      \\ \cmidrule(lr){2-5} \cmidrule(lr){6-9} 
		\multicolumn{1}{c}{Set}                     & \multicolumn{1}{l}{$K$-FVS} & \multicolumn{1}{l}{Incr.Index} & \multicolumn{1}{l}{Indegree} & \multicolumn{1}{l}{Total degree} & \multicolumn{1}{l}{$K$-FVS} & \multicolumn{1}{l}{Incr.Index} & \multicolumn{1}{l}{Indegree} & Total degree \\
\cmidrule(lr){1-1}
\cmidrule(lr){2-2}
\cmidrule(lr){3-3}
\cmidrule(lr){4-4}
\cmidrule(lr){5-5}
\cmidrule(lr){6-6}
\cmidrule(lr){7-7}
\cmidrule(lr){8-8}
\cmidrule(lr){9-9}
		N16\_A0                                   & 10.4                      & 20.3                            & 12.0                       & 11.4                              & 0.00                          & 0.00                               & 0.00                             & 0.00            \\
		N32\_A0                                   & 74.9                      & 124.1                        & 73.2                          & 73.4                              & 0.00                          & 0.00                               & 0.00                             & 0.00            \\
		N64\_A0                                   & 283.7                    & 479.8                           & 229.5                         & 282.4                           & 0.01                       & 0.01                            & 0.01                          & 0.01         \\
		N128\_A0                                  & 1,281.5                   & 1,753.6                         & 933.9                        & 1,192.2                           & 0.03                       & 0.02                            & 0.02                          & 0.02         \\
		N256\_A0                                  & 6,270.3                    & 7,542.1                         & 4,038.2                    & 5,646.7                           & 0.22                       & 0.14                            & 0.15                          & 0.15         \\
		N512\_A0                                  & 29,805.3                   & 31,531.5                     & 16,158.9                      & 25,537.5                          & 1.62                       & 0.96                            & 1.44                          & 1.34         \\
		N1024\_A0                                 & -                          & 129,608.1                       & 65,816.1                      & 112,277.8                         & -                          & 8.35                            & 22.84                         & 42.40        \\
		N2048\_A0                                 & -                          & 513,423.8                      & 252,659.7                     & 467,831.7                         & -                          & 85.29                           & 521.25                        & 484.82       \\ \bottomrule
	\end{tabular}
	\end{adjustbox}
\end{table}

Even though, Incr.Index  reports the smallest solution time on average, the number of states throughout all MDDs double those of the Indegree selection rule.  This selection rule offers the best trade-off between scalability and computational performance.

\paragraph{Prioritization for instances in KEP\_N \newline}
Table \ref{tab:SolStructureG1} shows the solution composition under the two prioritization schemes for a sample of the instance sets in KEP\_N. Header names are defined as in Table \ref{tab:SolStructure} of the Results Section. It was worth noting that, except for 23 runs in the set N588\_A76, all altruistic donors were used to trigger chains. For those runs, the percentage of used altruistic donors was 98.7\%.

\begin{table}[h!]
\spacing
	\centering
	\caption{\cred Solution composition by execution modes}
	\label{tab:SolStructureG1}
	\cred
	\begin{adjustbox}{width=0.8\textwidth}
	\begin{tabular}{ccccccccc}
		\toprule
		&  &  & \multicolumn{3}{c}{Averages on CH mode}                                                        & \multicolumn{3}{c}{Average Ratio CH/CY}                                    \\ \cmidrule(lr){4-6} \cmidrule(lr){7-9} 
		\multicolumn{1}{c}{Set}                     & \multicolumn{1}{c}{$\Kcycle$}                           & \multicolumn{1}{c}{$\Lchain$}                           & \multicolumn{1}{l}{AvChainLen} & \multicolumn{1}{l}{nCyclesSol} & \multicolumn{1}{l}{AvCycleLen} & \multicolumn{1}{l}{AvChainLen} & \multicolumn{1}{l}{nCycleSol} & AvCycleLen \\
\cmidrule(lr){1-1}
\cmidrule(lr){2-2}
\cmidrule(lr){3-3}
\cmidrule(lr){4-4}
\cmidrule(lr){5-5}
\cmidrule(lr){6-6}
\cmidrule(lr){7-7}
\cmidrule(lr){8-8}
\cmidrule(lr){9-9}
		N147\_A19                                 & 3                                               & 3                                               & 2.97                           & 17.22                           & 2.81                           & 1.50                           & 0.71                           & 1.02      \\
		N147\_A19                                 & 3                                               & 4                                               & 2.99                           & 17.00                           & 2.83                           & 1.51                           & 0.71                           & 1.02      \\
		N147\_A19                                 & 3                                               & 5                                               & 3.00                             & 17.11                           & 2.79                           & 1.55                           & 0.71                           & 1.00      \\
		N147\_A19                                 & 3                                               & 6                                               & 3.06                           & 16.50                           & 2.83                           & 1.58                           & 0.68                           & 1.02      \\
		N147\_A19                                 & 4                                               & 3                                               & 3.00                              & 17.21                           & 2.81                           & 1.49                           & 0.71                           & 1.03      \\
		N147\_A19                                 & 4                                               & 4                                               & 3.00                              & 17.00                           & 2.83                           & 1.52                           & 0.71                           & 1.02      \\
		N147\_A19                                 & 4                                               & 5                                               & 3.00                              & 17.12                           & 2.79                           & 1.55                           & 0.71                           & 1.00      \\
		N147\_A19                                 & 4                                               & 6                                               & 3.00                              & 16.55                           & 2.83                           & 1.58                           & 0.68                           & 1.02      \\
		N294\_A38                                 & 3                                               & 3                                               & 3.00                              & 30.63                           & 2.84                           & 1.76                           & 0.63                           & 1.02      \\
		N294\_A38                                 & 3                                               & 4                                               & 3.00                              & 30.43                           & 2.85                           & 1.67                           & 0.64                           & 1.02      \\
		N294\_A38                                 & 3                                               & 5                                               & 3.00                              & 30.57                           & 2.84                           & 1.58                           & 0.67                           & 1.01      \\
		N294\_A38                                 & 3                                               & 6                                               & 3.00                              & 30.70                           & 2.83                           & 1.61                           & 0.66                           & 1.01      \\
		N294\_A38                                 & 4                                               & 3                                               & 3.00                              & 30.67                           & 2.84                           & 1.74                           & 0.63                           & 1.02      \\
		N294\_A38                                 & 4                                               & 4                                               & 3.00                              & 30.44                           & 2.85                           & 1.68                           & 0.64                           & 1.03      \\
		N294\_A38                                 & 4                                               & 5                                               & 3.00                              & 30.52                           & 2.84                           & 1.58                           & 0.67                           & 1.01      \\
		N294\_A38                                 & 4                                               & 6                                               & 3.00                              & 30.71                           & 2.83                           & 1.60                           & 0.66                           & 1.02      \\
		N588\_A76                                 & 3                                               & 3                                               & 3.00                              & 61.30                           & 2.93                           & 1.54                           & 0.69                           & 1.03      \\
		N588\_A76                                 & 3                                               & 4                                               & 3.00                              & 61.30                           & 2.93                           & 1.56                           & 0.68                           & 1.03      \\
		N588\_A76                                 & 3                                               & 5                                               & 3.00                              & 61.30                           & 2.93                           & 1.36                           & 0.75                           & 1.04      \\
		N588\_A76                                 & 3                                               & 6                                               & 3.00                              & 61.30                           & 2.93                           & 1.41                           & 0.73                           & 1.03      \\
		N588\_A76                                 & 4                                               & 3                                               & 3.00                              & 61.30                           & 2.93                           & 1.53                           & 0.69                           & 1.03      \\
		N588\_A76                                 & 4                                               & 4                                               & 3.00                              & 61.30                           & 2.93                           & 1.50                           & 0.70                           & 1.03      \\
		N588\_A76                                 & 4                                               & 5                                               & 3.00                              & 61.30                           & 2.93                           & 1.36                           & 0.75                           & 1.04      \\
		N588\_A76                                 & 4                                               & 6                                               & 3.00                              & 61.30                           & 2.93                           & 1.41                           & 0.73                           & 1.03      \\ \bottomrule
	\end{tabular}
	\end{adjustbox}
\end{table}
}

%\section{Column Generation Algorithm}
%\label{CGDiagram}
%
%Figure \ref{fig:LogicDiagram} depicts the column generation algorithm.

% new section for proofs
\end{APPENDIX}
% CASE 2: BiBTeX used to generate mypaper.bbl (to be further fine tuned)
%\input{mypaper.bbl} % outcomment this line in Case 2

\end{document}